\documentclass[reqno]{amsart}
\usepackage{mathrsfs}
\usepackage{amsmath,amsfonts,amssymb,amsthm}
\usepackage{upgreek}
\usepackage[usenames,dvipsnames]{xcolor}
\usepackage{enumitem}
\usepackage{graphicx}
\usepackage{overpic}
\usepackage{cancel}
\usepackage[outdir=./]{epstopdf}
\usepackage[colorlinks=true]{hyperref} 
\usepackage{url}
\usepackage[utf8]{inputenc}
\usepackage{booktabs,multirow}
\usepackage[font=footnotesize,width=\textwidth]{caption}
\hypersetup{linkcolor=Violet,citecolor=PineGreen}
\newtheorem{teor}{Theorem}[section]

\newtheorem{prop}[teor]{Proposition}

\theoremstyle{definition}

\newtheorem{nota}[teor]{Remark}

\numberwithin{equation}{section}

\newcommand{\R}{\mathbb R}

\newcommand{\N}{\mathbb{N}}

\newcommand{\mA}{\mathcal{A}}
\newcommand{\mB}{\mathcal{B}}
\newcommand{\mC}{\mathcal{C}}
\newcommand{\mD}{\mathcal{D}}

\newcommand{\mK}{\mathcal{K}}
\newcommand{\mM}{\mathcal{M}}

\newcommand{\mP}{\mathcal{P}}
\newcommand{\mR}{\mathcal{R}}

\newcommand{\mU}{\mathcal{U}}
\newcommand{\mV}{\mathcal{V}}
\newcommand{\ep}{\varepsilon}

\newcommand{\mI}{\mathcal{I}}
\newcommand{\mJ}{\mathcal{J}}

\newcommand{\mQ}{\mathcal{Q}}

\newcommand{\W}{\Omega}
\newcommand{\WR}{\W\times\R}
\newcommand{\RR}{\R\times\R}
\newcommand{\w}{\omega}
\newcommand{\lb}{\lambda}

\newcommand{\ma}{\mathfrak{a}}
\newcommand{\mb}{\mathfrak{b}}
\newcommand{\mc}{\mathfrak{c}}
\newcommand{\md}{\mathfrak{d}}
\newcommand{\mr}{\mathfrak{r}}
\newcommand{\ml}{\mathfrak{l}}
\newcommand{\mm}{\mathfrak{m}}
\newcommand{\muk}{\mathfrak{u}}
\newcommand{\mf}{\mathfrak{f}}
\newcommand{\mg}{\mathfrak{g}}
\newcommand{\mh}{\mathfrak{h}}
\newcommand{\mpp}{\mathfrak{p}}
\newcommand{\tma}{\tilde\ma}
\newcommand{\tmr}{\tilde\mr}
\newcommand{\tml}{\tilde\ml}
\newcommand{\tmm}{\tilde\mm}
\newcommand{\tmuk}{\tilde\muk}

\newcommand{\upalfa}{$\upalpha$}
\newcommand{\upomeg}{$\upomega$}

\newcommand{\G}{\Gamma}

\newcommand{\wit}{\widetilde}

\newcommand{\ws}{\w{\cdot}s}
\newcommand{\wt}{\w{\cdot}t}
\newcommand{\pt}{p{\cdot}t}
\newcommand{\bwt}{\bar\w{\cdot}t}
\newcommand{\bws}{\bar\w{\cdot}s}

\newcommand{\n}[1]{\left\|#1\right\|}

\newcommand{\lsm}{\left[\begin{smallmatrix}}
\newcommand{\rsm}{\end{smallmatrix}\right]}

\newcommand{\merg}{\mathfrak{M}_\mathrm{erg}(\W,\sigma)}

\begin{document}
\title[Saddle-node bifurcations for nonautonomous ODEs]
{Saddle-node bifurcations for concave in measure and d-concave in measure
skewproduct flows with applications to population dynamics and circuits}
\author[J. Due\~{n}as]{Jes\'{u}s Due\~{n}as}
\author[C. N\'{u}\~{n}ez]{Carmen N\'{u}\~{n}ez}
\author[R. Obaya]{Rafael Obaya}
\address{{\rm(J. Due\~{n}as, C.N\'{u}\~{n}ez, R. Obaya)}. Departamento de Matem\'{a}tica Aplicada, Universidad de Va\-lladolid, Paseo Prado de la Magdalena 3-5, 47011 Valladolid, Spain. The authors are members of IMUVA: Instituto de Investigaci\'{o}n en Matem\'{a}ticas, Universidad de Valladolid.}
\email[J.~Due\~{n}as]{jesus.duenas@uva.es}
\email[C.~N\'{u}\~{n}ez]{carmen.nunez@uva.es}
\email[R.~Obaya]{rafael.obaya@uva.es}
\thanks{All the authors were supported by Ministerio de Ciencia, Innovaci\'{o}n y Universidades (Spain)
under project PID2021-125446NB-I00 and by Universidad de Valladolid under project PIP-TCESC-2020.
J.~Due\~{n}as was also supported by Ministerio de Universidades (Spain) under programme FPU20/01627.}
\date{}
\begin{abstract}
Concave in measure and d-concave in measure nonautonomous scalar ordinary differential equations given by coercive and
time-compactible maps have similar properties to equations satisfying considerably more restrictive hypotheses.
This paper describes the generalized simple or double saddle-node bifurcation diagrams for one-parametric
families of equations of these types, from which the dynamical possibilities for each of the equations
follow. This new framework allows the analysis of ``almost stochastic" equations, whose
coefficients vary in very large chaotic sets. The results also apply to the analysis of the
occurrence of critical transitions for a range of models much larger than in previous approaches.
\end{abstract}
\renewcommand{\subjclassname}{\textup{2020} Mathematics Subject Classification}

\maketitle
\section{Introduction}
This paper deals with the dynamical properties of certain types of scalar skewproduct
flows generated by nonautonomous ordinary differential equations. Classically,
these flows are called concave when they are generated by a vector field that is concave
with respect to the state variable. Analogously, they are called d-concave when they
are generated by a vector field with a concave derivative with respect to the state component.
These properties are relaxed in the paper Due\~{n}as {\em et al.}~\cite{dno4}, where the
applicability of the new hypotheses to analyze the occurrence of critical transitions
for mathematical models not fitting in the previously analyzed concave or d-concave settings is shown.
In this new work, we provide an exhaustive analysis of the dynamical properties of these scalar flows,
which we will refer to as {\em concave in measure\/} and {\em d-concave in measure}:
we describe the possibilities for their global dynamics,
obtain the global bifurcation diagrams of some suitable parametric families,
pose some models of population dynamics and electrical circuits to show the interest of the
obtained conclusions, and show the applicability of these new results in the analysis of
the occurrence of critical transitions in a considerably wider range of cases than in previous approaches.

Let $\sigma\colon\R\times\W\to\W$ define a
(general) continuous global flow on a compact metric space $\W$, let us denote $\wt:=\sigma(t,\w)$,
and let $\mh\colon\WR\to\R$ be continuous. The solutions of the family of equations
$\{x’=\mh(\wt,x)\,,\;\w\in\W\}$ induce a (possibly local) skewproduct flow $\tau\colon\mV\subseteq\R\times\WR\to\WR$:
$\tau(t,\w,x)=(\wt,v(t,\w,x))$, where $v(t,\w,x)$ satisfies
$x’=\mh(\wt,x)$ with value $x$ at $t=0$. We assume enough regularity on $\mh$, represent by
$\mh_x$ and $\mh_{xx}$ its first and second derivatives with respect to $x$, and say that
the flow $\tau$ is {\em concave in measure\/} (resp. {\em d-concave in measure}) if
$\{\w\in\W\,|\;x\mapsto\mh(\w,x)\text{ is concave (resp. $x\mapsto \mh(w,x)$ is d-concave)}\}$ has complete measure (i.e., measure 1 for
each ergodic measure on $\W$) and, in addition, for each compact interval $\mJ\subset\R$, the set
$\{\w\in\W\,|\;x\mapsto\mh_x(\w,x)\text{ (resp. $x\mapsto\mh_{xx}(\w,x)$) is}$  $\text{strictly decreasing on }\mJ\}$
has positive measure for every ergodic measure.

Skewproduct flows have been used to analyze the dynamics of a single nonautonomous equation $x'=h(t,x)$ during the
last 70 years, following the pioneer results of Bebutov \cite{bebu}. In these cases,
$\W=\W_h$ is given by the uniform hull of $h$, i.e., the closure of the set of time shifts $h_s(t,x):=
h(t+s,x)$ in the compact-open topology of $C(\RR,\R)$, which is compact under suitable conditions on $h$;
and $\mh(\w,x):=\w(0,x)$. As explained in \cite[Lemma 2.4 and proof of Lemma 4.3]{dno4}, the set $\W_h$ is composed
by three (perhaps equal) sets: $\{h_s\,|\;s\in\R\}$, and the corresponding \upalfa-limit set
$\W_h^\alpha$ and \upomeg-limit set $\W_h^\omega$; and the ergodic measures on $\W_h$ are concentrated either in
$\W_h^\alpha$ or in $\W_h^\omega$. So, in this hull setting, it is enough to restrict ourselves to
$\W_h^\alpha\cup\W_h^\omega$ to get the required concavity or d-concavity hypotheses.

The main results in \cite{dno4} focus on critical transitions, in the line of \cite{alas}, \cite{altw},
\cite{aspw}, \cite{awvc}, \cite{hcw1}, \cite{scheffer2}, \cite{scheffer1} and \cite{sebastian1}.
For the results on critical transitions in \cite{dno4}, and
for the examples of application there provided, the required concavity or d-concavity
in measure is obtained for an equation $x'=h(t,x)$ given by a regular enough map $h(t,x)$ that
can be suitably approached by two maps $h_\pm(t,x)$ as $t\to\pm\infty$, assuming global concavity
or d-concavity of these maps as well as the existence of three hyperbolic solutions for each equation
$x'=h_\pm(t,\w)$. So, in that paper we assume $h$ to be asymptotically concave or d-concave,
understand $x'=h(t,x)$ as a transition between its \upalfa-limit set and its \upomeg-limit set, and assume
also a strong hyperbolic structure for the dynamics over these sets.
The relations $\W_h^\alpha=\W_{h_-}^\alpha$ and $\W_h^\omega=\W_{h_+}^\omega$ combined with those global
properties are the key points to show the concavity or d-concavity in measure. But, even in the hull
framework, the range of applicability is much wider than this, as the next two examples show.

\begin{figure}[h]
\centering
\includegraphics[width=\textwidth]{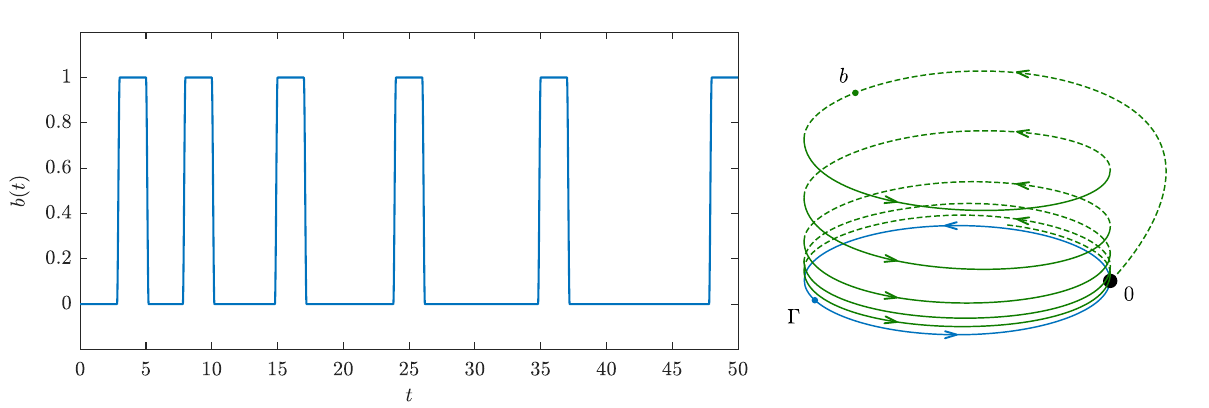}
\caption{In the left panel, the function $b(t):=\sum_{n=2}^\infty\Gamma(t-n^2)$, where
$\G$ is defined as the unique $C^1$ cubic spline which takes value $1$ on $[-1,1]$ and $0$
outside $[-1.2,1.2]$.
In the right panel, a sketch of the structure of the hull $\W_b$
of $b$. The unique critical point $0$ is denoted by a black large point.
In green, a sketch of the orbit $\{b_s\,|\;s\in\R\}$ of $b$ and, in blue, of that of $\G$.
The hull $\W_\G$ of $\G$ is composed by its orbit and $\{0\}$ as its \upalfa-limit set and \upomeg-limit set,
and it is the \upomeg-limit set $\W_b^\omega$ of $b$; and the \upalfa-limit set is $\W_b^\alpha=\{0\}$.}
\label{fig:tophull}
\end{figure}

We consider a nonautonomous scalar equation $x'=h(t,x)$ representing the
dynamics of a single-species population subject to increasingly spaced hunting periods,
\begin{equation}\label{eq:huntingintro}
 x'=\frac{1}{6}\,x\,(1-x)-b(t)\,\frac{x^2}{1+x^2}\,,
\end{equation}
where $b\colon\R\to[0,1]$ is the bounded and uniformly continuous map depicted
in the left panel of Figure \ref{fig:tophull}.
It can be checked that $x\mapsto h(t,x)$ is strictly concave if $b(t)<2/3$ and not concave if $b(t)>2/3$.
In this case, $\W_h$ can be identified with $\W_b$, where $\W_b$ is the (compact)
closure of the set of time shifts $b_s(t):=b(t+s)$
in the compact-open topology of $C(\R,\R)$. Equation \eqref{eq:huntingintro} is embedded in the family
$x'=x\,(1-x)/6-\w(t)\,x^2/(1+x^2)$ for $\w\in\W_b$; $\W_b$ is the union of $\{b_s\,|\;s\in\R\}$,
$\W_b^\alpha=\{0\}$ and $\W_b^\omega\supsetneq\{0\}$; and the unique ergodic measure on $\W_b$ is the
Dirac measure concentrated on $\{0\}$ (see the right panel and the caption of Figure \ref{fig:tophull}).
So, \eqref{eq:huntingintro} is concave in measure (i.e., the corresponding skewproduct is),
although $h$ is not asymptotically concave. The results we present in this paper will describe all the
possibilities for its global dynamics, as well the global bifurcation diagram of
$x'=h(t,x)+\lambda$. Similarly, the equation
\begin{equation}\label{eq:huntingintrod}
 x'=\frac{1}{6}\,x\,(1-x)(x-2)-b(t)\,\frac{x^2}{1+x^2}\,,
\end{equation}
which incorporates the multiplicative Allee effect on the population, is d-concave in measure
but not asymptotically d-concave.

The situation of these examples is quite general in the hull setting: it is
rather common that all the ergodic measures for these \upalfa-limit and \upomeg-limit sets
are concentrated in topologically small sets (with void interior or even of the first Baire category).
So, the results that we present are valid under conditions on concavity or d-concavity imposed for the elements of
very small subsets of the hull. That is, even in the hull framework, the general
definition and the exhaustive description of the dynamical properties of concave in measure and d-concave
in measure scalar flows that we provide in this work are justified by the much wider range of
applicability with respect to previous approaches, as those of \cite{nuos4} and \cite{dno1}.

In fact, in
those papers, the flow on $\W_h$ is assumed to be minimal: every $\sigma$-orbit is dense.
And this property results in, roughly speaking, ``the same dynamics" for each one of the equations of the family.
For instance, under this hypothesis, in the concave and coercive case, all the equations have
simultaneously either two hyperbolic solutions, bounded but not hyperbolic solutions, or no bounded solutions at all.
The dynamical structure for all the equations of the family is also common in the d-concave
and coercive minimal case. Here, the maximum number of hyperbolic solutions is three, and many other possibilities arise:
the existence of only one bounded (hyperbolic) solution, the existence of a unique hyperbolic solution uniformly
separated from a bounded nonhyperbolic one, the existence of two nonseparated hyperbolic solutions,
and the absence of hyperbolic solutions with two bounded separated solutions.
In the analysis performed in this new work, this minimality condition is removed.
In fact, the flow on $\W$ is not even assumed to be transitive, i.e., to have a dense orbit,
as it happens in the general hull setting.
The situation is much more complex in this much more general setting, despite of which a global description
is possible.

This paper presents a contribution to nonautonomous bifurcation theory, focusing on
saddle-nodes: see, e.g., \cite{anapotmar}, \cite{kloeden1}, \cite{nuob6}.
In fact, the description of the dynamical possibilities for each equation $x'=\mh(\wt,x)$,
made under an additional condition on coercivity for $\mh$, is intrinsically related
to the obtention of the global bifurcation diagram for parametric families that include $x'=\mh(\wt,x)+\lb$.
In the concave case, analyzed under more restrictive hypotheses in
\cite{alob3}, \cite{anja}, \cite{lnor}, \cite{lno3} and \cite{remo1},
we prove the existence of a unique bifurcation value $\lb_\w$ for each $\w\in\W$
such that the equation has no bounded solutions for $\lb<\lb_\w$, bounded but no hyperbolic solutions for
$\lb=\lb_\w$, and two hyperbolic solutions composing an attractor-repeller pair for $\lb>\lb_\w$ whose distance
decreases to $0$ as $\lb\downarrow\lb_\w$. The dynamics of $x'=\mh(\wt,x)$ depends, hence, on the sign of $\lb_\w$.
Contrary to the minimal case, the value of $\lb_\w$ may vary with $\w$, even in the transitive case.
However, there exists an interval $[\lb_*,\lb_-]$ containing all these values $\lb_\w$,
with $\lb_*=\inf_{\w\in\W}\lb_\w$ and $\lb_-=\sup_{\w\in\W}\lb_\w$. So, the family of equations
$x'=\mh(\wt,x)$ has an attractor-repeller pair of copies of the base
if $0>\lb_-$, no bounded orbits if $0<\lb_*$, and a $\w$-dependent dynamics if $0\in[\lb_*,\lb_-]$.
That is, we have a global saddle-node nonautonomous bifurcation diagram
for each equation, while the bifurcation diagram for the family can be understood as
a generalized global saddle-node one. Figure \ref{fig:diagramasconcavos} and its caption
explain these diagrams in detail.

In the d-concave in measure case, analyzed in \cite{dno1} in the minimal case,
we focus on the case of existence of $\lb_0$ such that the flow given by the family
$\{x'=\mh(\wt,x)+\lb_0\,|\;\w\in\W\}$ has three hyperbolic copies of the base (which is the largest possible number).
We show that this situation persists for an interval $(\lb_-,\lb^+)\ni\lb_0$, and that the global attractor is a
hyperbolic copy of the base for $\lb\in(-\infty,\lb_*)\cup(\lb^*,\infty)$, with $\lb_*\le\lb_-$ and $\lb^*\ge\lb^+$.
These two inequalities are equalities in the minimal
case while, in general, $\lb_*=\inf_{\w\in\W}(\lb_-)_\w$, $\lb_-=\sup_{\w\in\W}(\lb_-)_\w$,
$\lb^+=\inf_{\w\in\W}(\lb^+)_\w$ and $\lb^*=\sup_{\w\in\W}(\lb^+)_\w$, where $((\lb_-)_\w,(\lb^+)_\w)$ is the
interval of persistence of three hyperbolic solutions for each $\w\in\W$.
So, the dynamics of $x'=\mh(\wt,x)$ is that of existence of three hyperbolic copies of the base
if $0\in(\lb_-,\lb^+)$, of existence of a unique (hyperbolic) copy of the base if $0\notin[\lb_*,\lb^*]$, and
$\w$-dependent if $0\in[\lb_*,\lb_-]\cup[\lb^+,\lb^*]$.
For the particular equation corresponding to a point $\w\in\W$, there are also several possibilities
outside the interval $((\lb_-)_\w,(\lb^+)_\w)$: a
unique bounded (and hyperbolic) solution for all $\lb<(\lb_-)_\w$ (which is the situation in the recurrent case),
or just for $\lb<(\lb_\diamond)_\w$ for a new value $(\lb_\diamond)_\w\in[\lb_*,(\lb_-)_\w)$, in which case there is
only one hyperbolic solution but more bounded solutions for $\lb\in[(\lb_\diamond)_\w,(\lb_-)_\w]$; and
a similar structure to the right of $(\lb^+)_\w$. In the d-concave case it makes sense saying that we have a double
generalized saddle-node nonautonomous bifurcation both for the family and each process.
Figure \ref{fig:diagramasdconcavos} and its caption explain these diagrams.

In one of the most usual stochastic formulations, an equation is subject to the variation of one of the involved
time-dependent coefficients, say $p(t)$, in the set $\mQ:=C(\R,\R^n)$. With the aim of analyzing the dependence with respect to $p$,
the evolution of the orbits is described in probabilistic terms, once a probabilistic measure
is fixed in the sigma-algebra of the Borel sets of $\mQ$ (see, for instance, Arnold \cite{arnol} and Oksendal \cite{oskendal}).
Our results provide a kind of deterministic counterpart of this type of analysis: despite the required concavity,
the set $\W$ can be a very large set of maps, which contain random bounded perturbations of many of their elements.
Let us illustrate this assertion with the help of the example given by the family
\begin{equation}\label{eq:introp}
 x'=h(t,x)+p(t)\,,
\end{equation}
where $p$ varies in a set $\mP$ of uniformly bounded and uniformly Lipschitz real continuous maps, i.e.,
$\mP:=\{p\in C(\R,\R)\,|\;\|p\|_\infty\leq k_1\,,\;\mathrm{Lip}(p)\leq k_2\}$; and where
$h$ is chosen to ensure concavity or d-concavity in measure on the skewproduct flow defined
over $\W:=\W_h\times\mP$. The continuous flow on $\mP$ defined by the time shift $p_s(t):=p(t+s)$,
i.e., the second component of the flow on $\W$, turns out to be chaotic in the sense of Devaney.
Combining methods of ergodic theory,
topological dynamics, and numerical analysis allows us to analyze the global dynamics of the
skewproduct flow, and to focus on the variation of the solutions of \eqref{eq:introp} with respect to $p$.
We show with suitable examples that the chaotic behavior has an important influence on the variation
of the behaviour of the trajectories of some mathematical models involving this type of coefficients
in population dynamics and electrical circuits.

We complete this introduction by briefly describing the paper structure.
Section \ref{sec:2} recalls some concepts and properties regarding skewproduct flows
and admissible processes, required throughout the paper. The core of the paper is
contained in Sections \ref{sec:3} and \ref{sec:4}, which respectively deal with the
concave and d-concave cases (both in measure) and have a similar structure. First, several previous results
prepare the way to the corresponding main bifurcation theorems (Theorems \ref{th:3Cbifur}
and \ref{th:4Dbifur}), formulated for skewproduct flows and then particularized to processes
(in Theorems \ref{th:3Cbifurproceso} and \ref{th:4Dbifurproceso}); and second, some applications
to population dynamics (in Section \ref{subsec:32}) and electrical circuit models (in Section \ref{subsec:41})
are described. Section \ref{subsec:31} proves that the the shift flows on the spaces $\mP$
of uniformly bounded and uniformly Lipschitz real continuous maps are chaotic, which plays
an essential role in demonstrating the wide range of applicability of our results.
The applications also show that the theoretical result of this paper
contribute to the understanding of critical transitions under less restrictive
conditions than in previous approaches. In particular, in Section \ref{subsec:42}, we
come back to the population dynamics model given a parametric family of equations including
equation \eqref{eq:huntingintrod}, which are asymptotically d-concave in measure,
in order to deduce the occurrence of a critical transition from the existence of a bifurcation value.

Each section begins with a more detailed description of its structure.
\section{Basic concepts and preliminary results}\label{sec:2}
We will proceed without recalling the definitions of the basic concepts of
maximal solutions of ordinary differential equations,
flows, orbits, invariant and minimal sets, \upalfa-limit and \upomeg-limit sets,
and invariant and ergodic measures. They can be found in \cite{cano}, \cite{dno1},
\cite{dno4}, and references therein.
Instead, we recall some less known concepts, also fundamental in this paper.
Section \ref{subsec:skew} concerns skewproduct flows, and Section \ref{subsec:process}
focuses on a particular process.
\subsection{Skewproduct flows given by families of scalar nonautonomous ODEs}
\label{subsec:skew}
Let $\sigma\colon\R\times\W\to\W$, $(t,\w)\mapsto\sigma(t,\w)=:\w{\cdot}t$ be a global
continuous flow on a compact metric space $\W$.
The set of continuous functions $\mh\colon\WR\to\R$ for which the
derivative $\mh_x$ with respect to the second variable exists and is continuous is
$C^{0,1}(\WR,\R)$, and $C^{0,2}(\WR,\R)$ is the subset of maps $\mh$ for which
$\mh_{xx}$ exists and is continuous. Each $\mh\in C^{0,1}(\WR,\R)$ provides
a family of scalar nonautonomous differential equations
\begin{equation}\label{eq:2fam}
x'=\mh(\w{\cdot}t,x)\,,\quad\w\in\W\,,
\end{equation}
with maximal solutions $(\alpha_{\w,x},\beta_{\w,x})\to\R,\,t\mapsto v(t,\w,x)$
satisfying $v(0,\w,x)=x$. So, $v(t+s,\w,x)=v(t,\ws,v(s,\w,x))$ whenever the right-hand term is defined, and hence,
if $\mV:=\bigcup_{(\w,x)\in\WR}((\alpha_{\w,x},\beta_{\w,x})\times\{(\w,x)\})$, then
\begin{equation}\label{def:2tau}
 \tau\colon\mV\subseteq\R\times\WR\to\WR\,,\;\;(t,\w,x)\mapsto (\wt,v(t,\w,x))
\end{equation}
defines a (possibly local) continuous flow on $\WR$, of {\em skewproduct} type.
The interested reader can consult in \cite[Section~2.2]{dno4}
how families of this type appear in a natural way
from a suitable single nonautonomous equation, like those considered in
Section \ref{subsec:process}, by means of the hull construction.

A {\em $\tau$-equilibrium} is a map $\mb\colon\W\to\R$ whose graph is  $\tau$-invariant;
i.e., with $v(t,\w,\mb(\w))=\mb(\wt)$ for all $\w\in\W$ and $t\in\R$.
A {\em $\tau$-copy of the base\/} or {\em $\tau$-copy of $\W$} is the compact graph
of a continuous $\tau$-equilibrium $\mb$. We represent it by $\{\mb\}$. Sometimes, a
copy of the base is identified with the equilibrium that provides it. The prefix $\tau$
will be often omitted.

Given a bounded $\tau$-invariant set $\mB\subset\W\times[-k,k]\subset\WR$ projecting onto $\W$ (with $k>0$),
the maps $\w\mapsto\inf\{x\in\R\,|\;(\w,x)\in\mB\}$ and $\w\mapsto\sup\{x\in\R\,|\;(\w,x)\in\mB\}$
define $\tau$-equilibria. We will refer to these maps as
the {\em lower\/} and {\em upper equilibria of $\mB$}. They are respectively
lower and upper semicontinuous if $\mB$ is compact, and hence they are
$m$-measurable for any ergodic measure in $\W$.

Two compact subsets $\mK_1$ and $\mK_2$ of $\W\times\R$ are {\em ordered\/} with
$\mK_1<\mK_2$ if $x_1<x_2$ whenever there exists $\w\in\W$ such that
$(\w,x_1)\in\mK_1$ and $(\w,x_2)\in\mK_2$. Two (ordered)
bounded equilibria $\mb_1,\mb_2$ with $\mb_1<\mb_2$ (i.e., with $\mb_1(\w)<\mb_2(\w)$
for all $\w\in\W$) are {\em uniformly separated\/} if
$\inf_{\w\in\W}(\mb_2(\w)-\mb_1(\w))>0$.

If there exists a compact $\tau$-invariant set $\mA\subset\WR$ with
$\lim_{t\to\infty} \text{dist}(\mC{\cdot}t,\mA)=0$ for every bounded set $\mC$,
where $\mC{\cdot}t=\{(\wt,v(t,\w,x))\,|\;(\w,x)\in\mC\}$ and
\[
\text{dist}(\mC_1,\mC_2)=\sup_{(\w_1,x_1)\in\mC_1}\left(\inf_{(\w_2,x_2)\in\mC_2}
\big(\mathrm{dist}_{\WR}((\w_1,x_1),(\w_2,x_2))\big)\right)\,,
\]
then $\mA$ is the (unique) {\em global attractor for $\tau$}.

A $\tau$-copy of the base $\{\mb\}$ is
{\em hyperbolic attractive\/} if it is uniformly exponentially stable (on the fiber) as time increases;
i.e., if there exists $\rho>0$, $k\ge 1$ and $\gamma>0$ such that: if, for any $\w\in\W$,
$|\mb(\w)-x|<\rho$, then $v(t,\w,x)$ is defined for all $t\ge 0$, and in addition
$|\mb(\wt)-v(t,\w,x)|\le k\,e^{-\gamma\,t}\,|\mb(\w)-x|$ for $t\ge 0$. In this case,
$\rho$ is a {\em radius of uniform exponential stability} and $(k,\gamma)$ is a
{\em dichotomy constant pair} of $\mb$. Changing $t\ge 0$ by $t\le 0$ and $\gamma$ by $-\gamma$
provides the definition of {\em repulsive hyperbolic} $\tau$-copy of the base. We will usually
denote $\tilde\mb:=\mb$ if $\{\mb\}$ is a hyperbolic $\tau$-copy of the base.

Let $\mK\subset\WR$ be $\tau$-invariant compact set projecting onto $\W$,
and let $\mathfrak{M}_\mathrm{inv}(\mK,\tau)$ and $\mathfrak{M}_\mathrm{erg}(\mK,\tau)$
be the (nonempty) sets of the $\tau$-invariant and $\tau$-ergodic measures on $\mK$.
A value $\gamma\in\R$ is a \emph{Lyapunov exponent of $\mK$} if there exists
$(\w,x)\in\mK$ such that $\gamma=\lim_{t\to\pm\infty}(1/t)\int_0^t\mh_x(\tau(r,\w,x))\,dr$.
In this case, there exists $\nu\in\mathfrak{M}_\mathrm{inv}(\mK,\tau)$ such that
$\gamma=\int_\mK \mh_x(\w,x)\,d\nu$: this fact can be deduced from Riesz Representation
Theorem and Kryloff-Bogoliuboff’s Theorem.
In addition, Birkhoff's Ergodic Theorem ensures that $\gamma(\mK,\nu):=\int_\mK \mh_x(\w,x)\,d\nu$
is a Lyapunov exponent of $\mK$ for each $\nu\in\mathfrak{M}_\mathrm{erg}(\mK,\tau)$.
Since the ergodic measures are the extreme points in the set of invariant measures,
the upper and lower Lypaunov exponents of $\mK$ are $\gamma(\mK,\nu^u)$ and
$\gamma(\mK,\nu^l)$ for suitable measures $\nu^u,\nu^l\in \mathfrak{M}_\mathrm{erg}(\mK,\tau)$.
According to \cite[Theorem 4.1]{furstenberg1} and \cite[Theorem~1.8.4]{arnol},
if $\nu\in\mathfrak{M}_\mathrm{erg}(\mK,\tau)$ projects onto $m\in \merg$,
then there exists an $m$-measurable $\tau$-equilibrium $\mb\colon\W\to\R$ with graph in $\mK$
such that $\gamma(\mK,\nu)=\int_\W \mh_x(\w,\mb(\w))\,dm$.
Therefore, there exist $m^l,m^u\in\merg$, an $m^l$-measurable
equilibrium $\mb^l\colon\W\to\R$ and an $m^u$-measurable equilibrium $\mb^u\colon\W\to\R$
such that the lower and upper Lyapunov exponents of $\mK$ are given by
$\int_\W \mh_x(\w,\mb^l(\w))\,dm^l$ and $\int_\W \mh_x(\w,\mb^u(\w))\,dm^u$, respectively.
Finally, if $m\in\merg$ and $\mb\colon\W\to\R$ is an $m$-measurable $\tau$-equilibrium
with graph in $\mK$, then $\int_\W \mh_x(\w,\mb(\w))\,dm$ is one of the
Lyapunov exponents of $\mK$.

There is a strong connection between the hyperbolicity of a copy of the base $\{\mb\}$
and the sign of all its Lyapunov exponents: $\{\mb\}$ is hyperbolic attractive (resp.~repulsive)
if and only if all its Lyapunov exponents are strictly negative (resp.~positive).
This property is a consequence of \cite[Theorem 2.8]{dno4}. We repeat here its statement,
since it will be fundamental in some key points of the proofs.
\begin{teor}\label{th:2copia}
Let $\mK\subset\WR$ be a $\tau$-invariant compact set projecting
onto $\W$. Assume that its upper and lower equilibria
coincide (at least) on a point of each minimal subset $\mM\subseteq\W$.
Then, all the Lyapunov exponents of $\mK$
are strictly negative (resp.~positive) if and only if $\mK$ is
an attractive (resp.~repulsive) hyperbolic $\tau$-copy of the base.
In addition, if either $\mK$ (and hence $\W$) is minimal or its upper and lower equilibria
coincide on a $\tau$-invariant subset $\W_0\subseteq\W$
with $m(\W_0)=1$ for all $m\in\merg$, then the condition on its upper and lower equilibria holds.
\end{teor}
We will very often restrict the family \eqref{eq:2fam} to the elements of a $\sigma$-invariant
compact subset $\W_0\subseteq\W$. It is easy to check that $C^{0,1}(\WR,\R)\hookrightarrow C^{0,1}(\W_0\times\R,\R)$, so
that it makes perfect sense to consider the restricted skewproduct flow on $\W_0\times\R$,
whose definition is the same. In many of these cases, $\W_0$ will be the closure $\W_{\bar\w}$
of the $\sigma$-orbit $\{\bwt\,|\;t\in\R\}$ of a point $\bar\w\in\W$. At this point, we recall that the
base flow is {\em transitive\/} if there exists $\bar\w\in\W$ with $\W_{\bar\w}=\W$, and {\em minimal\/}
if $\W_\w=\W$ for all $\w\in\W$.

\subsection{A particular process}\label{subsec:process}
A map $h\colon\RR\to\R$ is {\em admissible\/}, or $h\in C^{0,0}(\RR,\R)$,
if the restriction of $h$ to $\R\times\mJ$ is bounded
and uniformly continuous for any compact set $\mJ\subset\R$. The map $h$
is {\em $C^1$-admissible\/} (resp.~{\em $C^2$-admissible\/}), or
$h\in C^{0,1}(\RR,\R)$ (resp.~$h\in C^{0,2}(\RR,\R)$), if it is admissible and its
derivative $h_x$ with respect to the second variable exists and is admissible
(resp.~$h_x$ and $h_{xx}$ exist and are admissible).

Let us take $h\in C^{0,1}(\RR,\R)$ and let $x_h(t,s,x)$ be the maximal solution of
\begin{equation}\label{eq:2proceso}
 x'=h(t,x)
\end{equation}
with $x_h(s,s,x)=x$. By uniqueness of solutions, $x_h(t,s,x_h(s,r,x))=x_h(t,r,x)$
whenever the left-hand term is defined. The map
$(t,s,x)\mapsto x_h(t,s,x)$ is often called a {\em process}, or an {\em admissible process}.

Two bounded (and hence globally defined) solutions $b_1(t)$ and
$b_2(t)$ of \eqref{eq:2proceso} are {\em uniformly separated\/} if
$\inf_{t\in\R}|b_1(t)-b_2(t)|>0$.

A bounded solution
$\tilde b(t)$ of \eqref{eq:2proceso} is {\em hyperbolic attractive\/}
(resp.~{\em hyperbolic repulsive\/}) if
there exist $k\ge 1$ and $\gamma>0$ such that
$\exp\Big(\int_s^t h_x(r,\tilde b(r))\,dr\Big)\le k\,e^{-\gamma(t-s)}$ whenever $t\ge s$
(resp.~$\exp\Big(\int_s^t h_x(r,\tilde b(r))\,dr\Big)\le k\,e^{\gamma (t-s)}$ whenever $t\le s$).
An attractive (resp.~repulsive) hyperbolic solution is uniformly
exponentially stable at $+\infty$ (resp.~at $-\infty$),
as deduced from the First Approximation Theorem (see Theorem~III.2.4 of \cite{hale});
i.e., there exist a {\em radius of uniform exponential stability\/} $\rho>0$ and a {\em dichotomy
constant pair\/} $(k,\gamma)$ with $k\ge 1$ and $\gamma>0$ such that,
if $|\tilde b(s)-x|\le\rho$, then $|\tilde b(t)-x_h(t,s,x)|
\le k\,e^{-\gamma(t-s)}|\tilde b(t)-x|$ for all $t\ge s$
(resp.~$|\tilde b(t)-x_h(t,s,x)|\le k\,e^{\gamma(t-s)}|\tilde b(s)-x|$
for all $t\le s$).

Let us take as starting point the family \eqref{eq:2fam}, with $\mh\in C^{0,1}(\WR,\R)$.
We fix $\bar\w\in\W$ and define $h_{\bar\w}(t,x):=\mh(\bwt,x)$. It is easy to check that
$h_{\bar\w}\in C^{0,1}(\RR,\R)$, and hence
\begin{equation}\label{eq:2una}
 x'=h_{\bar\w}(t,x) \qquad (\text{i.e., }x'=\mh(\bwt,x)\,)
\end{equation}
defines an admissible process. In addition, if $x_{\bar\w}(t,s,x)$ is the solution of \eqref{eq:2una} 
with $x_{\bar\w}(s,s,x)=x$, then
$x_{\bar\w}(t,s,x):=v(t-s,\bws,x)$. We will use often this fundamental relation, as well as the next
result, which reproduces \cite[Proposition 2.7]{dno4}:
\begin{prop}\label{prop:2extiende}
Let $\mh\in C^{0,1}(\WR,\R)$, and let $\tilde\mb\colon\W\to\R$ determine an attractive
(resp.~repulsive) copy of the base for \eqref{eq:2fam}. For any $\bar\w\in\W$, the function
$\tilde b_{\bar\w}$ defined by $\tilde b_{\bar\w}(t):=\tilde\mb(\bwt)$ is an attractive (resp.~repulsive)
hyperbolic solution of \eqref{eq:2una}.
\end{prop}

For the reader's convenience, we state the next fundamental persistence result,
very often used in the proofs of our results. The first assertion, concerning processes,
is \cite[Theorem 2.3]{dno4}, and the interested reader can find there suitable references for its proof.
This proof can be easily adapted to prove the second assertion, about skewproduct flows.
(The interested reader can find the details in \cite[Section 1.3.2]{duen}.)
For $h\in C^{0,1}(\RR,\R)$, we
denote $\n{h}_{1,\kappa}:=\sup_{(t,x)\in\R\times[-\kappa,\kappa]} |h(t,x)|+
\sup_{(t,x)\in\R\times[-\kappa,\kappa]} |h_x(t,x)|$; and similarly, for $\mh\in C^{0,1}(\WR,\R)$, we denote
$\n{\mh}_{1,\kappa}:=\max_{(\w,x)\in\W\times[-\kappa,\kappa]} |\mh(\w,x)|+
\max_{(\w,x)\in\W\times[-\kappa,\kappa]} |\mh_x(\w,x)|$.
\begin{teor}\label{th:2persistencia}
Let $h\in C^{0,1}(\RR,\R)$ be fixed, let $\tilde b_h$ be an attractive (resp. repulsive) hyperbolic solution
of \eqref{eq:2proceso} with dichotomy constant pair $(k_0,\gamma_0)$, and take $\kappa>\sup_{t\in\R}|\tilde b_h(t)|$.
Then, for every $\gamma\in(0,\gamma_0)$ and $\ep>0$, there exists $\delta_\ep>0$ and $\rho_\ep>0$ such that,
if $g$ is $C^1$-admissible and $\n{h-g}_{1,\kappa}<\delta_\ep$, then
there exists an attractive (resp.~repulsive) hyperbolic solution $\tilde b_g$ of $x'=g(t,x)$
with radius of uniform exponential stability dichotomy $\rho_\ep$ and dichotomy constant pair $(k_0,\gamma)$
such that $\sup_{t\in\R}|\tilde b_h(t)-\tilde b_g(t)|<\ep$.

Let $\mh\in C^{0,1}(\WR,\R)$ be fixed, let $\{\tilde\mb_\mh\}$ be an attractive (resp. repulsive) hyperbolic
copy of the base for \eqref{eq:2fam} with dichotomy constant pair $(\gamma_0,k_0)$, and take
$\kappa>\max_{\w\in\W}|\tilde\mb_\mh(\w)|$.
Then, for every $\gamma\in(0,\gamma_0)$ and $\ep>0$, there exists $\delta_\ep>0$ and $\rho_\ep>0$ such that,
if $\mg\in C^{0,1}(\WR,\R)$ and $\|\mh-\mg\|_{1,\kappa}<\delta_\ep$, then
there exists an attractive (resp. repulsive) hyperbolic copy of the base $\{\tilde\mb_\mg\}$
of $x'=\mg(\wt,x)$ with radius of uniform exponential stability $\rho_\ep$ and dichotomy constant pair $(k_0,\gamma)$
such that $\max_{\w\in\W}|\tilde\mb_\mh(\w)-\tilde\mb_\mg(\w)|<\ep$.
\end{teor}

\section{A global bifurcation diagram in the concave in measure case}\label{sec:3}
\begin{figure}[h]
\centering
\includegraphics[width=\textwidth]{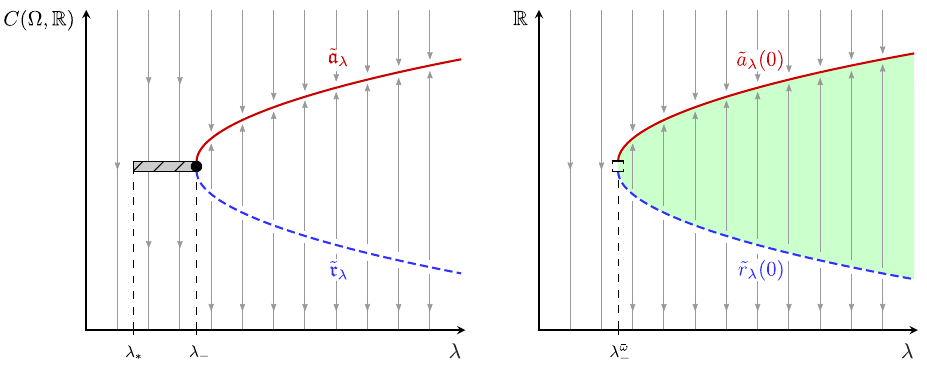}
\caption{In the left panel, the generalized saddle-node bifurcation diagram of the $\lb$-parametric
family $x'=\mh(\wt,x,\lb)$, $\w\in\W$, described in Theorem \ref{th:3Cbifur}. The strictly increasing solid red curve
(resp. strictly decreasing dashed blue curve) over $(\lb_-,\infty)$ represents the family of attractive
(resp. repulsive) hyperbolic copies of the base $\tilde\ma_\lb$ (resp. $\tilde\mr_\lb$). At $\lb_-$,
the large black point represents the maps $\ma_{\lb_-}:=\lim_{\lb\downarrow\lb_-}\tilde\ma_\lb$ and $\mr_{\lb_-}:=\lim_{\lb\downarrow\lb_-}\tilde\mr_\lb$, which may not belong to $C(\W,\R)$:
they are upper and lower semicontinuous respectively and satisfy $\inf_{\w\in\W}(\ma_{\lb_-}(\w)-\mr_{\lb_-}(\w))=0$.
There are no bounded solutions for $\lb<\lb_*\le\lb_-$. If $\lb_*<\lb_-$ (as drawn here) and $\lb\in[\lb_*,\lb_-)$,
the domain $\W_\lb^b$ of the semicontinuous maps $\ma_\lb$ and $\mr_\lb$ is a proper subset of $\W$ and,
for any $\w\in\W$, $\lb\mapsto a_\lb(\w),-r_\lb(\w)$ are defined and strictly increasing on a closed
positive $\w$-dependent halfline. These facts are depicted by the gray striped rectangle over $[\lb_*,\lb_-)$. The
light gray arrows partly depict the dynamics of the rest of the orbits.
\vspace{.15cm}\newline
In the right panel, the saddle-node bifurcation diagram of the $\lb$-parametric family of
nonautonomous equations $x'=\mh(\bar\w{\cdot}t,x,\lb)$
described in Theorems \ref{th:3Cbifuruna} and \ref{th:3Casint}, with $\lb_-^{\bar\w}\in[\lb_*,\lb_-]$
as unique bifurcation point. The strictly increasing solid red curve (resp. strictly decreasing dashed blue curve) on
$(\lb_-^{\bar\w},\infty)$ represents the value at $t=0$ of the upper (resp.~lower) bounded solution
$\tilde a_\lb$ (resp.~$\tilde r_\lb$), which is hyperbolic attractive (resp.~repulsive).
The solutions $a_{\lb_-^{\bar\w}}$ and $r_{\lb_-^{\bar\w}}$ are globally defined, can coincide, and
are never uniformly separated: the vertical interval $[r_{\lb_-^{\bar\w}}(0),a_{\lb_-^{\bar\w}}(0)$] over
$\lb_-^{\bar\w}$ can either be nondegenerate (as drawn here) or reduce to a point.
There are no bounded solutions for $\lb<\lb_-^{\bar\w}$.
The green-shadowed area represents the initial data of bounded solutions, and
the light gray arrows give some clues about the dynamics of the rest of the solutions.}
\label{fig:diagramasconcavos}
\end{figure}

This section focuses on bifurcation diagrams for $\lambda$-parametric families
of skewproduct flows and of admissible processes determined by scalar ODEs for which the
corresponding laws satisfy certain concavity properties: they are concave in measure.
The dynamical possibilities for a given skewproduct flow or process follow from the
bifurcation diagram. The main theoretical results are
presented in Theorems \ref{th:3Cbifur}, \ref{th:3Cbifuruna} and \ref{th:3Casint}, and
graphically explained in Figure \ref{fig:diagramasconcavos}. In Section \ref{subsec:32},
these results are applied to some population dynamics models, showing their high applicability in
two senses: they work for deterministic models for which the large size of the set
in which their coefficients vary (described in Section \ref{subsec:31}) allows
us to consider them as ``almost stochastic''; and they can be used to analyze the occurrence
of critical transitions.

Let $(\W,\sigma)$ be a global continuous real flow on a compact metric space,
and let us consider the family of scalar ordinary differential equations
\begin{equation}\label{eq:3ini}
 x'=\mh(\wt,x)\,,\quad\w\in\W\,,
\end{equation}
where $\mh\colon\WR\to\R$ satisfies (all or part of) the next conditions:
\begin{enumerate}[leftmargin=20pt,label=\rm{\bf{c\arabic*}}]
\item\label{c1} $\mh\in C^{0,1}(\WR,\R)$,
\item\label{c2} $\limsup_{x\to\pm\infty}\mh(\w,x)<0$ uniformly on $\W$.
\item\label{c3} $m(\{\w\in\W\,|\;x\mapsto\mh(\w,x) \text{ is concave}\})=1$ for all $m\in\merg$,
\item\label{c4} $m(\{\w\in\W\,|\; x\mapsto\mh_x(\w,x)$ is strictly decreasing on $\mJ\})>0$
for all compact interval $\mJ\subset\R$ and all $m\in\merg$,
\end{enumerate}
where $\merg$ is the (nonempty) set of $\sigma$-ergodic measures on $\W$. Recall that
$f\colon\R\to\R$ is {\em concave\/} if $f(\alpha\,x_1+(1-\alpha)\,x_2)\ge \alpha\,f(x_1)+(1-\alpha)\,f(x_2)$
for all $\alpha\in[0,1]$ and $x_1,x_2\in\R$.
We refer to the case of the equation \eqref{eq:3ini} under conditions
\ref{c1} and \ref{c3} as the {\em concave in measure case} (or just concave case, to
simplify language). Note that \ref{c4} precludes the {\em linear case\/}
$\mh(\w,x)=\mc(\w)\,x+\md(\w)$.
\begin{nota}\label{rm:3Ctambien}
Let $\W_0\subset\W$ be a nonempty compact $\sigma$-invariant subset. Then, any
$m_0\in{\mathfrak{M}_\mathrm{erg}(\W_0,\sigma)}$ can be extended to $m\in\merg$
by $m(\mU)=m_0(\mU\cap\W_0)$. So, if $\mh$ satisfies \hyperlink{c1}{\bf cj} for
$j\in\{1,2,3,4\}$, also the restriction $\mh\colon\W_0\times\R\to\R$ satisfies
\hyperlink{c1}{\bf cj}.
\end{nota}

Let $\tau$ be the skewproduct flow defined by \eqref{def:2tau}, with $\tau(t,\w,x)=(\wt,v(t,\w,x))$.
We say that {\em there exists an attractor-repeller pair $(\tma,\tmr)$
of copies of the base\/} (or {\em copies of\/} $\W$)
{\em for $\tau$} (or {\em for \eqref{eq:3ini}}), or that {\em $(\tma,\tmr)$ is an
attractor-repeller pair of $\tau$-copies of the base},
if $\{\tma\}$ is an attractive hyperbolic $\tau$-copy of $\W$ and $\{\tmr\}$ is a repulsive
hyperbolic $\tau$-copy of $\W$. Its existence under conditions \ref{c1},
\ref{c3} and \ref{c4} is characterized in
\cite[Theorem 3.3]{dno4} in terms of the existence of two disjoint and
ordered $\tau$-invariant compact sets $\mK_1<\mK_2$ projecting onto $\W$.
This result also proves that $\tmr<\tma$.

In addition, assuming \ref{c1} and \ref{c2},
\cite[Proposition 3.5]{dno4} proves several properties of
the (possibly empty and otherwise bounded) set of globally bounded orbits
\[
 \mB:=\big\{(\w,x)\,|\;\sup_{t\in\R}|v(t,\w,x)|<\infty\big\}\,,
\]
its (possibly empty and otherwise compact) projection $\W^b$ over $\W$, and the maps
$\mr,\ma\colon\W^b\to\R$ such that
\[
 \mB=\{(\w,x)\,|\;\w\in\W^b\text{ and }\mr(\w)\le x\le\ma(\w)\}\,,
\]
which are lower and upper semicontinuous equilibria for the restriction of
$\tau$ to $\W^b\times\R$. It also establishes some comparison properties which
will play a fundamental role in the proofs of our main results.

We will analyze the $\lb$-bifurcation diagram for a parametric family
\begin{equation}\label{eq:3Chlb}
 x'=\mh(\wt,x,\lb)\,,\quad\w\in\W\,,
\end{equation}
where $\lb$ varies in an open interval $\mI$ of $\R$ (often $\mI=\R$).
We will represent by \eqref{eq:3Chlb}$^\lb$ the family for a fixed value of the parameter,
by \eqref{eq:3Chlb}$_\w$ the $\lb$-parametric family of processes corresponding to a fixed point
$\w\in\W$, and by \eqref{eq:3Chlb}$^\lb_\w$ a particular process. (We will do the same
with the remaining families of equations of the paper.) Recall that
$\W_\w\subseteq\W$ is the closure of the $\sigma$-orbit $\{\wt\,|\:t\in\R\}$ for
each $\w\in\W$.

For now, we assume that the map $(\w,x)\mapsto\mh(t,\w,\lb)$
satisfies \ref{c1} and \ref{c2} for all $\lb\in\mI$.
Let $\tau_\lb\colon\mV_\lb\subseteq\R\times\WR\to\WR,\,(t,\w,x)\mapsto(\wt,v_\lb(t,\w,x))$ be the
corresponding skewproduct flow, let $\mB_\lb\subset\WR$ be the
set of bounded $\tau_\lb$-orbits, and, if it is nonempty,
let $\W_\lb^b$ be its (compact and $\tau_\lb$-invariant)
projection onto $\W$ and let $\mr_\lb,\ma_\lb\colon\W_\lb^b\to\R$ be
the bounded ${\tau_\lb}|_{\W_\lb^b}$-equilibria such that
$\mB_\lb=\bigcup_{\w\in\W_\lb^b}\big(\{\w\}\times[\mr_\lb(\w),\ma_\lb(\w)]\big)$.
For further purposes,
we point out that if \ref{c1} and \ref{c2} hold, then the set of
bounded trajectories of the process \eqref{eq:3Chlb}$_{\bar\w}^\lb$ is nonempty if
and only if $\bar\w\in\W_\lb^b$, in which case the upper and lower bounded solutions
are $t\mapsto\ma_\lb(\bwt)$ and $t\mapsto\mr_\lb(\bwt)$.
The set $C(\W,\R)$ is endowed with the uniform topology.

Theorem \ref{th:3Cbifur} provides an
extension of \cite[Theorem 3.1]{nuob6} to our more general
map $\mh$ and, most significantly, to a family of skewproduct flows~$(\WR,\tau_\lb)$
over a base $(\W,\sigma)$ which is not necessarily minimal or even transitive.
Proposition \ref{prop:3Csolocoer}, whose hypotheses are less restrictive,
shows that the set of bounded solutions increases as the law increases.
\begin{prop}\label{prop:3Csolocoer}
Let $(\w,x)\mapsto\mh(\w,x,\lb)$ satisfy \ref{c1} and \ref{c2} for all
$\lb\in\mI\subseteq\R$, and assume that $\mh(\w,x,\lb_1)<\mh(\w,x,\lb_2)$ for all $(\w,x)\in\WR$
if $\lb_1,\lb_2\in\mI$ and $\lb_1<\lb_2$. If $\mB_{\lb_1}$ is nonempty,
$\lb_1,\lb_2\in\mI$ and $\lb_1<\lb_2$, then
$\mB_{\lb_1}\subsetneq\mB_{\lb_2}$. In fact,
\begin{equation}\label{des:3Cchain}
 \mr_{\lb_2}(\w)<\mr_{\lb_1}(\w)\le\ma_{\lb_1}(\w)<\ma_{\lb_2}(\w)\quad
 \text{if $\w\in\W^b_{\lb_1}$, $\lb_1,\lb_2\in\mI$ and $\lb_1<\lb_2$}\,.
\end{equation}
Assume that, in addition, 
there exist $\lb^-,\lb^+\in\mI$ such that $\mB_{\lb^-}$ is empty and $\mB_{\lb^+}$
is nonempty, and that $\WR\times\mI\to\R,\,(\w,x,\lb)\mapsto \mh(\w,x,\lb)$ is
jointly continuous. Then, $\lb_*:=\inf\{\lb\in\mI\,|\;\mB_\lb\text{ is nonempty}\}$
belongs to $\mI$, and $\mB_{\lb_*}$ is nonempty.
\end{prop}
\begin{proof}
All the parameters are taken in $\mI$.
We assume the existence of $x_1\in(\mB_{\lb_1})_{\w_1}$, call $b(t):=v_{\lb_1}(t,\w_1,x_1)$, take $\lb_2>\lb_1$,
observe that $b'(t)<\mh(\wt,b(t),\lb_2)$, and deduce from \cite[Proposition 3.5(v)]{dno4} that $\w_1\in\W^{\lb_2}_b$ and
$\mr_{\lb_2}(\w_1)<b(0)=x_1<\ma_{\lb_2}(\w_1)$; that is, \eqref{des:3Cchain} holds, and it
ensures that $\mB_{\lb_1}\subsetneq\mB_{\lb_2}$.
Under the additional conditions, it is clear that
$\lb_*\ge\lb^-$ belongs to $\mI$. To check the last assertion, we take a sequence $(\lb_n)\downarrow\lb_*$.
For each $n\in\N$, we take $(\w_n,x_n)\in\mB_{\lb_n}$, and deduce from \eqref{des:3Cchain}
the existence of the limit $(\w_*,x_*)$ of a suitable subsequence of $(\w_n,x_n)$.
It is easy to deduce from \eqref{des:3Cchain} and the joint continuity of $\mh$
that $v_{\lb_*}(t,\w_*,x_*)$ is bounded.
\end{proof}
Some points of the proof of the next theorem, and even its statement (vi), require
the information provided by Remark~\ref{rm:3Ctambien}, as well as the fact that the
assumed properties on $\mh$ also hold for the restriction of the maps to
$\W_0\times\R\times\mI$ for any nonempty compact $\sigma$-invariant subset $\W_0\subseteq\W$.
\begin{teor}\label{th:3Cbifur}
Assume the next conditions on $\mh\colon\WR\times\R\to\R$:
\begin{itemize}[leftmargin=20pt]
\item[-] $(\w,x)\mapsto\mh(\w,x,\lb)$ satisfy \ref{c1}, \ref{c2},
\ref{c3} and \ref{c4} for all $\lb\in\R$;
\item[-] $\mh$ and $\mh_x$ are jointly continuous on $\WR\times\R$;
\item[-] $\inf_{(\w,x,\lb_1,\lb_2)\in\W\times[-j,j]\times\R^2,\,\lb_1<\lb_2}
(\mh(\w,x,\lb_2)-\mh(\w,x,\lb_1))/(\lb_2-\lb_1)=:c_j>0$ for each $j\in\N$.
\end{itemize}
Then, there exist real values $\lb_*,\,\lb_-$ with $\lb_*\le\lb_-$
such that
\begin{itemize}[leftmargin=20pt]
\item[\rm(i)] $\W_\lb^b=\W$ and $(\tma_\lb,\tmr_\lb):=(\ma_\lb,\mr_\lb)$ is an attractor-repeller
    pair of copies of the base for \eqref{eq:3Chlb}$^\lb$ if and only if $\lb>\lb_-$.
    Besides,
    the maps $(\lb_-,\infty)\to C(\W,\R),\lb\mapsto-\tmr_\lb,\,\tma_\lb$ are continuous, strictly increasing,
    and $\lim_{\lb\to\infty}\tmr_\lb=\lim_{\lb\to\infty}\tma_\lb=\infty$ uniformly on $\W$.
\item[\rm(ii)] $\W^b_{\lb_-}=\W$; $\mr_{\lb_-}=\lim_{\lb\to(\lb_-)^+}\tmr_\lb$
    and $\ma_{\lb_-}=\lim_{\lb\to(\lb_-)^+}\tma_\lb$ pointwise on $\W$;
    and $\inf_{\w\in\W}(\ma_{\lb_-}(\w)-\mr_{\lb_-}(\w))=0$.
\item[\rm(iii)] $\W_{\lb}^b$ is empty if and only if $\lb<\lb_*$.
\item[\rm(iv)] If $\lb_*<\lb_-$ and $\lb\in[\lb_*,\lb_-)$, then $\W_{\lb}^b$ is a proper subset of $\W$.
\end{itemize}
In addition,
\begin{itemize}[leftmargin=20pt]
\item[\rm(v)] if $\W$ is minimal, then $\lb_*=\lb_-$.
\item[\rm(vi)] For each $\w\in\W$, let $\lb_-^\w$ be the upper special value of the parameter associated
    to the restriction of the $\lb$-parametric family \eqref{eq:3Chlb} to the closure $\W_\w$ of
    $\{\wt\,|\;t\in\R\}$. Then, $\lb_*=\inf_{\w\in\W}\lb_-^\w$ and $\lb_-=\sup_{\w\in\W}\lb_-^\w$.
\end{itemize}%
\end{teor}
\begin{proof}
(i)-(iv) The conditions on
$\mh$ ensure $\mh(\w,x,\lb_1)<\mh(\w,x,\lb_2)$ if $\lb_1<\lb_2$ and that, for each $j\in\N$,
$\mh(\w,x,\lb)\ge -d_j+\lb\,c_j$ if $x\in[-j,j]$ and $\lb>0$ and
$\mh(\w,x,\lb)\le d_j+\lb\,c_j$ if $x\in[-j,j]$ and $\lb<0$,
where $d_j:=\max_{(\w,x)\in\W\times[-j,j]}|\mh(\w,x,0)|$.
We take $j\in\N$ and $\lb_j$ such that $-d_j+c_j\lb>0$ for all $\lb\ge\lb_j$.
This means that $\mh(\w,\pm j,\lb)>0=(\pm j)'$ for all $\w\in\W$ if $\lb\ge\lb_j$.
So, \cite[Proposition 3.5(v)]{dno4} ensures that $\W_\lb^b=\W$ and that
$\mr_\lb(\w)<-j<j<\ma_\lb(\w)$ for all $\w\in\W$. Hence, the closures of
$\{(\w,\mr_\lb(\w))\,|\;\w\in\W\}$ and $\{(\w,\ma_\lb(\w))\,|\;\w\in\W\}$
define two disjoint and ordered compact $\tau_\lb$-invariant sets.
Therefore, \cite[Theorem 3.3]{dno4}  shows that $(\tma_\lb,\tmr_\lb):=(\ma_\lb,\mr_\lb)$ is an
attractor-repeller pair of copies of the base for \eqref{eq:3Chlb}$^\lb$. Note that
we have also proved that $\lim_{\lb\to\infty}\tmr_\lb(\w)=-\infty$ and
$\lim_{\lb\to\infty}\tma_\lb(\w)=\infty$ uniformly on $\W$.

In addition, Proposition \ref{prop:3Csolocoer} ensures \eqref{des:3Cchain}.
So, for any $\lb<0$, $\mB_\lb\subseteq\mB_{0}\subset\W\times[-k,k]$ for a large enough $k\in\N$.
We take $\lb_0<0$ such that $\mh(\w,x,\lb_0)\le d_k+\lb_0\,c_k\le-1$ if $x\in[-k,k]$,
deduce that any solution of \eqref{eq:3Chlb}$^{\lb_0}$ leaves the set $[-k,k]$,
and conclude that $\mB_{\lb_0}$ is empty.
So, all the hypotheses of Proposition \ref{prop:3Csolocoer} hold.
Let us define $\mI:=\{\lb\in\R\,|\;\W^b_\lb=\W\text{ and }(\ma_\lb,\mr_\lb)$
is an attractor-repeller pair of copies of the base for \eqref{eq:3Chlb}$^\lb\}$.
Using again \eqref{des:3Cchain} and \cite[Theorem 3.3]{dno4} , we check that
$[\lb,\infty)\subset\mI$ if $\lb\in\mI$. We define $\lb_*>\lb_0$ as in
Proposition \ref{prop:3Csolocoer}, and $\lb_-:=\inf\mI\ge\lb_*$. These are the constants of the
statement. The continuity of the maps is a consequence of the robustness of the hyperbolicity: see
Theorem \ref{th:2persistencia}.
For next purposes, we point out that this property also ensures that $\lb_-\notin\mI$.
The monotonicity properties follow from Proposition \ref{prop:3Csolocoer}, and this completes the proof of (i).
Note that (iii) is also proved by Proposition \ref{prop:3Csolocoer}.

It follows from \eqref{des:3Cchain} that $r_{\lb_-}:=\lim_{\lb\to(\lb_-)^+}\mr_\lb$
and $a_{\lb_-}:=\lim_{\lb\to(\lb_-)^+}\ma_\lb$
are globally defined and bounded maps, and that $[r_{\lb_-}(\w),a_{\lb_-}(\w)]\supseteq
[\mr_{\lb_-}(\w),\ma_{\lb_-}(\w)]$ for all $\w\in\W^b_{\lb_-}$. Clearly,
$v_{\lb_-}(t,\w,r_{\lb_-}(\w))=r_{\lb_-}(\wt)$ and $v_{\lb_-}(t,\w,a_{\lb_-}(\w))=
a_{\lb_-}(\wt)$ for all $(\w,t)\in\WR$, from where we conclude that $\W^b_{\lb_-}=\W$ and
$[r_{\lb_-}(\w),a_{\lb_-}(\w)]\subseteq [\mr_{\lb_-}(\w),\ma_{\lb_-}(\w)]$.
That is, $r_{\lb_-}=\mr_{\lb_-}$ and $a_{\lb_-}=\ma_{\lb_-}$.
To complete the proof of (ii), we assume for contradiction that $\delta:=\inf_{\w\in\W}(\ma_{\lb_-}(\w)-\mr_{\lb_-}(\w))>0$.
Let us check that $\bar x\ge\mr_{\lb_-}(\bar\w)+\delta$ whenever $(\bar\w,\bar x)$ is in the closure $\mK^a_{\lb_-}$ of
$\{(\w,\ma_{\lb_-}(\w))\,|\;\w\in\W\}$. We write $(\bar\w,\bar x)=\lim_{n\to\infty}(\w_n,\ma_{\lb_-}(\w_n))$
and assume without restriction the existence of $(\bar\w,x^*):=\lim_{n\to\infty}(\w_n,\mr_{\lb_-}(\w_n))$,
which belongs to the (closed) set $\mB_{\lb_-}$.
Then, $x^*\le\bar x-\delta<\mr_{\lb_-}(\bar\w)+\delta-\delta<\mr_{\lb_-}(\bar\w)$
if $\bar x< \mr_{\lb_-}(\bar\w)+\delta$, impossible.
As said in Section \ref{subsec:skew}, the upper Lyapunov exponent of
$\mK^a_{\lb_-}$ takes the form $\int_{\W}\mh_x(\w,\mb(\w),\lb_-)\,dm$ for a measure
$m\in\merg$ and an $m$-measurable $\tau$-equilibrium $\mb$ with graph in $\mK^a_{\lb_-}$;
hence $\mb>\mr_{\lb_-}$, which is also $m$-measurable due to its semicontinuity. Hence,
\cite[Theorem 3.2]{dno4} shows that all the Lyapunov exponents of $\mK^a_{\lb_-}$ are
strictly negative, and that the upper and lower equilibria of $\mK^a_{\lb_-}$
coincide on a $\sigma$-invariant set $\W_0$ with $m(\W_0)=1$ for all $m\in\merg$.
Hence, Theorem \ref{th:2copia} ensures that $\mK^a_{\lb_-}$ is an attractive hyperbolic copy of $\W$.
This fact and the previous property ensure that $\mK^a_{\lb_-}$
is strictly above the closure $\mK^r_{\lb_-}$ of
$\{(\w,\mr_{\lb_-}(\w))\,|\;\w\in\W\}$. Hence, \cite[Theorem 3.3]{dno4} ensures
that $\mK^r_{\lb_-}$ is a repulsive hyperbolic copy of $\W$.
That is, $\lb_-\in\mI$, which, as said before, is not the case. The proof of (ii) is complete.

Let us prove (iv) assuming for contradiction that there exists $\lb_1\in[\lb_*,\lb_-)$ such that
$\W^b_{\lb_1}=\W$, which according to Proposition \ref{prop:3Csolocoer} ensures the same for
all $\lb\ge\lb_1$. We will show that $\inf_{\w\in\W}(\ma_{\lb_-}(\w)-\ma_{\lb_1}(\w))>0$,
which combined with \eqref{des:3Cchain} shows that $\inf_{\w\in\W}(\ma_{\lb_-}(\w)-\mr_{\lb_-}(\w))>0$ and
contradicts (ii). Since $\ma_{\lb_-}>\ma_{\lb_1}$,
\[
\begin{split}
 \ma_{\lb_-}(\w{\cdot}1)-\ma_{\lb_1}(\w{\cdot}1)
 &=v_{\lb_-}(1,\w,\ma_{\lb_-}(\w))-v_{\lb_1}(1,\w,\ma_{\lb_1}(\w))\\
 &>v_{\lb_-}(1,\w,\ma_{\lb_1}(\w))-v_{\lb_1}(1,\w,\ma_{\lb_1}(\w))=y_\w(1)\,(\lb_--\lb_1)
\end{split}
\]
for
$y_\w(t):=(v_{\lb_-}(t,\w,\ma_{\lb_1}(\w))-v_{\lb_1}(t,\w,\ma_{\lb_1}(\w)))/(\lb_--\lb_1)$, which satisfies
$y_\w(0)=0$ and $y_\w(t)\ge 0$ for all $t>0$.  So, if suffices to check
that $y_\w(1)\ge k$, where $k>0$ is independent of $\w$. We take $j\in\N$ such that $\mB_{\lb_-}\subseteq\W\times[-j,j]$, so that
$\mB_{\lb}\subseteq\W\times[-j,j]$ for all $\lb\in[\lb_1,\lb_-]$. Hence,
$[v_{\lb_1}(t,\w,\ma_{\lb_1}(\w)),v_{\lb_-}(t,\w,\ma_{\lb_1}(\w))]\subseteq
[\ma_{\lb_1}(\wt),\ma_{\lb_-}(\wt)]\subseteq[-j,j]$
for all $t\ge 0$ and $\w\in\W$. By hypothesis,
\[
 \bar c_j:=\inf_{\w\in\W}\big(\mh(\w,v_{\lb_1}(1,\w,\ma_{\lb_1}(\w)),\lb_-)-
 \mh(\w,v_{\lb_1}(1,\w,\ma_{\lb_1}(\w)),\lb_1)\big)/(\lb_--\lb_1)>0\,.
\]
We call $d_j:=\inf_{(\w,x)\in\W\times[-j,j]}\mh_x(\w,x,\lb_-)$.
Then, if $t\ge 0$,
\[
\begin{split}
 (\lb_--\lb_1)\,(y_\w)'(t)&=\big(\mh(\wt,v_{\lb_-}(t,\w,\ma_{\lb_1}(\w)),\lb_-)-
 \mh(\wt,v_{\lb_1}(t,\w,\ma_{\lb_1}(\w)),\lb_-) \big)\\
 &\quad+\big(\mh(\wt,v_{\lb_1}(t,\w,\ma_{\lb_1}(\w)),\lb_-)
 -\mh(\wt,v_{\lb_1}(t,\w,\ma_{\lb_1}(\w)),\lb_1)\big)\,,
\end{split}
\]
from where it follows easily that $(y_\w)'(t)\ge d_j\,y_\w(t)+\bar c_j$.
Therefore, $y_\w(1)\ge \bar c_j(e^{d_j}-1)/d_j>0$, as asserted. This completes the proof of (iv).
\smallskip

(v)-(vi) Assertion (v) follows from (iv): we assume that $\W$ is minimal and, for contradiction,
that $\lb_*<\lb_-$. Then, the set $\W_{\lb}^b$ is a proper compact $\tau_\lb$-invariant set
of $\W$ for all $\lb\in(\lb_*,\lb_-)$, which contradicts the minimality.

Let us prove (vi). It follows from (i) and (ii) that, for any $\w\in\W$,
\eqref{eq:3ini}$_{\w}^\lb$ has (at least) a bounded solution if $\lb\ge\lb_-^\w$.
This property ensures that $\lb_*\le\lb_-^\w$,
and hence $\lb_*\le\inf_{\w\in\W}\lb_-^\w$. On the other hand, if $\lb\ge\lb_-$, then
\eqref{eq:3ini}$_{\w}^\lb$ has bounded solutions for any $\w\in\W_\w$,
so (iii) ensures that $\lb\ge\lb_-^\w$. Hence, $\lb_-\ge\sup_{\w\in\W}\lb_-^\w$.
To check the converse inequalities, we first assume for contradiction the existence of
$\lb\in(\lb_*,\inf_{\w\in\W}\lb_-^\w)$. Then, by (iv), $\W^b_\lb$ is a proper subset of $\W$.
Let $\mM_\lb$ be a minimal subset of $\W^b_\lb$, and let us take $\w_0\in\mM_\lb$. Then, there are
bounded solutions for \eqref{eq:3Chlb}$_{\w_0}^\lb$, which according to (v) ensures that
$\lb_-^{\w_0}\le\lb$, impossible. Finally, also for contradiction, we assume the existence of
$\lb\in(\sup_{\w\in\W}\lb_-^\w,\lb_-)$. Then, for any $\w\in\W$, there exists at least a
bounded solution $v_\lb(t,\w,x)$. This means that $\W_\lb^b=\W$, which contradicts (iv).
\end{proof}
So, this theorem describes a {\em generalized nonautonomous global saddle-node
bifurcation diagram\/}
for $x'=\mh(\wt,x,\lb)$: for large values of $\lb$, there exists an attractor-repeller
pair of hyperbolic copies of the base which are as large and separated from each other
as desired; as $\lb$ decreases they approach each other monotonically; there is a {\em bifurcation\/}
value $\lb_-$ of $\lb$ at which still there are bounded solutions for all the equations of the
family but the upper and lower bounded equilibria are no longer uniformly separated; and to the
left of this bifurcation value, either there are no bounded solutions, or they exist for part
(not all) of the equations of the family and are nonhyperbolic, until they finally disappear
to the left of a second special value of the parameter. The first one of these final two options
means a {\em classic nonautonomous saddle-node bifurcation diagram}, and it occurs if
$\W$ is minimal.

Usually, in applications, one has a single equation giving rise
to a process instead of a family of equations giving rise to a flow.
The remaining results of this section explain how to apply part of the information so far obtained
to this single case, in which we obtain a classic global saddle-node nonautonomous
bifurcation diagram (with possibly highly complex dynamics at the unique bifurcation point:
see, e.g.,~\cite{jnot,nuob6}). To our knowledge, any proof of the stated properties requires
the hull construction, which is the reason for which this is the right point
to state and prove these properties. Recall that $\W_\w\subseteq\W$ is the closure of $\{\wt\,|\;t\in\R\}$.
\begin{teor}\label{th:3Cbifuruna}
Assume all the hypotheses of Theorem {\rm\ref{th:3Cbifur}}.
Let us fix $\bar\w\in\W$, and let $r_\lb$ and $a_\lb$ represent the lower and upper
bounded solutions of $x'=\mh(\bwt,x,\lb)$ in the case of existence. Then,
there exists $\bar\lb\in\R$ such that
$x'=\mh(\bwt,x,\lb)$ has:
\begin{itemize}[leftmargin=20pt]
\item[-] two uniformly separated hyperbolic solutions
$\tilde r_\lb:=r_\lb$ and $\tilde a_\lb=:a_\lb$ for $\lb>\bar\lb$, with
$\tilde r_\lb<\tilde a_\lb$, $\tilde r_\lb$ repulsive and
$\tilde a_\lb$ attractive;
\item[-] bounded but neither hyperbolic solutions nor uniformly separated solutions
for $\lb=\bar\lb$;
\item[-] and no bounded solutions for $\lb<\bar\lb$;
\end{itemize}
In addition, $(\bar\lb,\infty)\to C(\R,\R)\,,\;\lb\mapsto-\tilde r_\lb,\,\tilde a_\lb$
are continuous for the uniform topology and strictly increasing
on $(\bar\lb,\infty)$; $a_{\bar\lb}=\lim_{\lb\to(\bar\lb)^+}\tilde a_\lb$ and
$r_{\bar\lb}=\lim_{\lb\to(\bar\lb)^+}\tilde r_\lb$ pointwise on $\R$;
and $\lim_{\lb\to\infty}\tilde a_\lb=
-\lim_{\lb\to\infty}\tilde r_\lb=\infty$ uniformly on $\R$.
And $\bar\lb$ is common for all $\w\in\W$ such that $\W_\w=\W_{\bar\w}$.
\end{teor}
\begin{proof}
Let $\bar\lb:=\lb_-^{\bar\w}$ be the upper special value of
the parameter associated by Theorem \ref{th:3Cbifur}
to the family $x'=\mh(\wt,x,\lb)$ for $\w\in\W_{\bar\w}$,
where $\W_{\bar\w}$ is the closure of $\{\bwt\,|\;t\in\R\}$.
We will check that this is the value of the statement,
which in turn implies that it is common for all the elements $\w$ of $\W$ with
$\W_\w=\W_{\bar\w}$.

First, let $\tmr_\lb,\tma_\lb\colon\W_{\bar\w}\to\R$ be the maps provided by
Theorem \ref{th:3Cbifur}(i) restricted to $\W_{\bar\w}$.
Proposition \ref{prop:2extiende} shows that, if $\lb>\bar\lb$, then $t\mapsto\tmr_\lb(\bwt)$ and
$t\mapsto\tma_\lb(\bwt)$ are the two hyperbolic solutions described
in the statement. The stated continuity,
monotonicity and limiting behaviour as $\lb\to(\bar\lb)^+$ and as $\lb\to\infty$
follow from Theorem \ref{th:3Cbifur}(i),(ii).

Second, we take $\lb<\bar\lb$ and assume for contradiction the
existence of a bounded solution $b(t)=v_\lb(t,\bar\w,b(0))$ of
$x'=\mh(\bwt,x,\lb)$.
Given any $\w\in\W_{\bar\w}$, we write $\w=\lim_{n\to\infty}\bwt_n$ for a
sequence $(t_n)$, assume without restriction the existence of $b_0:=\lim_{n\to\infty}
b(t_n)$, and deduce that $v_\lb(t,\w,b_0)=\lim_{n\to\infty}v_\lb(t,\bwt_n,b(t_n))=
\lim_{n\to\infty}v_\lb(t,\bwt_n,v_\lb(t_n,\bar\w,b(0)))=
\lim_{n\to\infty}v_\lb(t+t_n,\bar\w,b(0))=\lim_{n\to\infty}b(t+t_n)$. So,
$v_\lb(t,\w,b_0)$ is a bounded solution of $x'=\mh(\wt,x,\lb)$, which means that
$\W_{\lb_{\bar\w}}^b=\W_{\bar\w}$ and hence contradicts
Theorem \ref{th:3Cbifur}(iii),(iv).

And third, Theorem \ref{th:3Cbifur}(ii) ensures the existence of bounded solutions for
$x'=\mh(\bwt,x,\bar\lb)$.
If one of them were hyperbolic, there would also be a hyperbolic solution
for $\lb<\bar\lb$ close enough (see Theorem \ref{th:2persistencia}), which,
as just seen, is not the case. Now, we assume for contradiction that there exist
two uniformly separated solutions for $\bar\lb$.
Then, there exists $\delta>0$ such that the upper and lower equilibria of
$\mB_{\bar\lb}$ satisfy $\ma_{\bar\lb}(\bwt)-\mr_{\bar\lb}(\bwt)\ge\delta$ for all $t\in\R$,
which combined with the semicontinuity properties of $\ma_{\bar\lb}$ and $\mr_{\bar\lb}$
before mentioned ensures that $\ma_{\bar\lb}(\w)-\mr_{\bar\lb}(\w)\ge\delta$
for all $\w\in\W_{\bar\w}$. This contradicts Theorem \ref{th:3Cbifur}(ii).
\end{proof}

\begin{teor}\label{th:3Casint}
Assume the conditions of Theorem {\rm\ref{th:3Cbifuruna}} and let
$x_\lb(t,s,x)$ be the maximal solution of $x'=\mh(\bwt,x,\lb)$ with $x_\lb(s,s,x)=x$,
defined on $(\alpha_{s,x,\lb},\beta_{s,x,\lb})$. With the notation established in
Theorem {\rm\ref{th:3Cbifuruna}}: if $\lb\ge\bar\lb$, then
\begin{itemize}[leftmargin=20pt]
\item[-] $\lim_{t\to(\alpha_{s,x,\lb})^+}x_\lb(t,s,x)=\infty$
if and only if $x>a_\lb(s)$,
\item[-] and $\lim_{t\to(\beta_{s,x,\lb})^-}x_\lb(t,s,x)=-\infty$
if and only if $x<r_\lb(s)$;
\end{itemize}
if $\lb>\bar\lb$, then
\begin{itemize}[leftmargin=20pt]
\item[-] $\lim_{t\to-\infty}|x_\lb(t,s,x)-\tilde r_\lb(t)|=0$
if and only if $x<\tilde a_\lb(s)$,
\item[-] $\lim_{t\to\infty}|x_\lb(t,s,x)-\tilde a_\lb(t)|=0$ if and only if
$x>\tilde r_\lb(s)$;
\end{itemize}
and, if $\lb<\bar\lb$ and $x\in\R$, then
\begin{itemize}[leftmargin=20pt]
\item[-] $\lim_{t\to(\alpha_{s,x,\lb})^+}x_\lb(t,s,x)=\infty$ or
$\lim_{t\to(\beta_{s,x,\lb})^-}x_\lb(t,s,x)=-\infty$ (perhaps both).
\end{itemize}
\end{teor}
\begin{proof}
If $\lb\ge\bar\lb$, \cite[Proposition 3.5(vi)]{dno4} shows that,
if $x>a_\lb(s)$, then there exists $t_0\in(s,\beta_{s,x,\lb})$
such that $x_\lb(t_0,s,x)>\rho$, where $\rho$ satisfies $\mh(\w,x,\lb)<\delta<0$ for all $\w\in\W$ and
$x\le\rho$ (whose existence follows from \ref{c2}), and it is easy to deduce that
$\lim_{t\to(\alpha_{s,x,\lb})^+}x_\lb(t,s,x)=\infty$. A similar argument shows that
$\lim_{t\to(\beta_{s,x,\lb})^-}x_\lb(t,s,x)=-\infty$
if and only if $x<\tilde r_\lb(s)$. If $\lb<\bar\lb$, any solution is
unbounded as time decreases or increases (or both), and the same arguments work. The remaining
assertions for $\lb>\bar\lb$ are proved in \cite[Proposition 3.5(vii)]{dno4}.
\end{proof}

\begin{nota}\label{rm:3Cbifclasica}
Let us fix $\mh^0\colon\WR\to\R$ and define $\mh(\w,x,\lb):=\mh^0(\w,x)+\lb\,\mg(\w)$
for a continuous map $\mg\colon\W\to\R$
with $\inf_{\w\in\W}\mg(\w)>0$ and $\lb\in\R$, and let us call $\mh^\lb(\w,x):=\mh(\w,x,\lb)$.
Property \ref{c1} for $\mh^0$ is inherited by $\mh^\lb$ for all $\lb\in\R$.
If we replace \ref{c2} for $\mh^0$ by a stronger
condition, as $\lim_{x\to\pm\infty}\mh^0(\w,x)=-\infty$ uniformly on $\W$, then
$\mh^\lb$ satisfies \ref{c2} for all $\lb\in\R$.
So, these conditions ensure that
all the conclusions of Proposition \ref{prop:3Csolocoer} hold. In addition, $\mh^\lb$
also inherits \ref{c3} and \ref{c4} from $\mh$. Therefore, if
$\mh^0$ satisfies \ref{c1}, \ref{c3}, \ref{c4} and
the mentioned stronger version of \ref{c2}, then $\mh$
satisfies all the hypotheses of Theorems \ref{th:3Cbifur} and \ref{th:3Cbifuruna}.
\end{nota}
Let us finally establish conditions on $h(t,x)$ and $g(t)$
giving rise to a global saddle-node nonautonomous bifurcation diagram for the parametric family
$x'=h(t,x)+\lb\,g(t)$: for large values of $\lb$, the upper and lower bounded solutions form
an attractor-repeller pair of hyperbolic solutions which are as large and separated
from each other as desired; as $\lb$ decreases they approach each other monotonically;
there is a {\em bifurcation value\/} of $\lb$ at which they are no longer separated
and for which any possible bounded solution is nonhyperbolic; and there are no bounded
solutions to the left of this bifurcation value. Theorem \ref{th:3Cbifurproceso} provides
an extension of \cite[Theorem 3.6]{lnor} to a much more general type of equation.
\begin{teor}\label{th:3Cbifurproceso}
Let $h\colon\RR\to\R$ and $g\colon\R\to\R$ satisfy: the functions $h$ and $h_x$ are admissible;
$\lim_{x\to\pm\infty}h(t,x)=-\infty$ uniformly on $\R$;
for any $x_1<x_2$, there exists $\delta_{x_1,x_2}>0$ such that
$\liminf_{t\to\pm\infty}(h_x(t,x_1)-h_x(t,x_2))\ge\delta_{x_1,x_2}$; and $g$ is uniformly
continuous, bounded, and with $\inf_{t\in\R} g(t)>0$.
Then, there exists $\bar\lb\in\R$ such that all the conclusions of
Theorems {\rm\ref{th:3Cbifuruna}} and {\rm\ref{th:3Casint}} hold for $x'=h(t,x)+\lb\,g(t)$.
\end{teor}
\begin{proof}
We define $\W_{h,g}$ as the closure of $\{(h_s,g_s)\,|\;s\in\R\}$ in the compact-open
topology of $C(\RR,\R)\times C(\R,\R)$, where $h_s(t,x)=h(s+t,x)$ and $g_s(t)=g(s+t)$.
We also define $\mh(\w,x)=\w_1(0,x)$ and $\mg(\w)=\w_2(0)$ if $\w=(\w_1,\w_2)$.
Standard arguments (see \cite[Theorem~I.3.1]{shyi4} and
\cite[Theorem IV.3]{selltopdyn}) prove that $\W_{h,g}$ is a compact metric space, that
the time shift defines a continuous transitive flow, and that $\mh^\lb:=\mh+\lb\,\mg$
satisfies \ref{c1}.
If, in addition,
$\lim_{x\to\pm\infty} h(t,x)=-\infty$ uniformly on $\R$, then, for each
$\lb\in\R$, there exists $\rho_\lb>0$ such
that $h(t,x)+\lb g(t)\le-1$ if $|x|\ge\rho_\lb$ and $t\in\R$.
Since any $\w=(\w_1,\w_2)\in\W_{h,g}$ satisfies
$(\w_1(0,x),\w_2(0))=\lim_{n\to\infty}(h(t_n,x),g(t_n))$ for a sequence $(t_n)$, we have
$\mh(\w,x)+\lb\mg(\w)=\w_1(0,x)+\lb\,\w_2(0)\le-1$ if $|x|\ge\rho_\lb$,
and this ensures \ref{c2} for all the functions $\mh^\lb$.

The arguments of \cite[Lemma 4.3]{dno4} also prove that any ergodic measure is
concentrated in the union of the \upalfa-limit and \upomeg-limit sets of
the (dense) orbit of $(h,g)$ in $\W_{h,g}$. So, \ref{c3} and
\ref{c4} follow from the fact that the map $(\w,x)\mapsto(\mh^\lb)_x(\w,x)=\mh_x(\w,x)$
is strictly decreasing on $\R$ for each $\w$ in these limit sets.
We take $\w=(\w_1,\w_2)$ in the \upomeg-limit set, write $\w_1(0,x)=\lim_{n\to\infty}h(t_n,x)$
uniformly on the compact subsets of $\R$ for a suitable sequence $(t_n)\uparrow\infty$,
assume without restriction that $(\w_1)_x(0,x)=\lim_{n\to\infty}h_x(t_n,x)$ uniformly on the
compact subsets of $\R$, and deduce from the hypothesis on $\liminf_{t\to\infty}$ that
$\mh_x(\w,x_1)-\mh_x(\w,x_2)=(\w_1)_x(0,x_1)-(\w_1)_x(0,x_2)=
\lim_{n\to\infty}(h_x(t_n,x_1)-h_x(t_n,x_2))\ge\delta_{x_1,x_2}>0$
if $x_1<x_2$, which proves the assertion in this case. The other one is analogous.

So, the maps $(\w,x,\lb)\mapsto\mh(\w,x)+\lb\,\mg(\w)$ satisfy all the conditions of
Theorems \ref{th:3Cbifur} and \ref{th:3Cbifuruna}, from where the assertion follows.
\end{proof}
Observe that $h$ is not required to be concave in Theorem \ref{th:3Cbifurproceso}: for instance,
$h(t,x)=-x^2+x^3\,e^{-t^2}$ satisfies all the assumed conditions on $h$ in its statement. This fact
considerably widens the possibilities of application of the result with respect to the
previous approaches of \cite{lnor} and \cite{dlo1}.
\subsection{A Devaney-chaotic flow}\label{subsec:31}
As explained in the Introduction, the very general framework of the results of Sections
\ref{sec:3} and \ref{sec:4}
allows us to use the skewproduct techniques in order to analyze dynamics and bifurcation diagrams for ``almost
stochastic" families of equations. This subsection describes the ``very large" set of continuous
functions which serves to illustrate the previous assertion.
\par
Let us fix $k_1>0$ and $k_2>0$ and consider the set
\begin{equation}\label{def:3Deva}
 \mP:=\{p\in C(\R,\R)\,|\;\|p\|_\infty\leq k_1\,,\;\mathrm{Lip}(p)\leq k_2\}\,,
\end{equation}
where $\|p\|_\infty:=\sup_{t\in\R}|p(t)|$ and $\mathrm{Lip}(p):=\sup_{t\neq s}|p(t)-p(s)|/|t-s|$, and
endow it with the compact-open topology. It is easy to check that the space $\mP$ is completely
metrizable, with
\[
 d(p,q)=\sum_{n=1}^\infty\frac{1}{2^n}\frac{\sup_{t\in[-n,n]}|p(t)-q(t)|}{1+\sup_{t\in[-n,n]}|p(t)-q(t)|}\,.
\]
So, $\mP$ is a metric space. Let us check that it is compact, and hence separable.
We take a sequence $(p_n)$ in $\mP$ and call $p_n^m$ to the restriction of $p_n$ to $[-m,m]$
for each $m\in\N$. Since the sequence $(p_n^m)$ is bounded and equicontinuous, Arzel\`{a}-Ascoli's
theorem provides a subsequence $(p^m_{n^m_k})$ which converges as $k\to\infty$
to a limit $p^m$ uniformly on $[-m,m]$. Now we follow a Cantor diagonal process: we choose $(p^{m+1}_{n^{m+1}_k})$ as a subsequence of $(p^m_{n^m_k})$. Then, the diagonal subsequence
$(p^m_{n_m^m})$ converges as $m\to\infty$ to a function $p\in C(\R,\R)$ uniformly at
each interval $[-m,m]$, and hence in the compact-open topology of $C(\R,\R)$.
The pointwise convergence suffices to check that $p\in\mP$.
This shows the assertion.

It is easy to check that the shift $\sigma_\mP\colon\R\times\mP\to\mP$, $(t,p)\mapsto \pt$, where $\pt(s)=p(t+s)$ for all $s\in\R$, defines a continuous flow on $\mP$. The next results proves, among
other properties, the transitivity of this flow. More precisely, the existence of a dense forward semiorbit. Recall that, on a compact metric space, this property is equivalent to {\em topological transitivity\/} (see e.g. \cite[Lemma 3]{auyo}), which is the property usually required in the definition of Devaney's chaos.
\begin{teor}\label{th:3Deva}
The continuous flow $(\mP,\sigma_\mP)$ is chaotic in the sense of Devaney:
\begin{itemize}
 \item[{\rm (i)}] there exists a dense forward semiorbit;
 more precisely, $\mP$ is the \upalfa-limit set and the \upomeg-limit set of a residual subset of points;
 \item[{\rm (ii)}] the set of points with periodic orbit is dense in $\mP$; and
 \item[{\rm (iii)}] the flow has sensitive dependence on initial conditions, that is,
 there exists $\delta>0$ such that, for any $p\in\mP$ and $\ep>0$, there exists $p_0\in B_\mP(p,\ep)$ and $t_0>0$ such that $d(\pt_0,p_0{\cdot}t_0)>\delta$.
 \end{itemize}
\end{teor}
\begin{proof}
(i) The separability of $\mP$ ensures that there exists a countable dense
subset $\{p_m\,|\;m=1,2,\ldots\}$ of $\mP$. Our goal is to construct a function $q\in\mP$ with dense forward semiorbit orbit.
Let us call $k:=2k_1/k_2$. We consider the disjoint closed intervals $\mI_j=[a_j,b_j]$ for $j=1,2,\ldots$,
defined as follows: $a_{j+1}-b_j=k$ for all $j\ge 0$; $\mI_0=[0,2]$ is the unique interval of length $2$;
the two next two intervals, $\mI_1$ and $\mI_2$ have length $4$; the next three intervals,
$\mI_3,\,\mI_4$ and $\mI_5$, have length $6$; and so on, so that the $j$ intervals
$\mI_{(j-1)\,j/2+1}$, $\mI_{(j-1)\,j/2+2}$, \ldots,\,
$\mI_{(j-1)\,j/2+j}=\mI_{j\,(j+1)/2}$ have length $2j$.
In each one of these intervals, the map $q$ is defined in terms of one of the maps $p_m$, with this schedule: of $p_1$ on $\mI_1$; $p_1$ on $\mI_2$ and $p_2$ on $\mI_3$; $p_1$ on $\mI_4$, $p_2$ on $\mI_5$ and
$p_3$ on $\mI_6$; and so on, so that $p_1$ on $\mI_{(j-1)\,j/2+1}$, $p_2$ on $\mI_{(j-1)\,j/2+2}$, \ldots and $p_j$ on $\mI_{(j-1)\,j/2+j}$. And, if $p_m$ is the function assigned to $\mI_j=[a_j,a_j+2i_j]$ for $i_j\in\N$, (``$p_m$ acts on $\mI_j$", in what follows), then $q(t)=p_m(t-a_j-i_j)$.
That is, the graph of $q$ on
$[a_j,a_j+2i_j]$ reproduces that of $p_m$ on $[-i_j,i_j]$.
The definition of $q$ on $[0,\infty)$ is completed by linear interpolation: for $t\in[b_j,a_{j+1}]=[b_j,b_j+k]$, $q(t)=q(b_j)+(q(a_{j+1})-q(b_j))\,(t-b_j)/k$.
And $q$ is extended to $\R$ as
an even map: $q(-t)=q(t)$ for $t\ge 0$. With this construction,
\begin{itemize}[leftmargin=15pt]
  \item[-] $q\in\mP$: it is continuous, with $\n{q}_\infty\le k_1$ and Lip$(q)\le k_2$ (since
    $|q'(t)|=|(q(a_{j+1})-q(b_j))/k|\le k_2$ for $t\in(b_j,a_{j+1})$).
  \item[-] For each $m=1,2,\ldots$, $p_m$ acts on intervals of lengths $2m,\,2m+2,2m+4,\ldots$ As a sample: the first six intervals are $[0,2]$ (of length $2$), $[2+k,6+k]$ and $[6+2k,10+2k]$ (of length $4$), and $[10+3k,16+3k]$, $[16+4k, 22+4k]$ and $[22+5k, 28+5k]$ (of length $6$); and on them the definition of $q$ is given in terms of $p_1,\,p_1$ and $\,p_2$, and $\,p_1,\,p_2$ and $p_3$, respectively.
\end{itemize}
An analytic construction of $q$ is, of course, possible. We omit the details, since they do not add
clarity to the explanation.

Let us check that $\{q{\cdot}s\,|\;s\ge 0\}$ is dense on $\mP$. We take $p\in\mP$, $\ep>0$,
$m_0\ge 0$ such that $d(p_{m_0},p)\le\ep/2$, and $i_0\ge m_0$ such that
$\sum_{n=i_0+1}^\infty 1/2^n\le\ep/2$.
Our goal is to find $s>0$ large enough to ensure that $d(q{\cdot}s,p_{m_0})\le\ep/2$, which, by the definition of the distance, is achieved if we get $q{\cdot}s(t)=p_{m_0}(t)$ for $t\in[-i_0,i_0]$.
We look for an interval $\mI_j$ of length $2i_j\ge 2i_0$ on which $p_{m_0}$ acts,
say $\mI_j=[a_j,b_j]=[a_j,a_j+2i_j]$. Then, for $t\in[-i_j,i_j]$,
$q(t+a_j+i_j)=p_m(t)$, which means that $s_{m_0}=a_j+i_j$ satisfies the requirement. Since
$d(q{\cdot}s_{m_0},p)\le d(q{\cdot}s_{m_0},p_{m_0})+d(p_{m_0},p)\le\ep$, the proof is complete.
We point out that, for each $m_0\in\N$, we can take a sequence $\mI_{j_n}$ of intervals with $(j_n)\uparrow\infty$ of lengths $(2i_n)\uparrow\infty$ on which $p_{m_0}$ acts, and hence we can take $s_{m_0}$ as large as desired.
This means that $\mP$ is, in fact, the \upomeg-limit set of $q$.
The same ideas allow us to check that $\mP$ is also the \upalfa-limit set of $q$.
By reasoning as in \cite[Proposition I.11.4]{mane1}, we conclude that this happens for a
residual subset of points of $\mP$.
\smallskip

(ii) Let $q\in\mP$ be the function with dense orbit described in (i) and, for $n\in\N$, let
$q_n\colon\R\to\R$ be the $2n$-periodic function which satisfies $q_n(t)=q(t)$ for all
$t\in[-n,n]$. It is clear that $q_n\to q$ in $\mP$ as $n\to\infty$. We take $p\in\mP$ and $\ep>0$,
look for $s_0\in\R$ such that $d(q{\cdot}s_0,p)\le\ep/2$, look for $n_0\in\N$ such that $d(q_{n_0}{\cdot}s_0,q{\cdot}s_0)\le\ep/2$, and conclude that $d(q_{n_0}{\cdot}s_0,p)\leq\ep$.
So, $\{q_n\,|\;n\in\N\}$ is dense in $\mP$, which proves (ii).
\smallskip

(iii) We define $\delta:=k_1/(4+2\,k_1)$. Let $p\in\mP$ and $\ep>0$ be fixed.
We take $j\in\N$ such that $\sum_{n=j+1}^\infty1/2^n<\ep/2$ and $t_0=-j-k_1/k_2$, and consider
the (unique) continuous function $\bar p$ such that: for $t\le t_0$, $\bar p(t)=p(t_0)+k_1$
if $p(t)< 0$ and $\bar p(t)=p(t_0)-k_1$ if $p(t_0)\ge 0$,
$\bar p(t)=p(t)$ for $t\ge-j$, and $\bar p$ is linear on $[t_0,-j]$. It is easy to check
that $\bar p\in\mP$, and that
\[
 d(\bar p,p)=\sum_{n=j+1}^\infty\frac{1}{2^n}\frac{\sup_{t\in[-n,n]}|\bar p(t)-p(t)|}
 {1+\sup_{t\in[-n,n]}|\bar p(t)-p(t)|}\le \sum_{n=j+1}^\infty\frac{1}{2^n}<\frac{\ep}{2}\,.
\]
Let $q\in\mP$ be the function with dense orbit described in (i). Hence, $\bar p$ is the limit of
$(q{\cdot}s_n)$ for a suitable sequence $(s_n)$. Since $\lim_{n\to\infty}d(q{\cdot}s_n,\bar p)=0$
and $\lim_{n\to\infty}(q{\cdot}s_n(t_0)-\bar p(t_0))=0$, there exists
$s\in\R$ such that $d(q{\cdot}s,\bar p)\le\ep/2$ and $|q{\cdot}s(t_0)-\bar p(t_0)|\le k_1/2$.
Therefore, $d(q{\cdot}s,p)\le d(q{\cdot}s,\bar p)+d(\bar p,p)<\ep$ and
$|q{\cdot}s(t_0)-p(t_0)|\ge |\bar p(t_0)-p(t_0)|-|\bar p(t_0)-q{\cdot}s(t_0)|=
k_1-|\bar p(t_0)-q{\cdot}s(t_0)|\ge k_1/2$.
Since $x\mapsto x/(1+x)$ is strictly increasing for $x>0$ and  $\sup_{t\in[-1,1]}|q{\cdot}s{\cdot}t_0(t)-\pt_0(t)|\ge |q{\cdot}s{\cdot}t_0(0)-\pt_0(0)|
=|q{\cdot}s{\cdot}(t_0)-p(t_0)|\ge k_1/2$,
\[
 d(q{\cdot}s{\cdot}t_0,\pt_0)\geq\frac{1}{2}\,
 \frac{\sup_{t\in[-1,1]}|q{\cdot}s{\cdot}t_0(t)-\pt_0(t)|}
 {1+\sup_{t\in[-1,1]}|q{\cdot}s{\cdot}t_0(t)-\pt_0(t)|}\ge \frac{1}{2}\,\frac{k_1/2}{1+k_1/2}=\delta\,.
\]
Hence, $p_0=q{\cdot}s$ satisfies the condition in the statement of (iii).
\end{proof}
\subsection{Applications to population dynamics}\label{subsec:32}
The main purpose of this section is to illustrate the scope of the theoretical results obtained in the concave
in measure case with the help of a nonautonomous single-species population model that includes
predation and emigration, one of whose coefficients moves in a large set of continuous functions.
Before describing it, in Section \ref{subsubsec:321}
we explain some properties of the bifurcation value $\tilde\lb(\theta,p)$
associated by Theorem \ref{th:3Cbifuruna} to each one of the equations of a certain family that varies
on the compact set given by the product $\Theta\times\mP$, where $\mP$ is one of the sets described
in Section \ref{subsec:31} and $\Theta$ is a compact metric space which, in the examples, will
be obtained as the hull of a certain map.
\par
So, let us consider a global real flow on a compact metric space, $(\Theta,\sigma_{\Theta})$, call
$\W:=\Theta\times\mP$, with $\mP$ defined by \eqref{def:3Deva}, and define $\sigma(t,\w):=(\sigma_{\Theta}(t,\theta),
\pt)$ for $\w=(\theta,p)$. Since there is no risk of confusion, we also represent $\theta{\cdot}t=\sigma_{\Theta}(t,\theta)$.
For each $\w=(\theta,p)\in\W$, we consider the parametric family of equations
\begin{equation}\label{eq:3Deva-2}
 x'=\mh(\theta{\cdot}t,x,\lb)+\mpp(\pt)\,.
\end{equation}
Here, $\mpp\colon\mP\to\R$ is the continuous map defined by $\mpp(p):=p(0)$,
and $\wit\mh(\w,x,\lb)=\mh(\theta,x,\lb)+\mpp(p)$ is assumed to satisfy the conditions of
Theorem \ref{th:3Cbifur}. Note that all these conditions can be rewritten in
terms of conditions on $\mh$.  Theorem \ref{th:3Cbifuruna} provides a unique bifurcation value $\tilde\lb(\theta,p)$ for the
parametric family \eqref{eq:3Deva-2}. Recall that \eqref{eq:3Deva-2} has:
no bounded solutions for $\lb<\tilde\lb(\theta,p)$, bounded but not hyperbolic solutions for
$\lb=\tilde\lb(\theta,p)$, and two uniformly separated hyperbolic solutions for $\lb>\tilde\lb(\theta,p)$.
Our next goal is to analyze the properties of the map
\begin{equation}\label{def:3lb-2}
 \tilde\lb\colon\W\to\R\,,\;(\theta,p)\mapsto\tilde\lb(\theta,p)\,.
\end{equation}
\par
Recall that the global flow $(\Theta,\sigma_{\Theta})$ is {\em almost-periodic \/} if,
for all $\ep>0$, there exists $\delta_\ep>0$ such that,
if $\text{\rm dist}_{\Theta}(\theta_1,\theta_2)\le\delta_\ep$, then
$\text{\rm dist}_{\Theta}(\theta_1{\cdot}t,\theta_2{\cdot}t)\le\ep$ for all $t\in\R$. Or, equivalently,
if for all $\ep>0$ there exists a relatively dense set $\mQ_\ep\subseteq\R$ such that
$\text{\rm dist}_{\Theta}(\sigma_{\Theta}(t,\theta),\theta)\le\ep$ for all $\theta\in\Theta$ and $t\in\mQ_\ep$: see
\cite[Theorem I.2.9]{shyi4} and \cite[Theorem IV.2.8]{vrie}. Of course, $\text{\rm dist}_{\Theta}$ represents
the distance in $\Theta$. This almost-periodicity condition will be assumed in the second part of the next result.
To simplify the notation, we will fix $k_1=1$ in the definition \eqref{def:3Deva} of the set $\mP$.
The result is identical for any positive constant.
\begin{prop}\label{prop:3lambda}
Let $\mP$ and $\tilde\lb$ be defined by \eqref{def:3Deva} with $k_1=1$ and \eqref{def:3lb-2}. Then,
\begin{itemize}
\item[\rm(i)] $\tilde\lb(\theta,p_1)\ge\tilde\lb(\theta,p_2)$ if $p_1\le p_2$, and $\tilde\lb(\theta,p_1)>\tilde\lb(\theta,p_2)$ if $\inf_{t\in\R}(p_2(t)-p_1(t))>0$.
\item[\rm(ii)] $\tilde\lb(\theta,1)\le\tilde\lb(\theta,p)\le\tilde\lb(\theta,-1)$ for all $(\theta,p)\in\W$.
\item[\rm(iii)] The map $\tilde\lb$ is $\sigma$-invariant and lower semicontinuous,
and the map $p\mapsto\tilde\lb(\theta,p)$ is continuous at $p\in\mP$ (at least) if $\tilde\lb(\theta,p)=\tilde\lb(\theta,-1)$.
\item[\rm(iv)] If $p=\lim_{n\to\infty} p_n$ uniformly on $\R$,
    then $\tilde\lb(\theta,p)=\lim_{n\to\infty}\tilde\lb(\theta,p_n)$ for all $\theta\in\W$.
\end{itemize}
Assume that, in addition, the flow $(\Theta,\sigma_{\Theta})$ is almost periodic. Then,
\begin{itemize}
\item[\rm(v)] if $p=\lim_{n\to\infty} p_n$ uniformly on $\R$ and $\theta=\lim_{n\to\infty}\theta_n$ in $\Theta$,
    then $\tilde\lb(\theta,p)=\lim_{n\to\infty}\tilde\lb(\theta_n,p_n)$ for all $\theta\in\Theta$.
\item[\rm(vi)] For any $\theta\in\Theta$, $\mR_{\theta}:=\{p\in\mP\,|\;\tilde\lb(\theta,p)=\tilde\lb(\theta,-1)\}$ is
the residual set of continuity points of the map $p\mapsto\tilde\lb(\theta,p)$ and contains the dense subset $\mD\subset\mP$ of
points with dense $\sigma_\mP$-orbit.
\item[\rm(vii)] The set $\mP\setminus\mR_\theta$ is dense for all $\theta\in\Theta$.
\item[\rm(viii)] If $(\theta,p)\in\W$, then $[\tilde\lb(\theta,p),\tilde\lb(\theta,-1)]
\subseteq\tilde\lb(\theta,B_\mP(p,\delta))$ for all $\delta>0$.
\end{itemize}
\end{prop}
\begin{proof}
(i) We take a bounded solution $b\colon\R\to\R$ of $x'=\mh(\theta{\cdot}t,x,\tilde\lb(\theta,p_1))+p_1(t)$,
with $\sup_{t\in\R}|b(t)|=c$, and observe that
$b'(t)\le \mh(\theta{\cdot}t,b(t),\tilde\lb(\theta,p_1))+p_2(t)-\delta$ if $p_1(t)\le p_2(t)-\delta$ for a $\delta\ge 0$ and
all $t\in\R$. If $\delta=0$, this property and \cite[Proposition 3.5(v)]{dno4} yield
$\tilde\lb(\theta,p_2)\le\tilde\lb(\theta,p_1)$. If $\delta>0$, the compactness of $\Theta\times[-j,j]$ and
the continuity of $\mh$ allow us to find $\ep>0$ such that
$b'(t)\le\mh(\theta{\cdot}t,b(t),\tilde\lb(\theta,p_1)-\ep)+p_2(t)$,
and \cite[Proposition 3.5(v)]{dno4} yields $\tilde\lb(\theta,p_2)\le\tilde\lb(\theta,p_1)-\ep<\tilde\lb(\theta,p_1)$.
\smallskip\par
(ii) These inequalities follow from the first assertion in (i) and from $-1\le p(t)\le 1$ for all $t\in\R$
if $p\in\mP$.
\smallskip\par
(iii) The definition of $\tilde\lb$ given in
Theorem \ref{th:3Cbifurproceso} shows that it is $\sigma$-invariant.
We take a sequence $(\theta_n,p_n)$ in $\W$ with limit $(\theta,p)$ and a bounded solution
$b_n$ of $x'=\mh(\theta_n{\cdot}t,x,\tilde\lb(\theta_n,p_n))+p_n(t)$ for each $n\in\N$. Assume the existence of
$\lb_0:=\lim_{m\to\infty}\tilde\lb(\theta_n,p_n)$ for a subsequence $(\theta_m,p_m)$.
The coercivity property of $\mh$ ensures the existence of $\rho>0$ such that $\sup_{t\in\R}|b_n(t)|\le\rho$:
this can be checked, for instance, as in \cite[Proposition 3.5(ii)]{dno4}. From here and
Arzel\`{a}-Ascoli's Theorem, we deduce the existence of a subsequence $(b_k)$ of $(b_m)$ which
converges uniformly on the compact intervals of $\R$ to a map $b_0$. Since
$b_k(t)=b_k(0)+\int_{0}^t (\mh(\theta_k{\cdot}s,b_k(s),\tilde\lb(\theta_k,p_k))+p_k(s))\,ds$ for any $t\in\R$, Lebesgue's Dominated
Convergence Theorem ensures that $b_0(t)=b_0(0)+\int_{0}^t (\mh(\theta{\cdot}s,b_0(s),\lb_0)+p(s))\,ds$.
This means that $b_0$ is a bounded solution of $x'=\mh(\theta{\cdot}t,b_0(t),\lb_0)+p(t)$, which in turn yields
$\tilde\lb(\theta,p)\le\lb_0$. This shows the lower semicontinuity of the map. Now we fix
$\theta\in\Theta$, assume that $\tilde\lb(\theta,p)=\tilde\lb(\theta,-1)$ for a map $p\in\mP$,
take a sequence $(p_n)$ in $\mP$ with limit $p$, and deduce from the semicontinuity and from (ii) that
$\tilde\lb(\theta,-1)=\tilde\lb(\theta,p)\le\liminf_{n\to\infty}\tilde\lb(\theta,p_n)\le
\limsup_{n\to\infty}\tilde\lb(\theta,p_n)\le\tilde\lb(\theta,-1)$,
which proves the second assertion in (iii).
\smallskip\par
(iv) The proof of this point reproduces that of (v) by taking $\theta_n=\theta$ for all $n\in\N$,
which makes it unnecessary to use the almost-periodicity of the flow on $\Theta$.
\smallskip\par
(v) The hypotheses in (v) and the semicontinuity property (iii) ensure that
$\tilde\lb(\theta,p)\le\liminf_{n\to\infty}\tilde\lb(\theta_n,p_n)$.
We take any $\ep>0$ and call $\lb_\ep=\tilde\lb(\theta,p)+\ep$. Then, $x'=\mh(\theta{\cdot}t,x,\lb_\ep)+p(t)$ has
two hyperbolic and uniformly separated solutions, with absolute values bounded by a constant $\rho>0$.
Since the maps $(\theta,x)\mapsto\mh(\theta,x,\lb_\ep),\,\mh_x(\theta,x,\lb_\ep)$ are uniformly
continuous on $\Theta\times[-\rho,\rho]$, given $\ep>0$ there exists $\delta_\ep>0$ such that
$\sup_{x\in[-\rho,\rho]}|\mf(\theta_1,x,\lb_\ep)-\mf(\theta_2,x,\lb_\ep)|<\ep$ for $\mf=\mh,\mh_x$ if
$\text{\rm dist}_{\Theta}(\theta_1,\theta_2)<\delta_\ep$. Since the flow on $\Theta$ is almost-periodic,
$\text{\rm dist}_{\Theta}(\theta{\cdot}t,\theta_n{\cdot}t)<\delta_\ep$ for all $t\in\R$ if $n$ is large enough.
Altogether, we have $\sup_{(t,x)\in\mathbb{R}\times[-\rho,\rho]}|\mf(\theta{\cdot}t,x,\lb_\ep)-\mf(\theta_n{\cdot}t,x,\lb_\ep)|<\ep$
for $\mf=\mh,\mh_x$ if $n$ is large enough. Recall also that, by hypothesis,
$\lim_{n\to\infty}\sup_{t\in\R}|p_n(t)-p(t)|=0$. In these conditions,
Theorem \ref{th:2persistencia} shows the existence of two hyperbolic solutions
for $x'=\mh(\theta_n{\cdot}t,x,\lb_\ep)+p_n(t)$ (which implies $\tilde\lb(\theta_n,p_n)\le\lb_\ep$)
if $n$ is large enough.
Altogether, we have $\tilde\lb(\theta,p)\le\liminf_{n\to\infty}\tilde\lb(\theta_n,p_n)\le
\limsup_{n\to\infty}\tilde\lb(\theta_n,p_n)\le \tilde\lb(\theta,p)+\ep$, which proves (v).
\smallskip\par
(vi) Let us fix $\theta\in\Theta$. We take $p\in\mD$ and a sequence $(t_n)$ such that $-1=\lim_{n\to\infty}\pt_n$ in $\mP$, and
fix $\ep>0$. Property (v) provides $\delta_\ep>0$ such that
$|\tilde\lb(\theta,p)-\tilde\lb(\theta_2,p)|<\ep$ if $\text{\rm dist}_{\Theta}(\theta_1,\theta_2)<\delta_\ep$. Note also that,
according to (iii), $\tilde\lb(\theta,\pt)=\tilde\lb(\theta{\cdot}(-t),p)$ for all $t\in\R$.
The almost periodicity of the flow in $\Theta$ provides a bounded sequence $(s_n)$ such that
$\text{\rm dist}_{\Theta}(\theta{\cdot}(-t_n+s_n),\theta)<\delta_\ep$, and it is easy to check that
$-1=\lim_{n\to\infty}p{\cdot}(t_n-s_n)$ in $\mP$. These properties combined with the lower semicontinuity proved in (iii) and the
second bound of (ii) yield $\tilde\lb(\theta,-1)\le\liminf_{n\to\infty}\tilde\lb(\theta,p{\cdot}(t_n-s_n))=
\liminf_{n\to\infty}\tilde\lb(\theta{\cdot}(-t_n+s_n),p)\le
\liminf_{n\to\infty}(\tilde\lb(\theta,p)+\ep)=\tilde\lb(\theta,p)+\ep\le\linebreak\tilde\lb(\theta,-1)+\ep$.
The arbitrary choice of $\ep$ shows that $\tilde\lb(\theta,p)=\tilde\lb(\theta,-1)$, which means that
$\mD\subseteq\mR_{\theta}$.
\par
Now, we take a continuity point $p$ of the map $p\mapsto\tilde\lb(\theta,p)$, and a sequence $(p_n)$ in the dense
set $\mD$ with $p=\lim_{n\to\infty}p_n$ in $\mP$. Then, $\tilde\lb(\theta,p)=\lim_{n\to\infty}
\tilde\lb(\theta,p_n)=\lim_{n\to\infty}\tilde\lb(\theta,-1)=\tilde\lb(\theta,-1)$,
and hence $p\in\mR_{\theta}$. Conversely, as seen in (iii), $p$ is a continuity point
if $\tilde\lb(\theta,p)=\tilde\lb(\theta,-1)$, which completes the proof of (vi).
\smallskip\par
(vii) We fix $\theta\in\Theta$, take $p\in\mP$ and define $p_n:=\min(p+1/n,1)\in\mP$ for $n\in\N$. Then,
$p_n\ge-1+1/n$, and so (i) yields $\tilde\lb(\theta,p_n)<\tilde\lb(\theta,-1)$, which ensures that
$p_n\in\mP\setminus\mR_{\theta}$: see (vi). And $\sup_{t\in\R}|p(t)-p_n(t)|\le 1/n$, so $(p_n)$ converges to $p$
in~$\mP$.\hspace{-1cm}~
\smallskip\par
(viii) We fix $(\theta,p)\in\W$ and $\delta>0$, take $q\in B_\mP(p,\delta)\cap\mD$,
observe that (vi) ensures that $\tilde\lb(\theta,q)=\tilde\lb(\theta,-1)$,
and define $p_s:=s\,q+(1-s)\,p$ for $s\in[0,1]$. Then $p_s$ belongs
to the convex set $B_\mP(p,\delta)$ (the ball centered in $p$ with radius $\delta$),
and property (iv) ensures the continuity of the map
$[0,1]\to\R,\,s\mapsto\tilde\lb(\theta,p_s)$. The range of this map contains
$[\tilde\lb(\theta,p_0),\tilde\lb(\theta,p_1)]=[\tilde\lb(\theta,p),\tilde\lb(\theta,q)]=
[\tilde\lb(\theta,p),\tilde\lb(\theta,-1)]\subseteq\tilde\lb(\theta,B_\mP(p,\delta))$. This proves (viii).
\end{proof}
\begin{nota}\label{rm:3casiper}
If the flow $(\Theta,\sigma_\theta)$ is uniformly almost-periodic and minimal, and if $k\in\R$,
then the maps $\theta\mapsto\tilde\lb(\theta,k)$ are constant: we fix $\theta_0$, write any other
$\theta$ as $\lim_{n\to\infty}\theta_0{\cdot}t_n$ in $\Theta$ for a sequence $(t_n)$,
note that $k{\cdot}t_n=k$, and deduce from (v) and (iii) of Proposition \ref{prop:3lambda}
that $\tilde\lb(\theta,k)=\lim_{n\to\infty}\tilde\lb(\theta_0{\cdot}t_n,k{\cdot}t_n)=
\lim_{n\to\infty}\tilde\lb(\theta_0,k)=\tilde\lb(\theta_0,k)$. Therefore, adding minimality to the
hypotheses of Proposition \ref{prop:3lambda} guarantees that the constant values $\tilde\lb(\cdot,1)$ and
$\tilde\lb(\cdot,-1)$ provide the values $\lb_*$ and $\lb_-$ of Theorem \ref{th:3Cbifur} for the family
given by \eqref{eq:3Deva-2} for $(\theta,p)\in\W$, as
we deduce from combining the definition of $\tilde\lb$, Proposition \ref{prop:3lambda}(ii) and
Theorem \ref{th:3Cbifur}(vi).
\end{nota}
\subsubsection{Numerical simulations}\label{subsubsec:321}
Let us describe the mentioned population model, for which we will perform two numerical simulations
intended to illustrate the variation in $\mP$ of the bifurcation function
$p\mapsto\tilde\lb(\theta,p)$ for two different sets $\Theta$.
The initial $\lb$-parameter bifurcation problem that we consider is
\begin{equation}\label{eq:concavemodel}
 x'=r(t)\,x\,\left(1-\frac{x}{K(t)}\right)+\lb\,\G(t)\,\frac{x^2}{b(t)+x^2}+p(t)+a\,,\quad p\in\mP.
\end{equation}
The quadratic polynomial $r(t)\,x\,(1-x/K(t))$ captures the internal dynamics of the species, where $r(t)$
stands for the intrinsic growth rate at time $t$, and $K(t)$, which is not a solution of the quadratic equation
unless it is constant, is the threshold over which the per capita growth rate at time $t$ is negative.
Both $r$ and $K$ are quasiperiodic and positively bounded from below.
The Holling's type III functional response term $\lb\,\G(t)\,x^2/(b(t)+x^2)$ captures the contribution
of predation, where the nonnegative number $-\lb\,\G(t)$ represents the number of predators at time $t$ and
$b(t)$ is proportional to the average time between attacks of a predator at time $t$.
We assume that $\G$ is continuous and negatively bounded from above, work with positive values of the
bifurcation parameter $\lb$, and assume that $b$ is quasiperiodic and positively bounded from below.
The map $p$ lies in the set $\mP$ given by \eqref{def:3Deva} with $k_1=k_2=1$, and
the sum $p(t)+a$, with $a\in\R$, describes the emigration at time $t$ and takes values in $[a-1,a+1]$.
We will assume that $a<-1$ to represent a negative net migration balance for every $t\in\R$.
The time variation of the aforementioned functions allows the influence of factors such as the Earth's rotation,
changes in related ecological communities, seasonal alternation or climatic variability to be incorporated
into the model.
\par
Let $\Theta$ be the joint hull of $r$, $K$, $\G$ and $b$,
$\Theta:=\mathrm{closure}_{C(\R,\R^4)}\{(r_s,K_s,\G_s,b_s)\mid\,s\in\R\}$ (with $f_s(t):=f(t+s)$),
and let $\W:=\Theta\times\mP$. So, if $\mf_i\colon\Theta\to\R$ is defined by
$\mf_i(x_1,x_2,x_3,x_4):=x_i(0)$, $\mpp(p):=p(0)$, and $\theta_0:=(r,K,\G,b)$, we include
the family \eqref{eq:concavemodel} (which varies in $\mP$) in the varying-in-$\W$ family
\begin{equation}\label{eq:concavemodelhull}
x'=\mf_1(\theta{\cdot}t)\,x\,\left(1-\frac{x}{\mf_2(\theta{\cdot}t)}\right)+
\frac{\lb\,\mf_3(\theta{\cdot}t)\,x^2}{\mf_4(\theta{\cdot}t)+x^2}+\mpp(p{\cdot}t)+a\,,\quad (\theta,p)\in\W\,:
\end{equation}
\eqref{eq:concavemodel} is \eqref{eq:concavemodelhull} for $\theta=\theta_0$.
Let us rewrite \eqref{eq:concavemodelhull} as $x'=\mh(\theta{\cdot}t,x,\lb)+\mpp(p{\cdot}t)+a$.
Since $\mh(\theta,x,\lb)<\mh(\theta,x,0)$ for all $\lb>0$, \cite[Proposition 3.5(iv)]{dno4}
ensures that the set of bounded orbits  for $\lb>0$ is strictly contained in that corresponding to $0$.
The inequality $\sup_{t\in\R}(p(t)+a)<0$ makes it easy to check the existence of constants
$\delta>0$ and $m_2>m_1>0$ such that $\mh(\w,\rho,0)+\mpp(p{\cdot}t)+a\le-\delta$ for any $\rho\notin(m_1,m_2)$
and all $(\theta,p)\in\W$,
so \cite[Proposition 3.5(ii)]{dno4} and the previous property
ensure that any bounded solution of \eqref{eq:concavemodelhull}$_\lb$ is positively bounded from below by $m_1$
for all $\lb\ge 0$.
In addition, the incremental quotient with respect to $\lb$ that appears in the hypotheses of Theorem \ref{th:3Cbifur},
and that here takes the form $(\theta,x)\mapsto \mf_3(\theta)\,x^2/(\mf_4(\theta)+x^2)$,
is strictly negative on $\R\times\mK$ for any compact subset $\mK\subset(0,\infty)$,
since $\mf_3$ is strictly negative and $\mf_4$ is strictly positive. It is clear
that there exists an interval $[0,\lb_0]$ such that
$(\theta,x)\mapsto\mh(\theta,x,\lb)$ satisfies \ref{c3} and \ref{c4} for all $\lb\in[0,\lb_0]$.
Note that, then, $((\theta,p),x)\mapsto\mh(\theta,x,\lb)+\mpp(p)+a$
satisfies \ref{c1}-\ref{c4} for all $\lb\in[0,\lb_0]$.
The extra conditions required in Theorems \ref{th:3Cbifur} and \ref{th:3Cbifuruna} regarding
global regularity and (more than) strict monotonicity with respect to $\lb$ are also fulfilled,
but now the map decreases with $\lb$, since the coefficient accompanying $\lb$ in
\eqref{eq:concavemodelhull} is strictly negative. Assume also that
\eqref{eq:concavemodel}$^0$ (for which the map is globally strictly concave in $x$)
has two hyperbolic solutions for $p=-1$ (as can be numerically checked in the particular examples
analyzed below). Then, \cite[Proposition 3.5(v)]{dno4} shows the
existence of two uniformly separated solutions for any $p\in\mP$, and hence \cite[Theorem 3.6]{dno4}
ensures the existence of two hyperbolic solutions for any $p\in\mP$. Assume also that,
for $p=1\in\mP$, the two hyperbolic solutions disappear as $\lambda$ increases from 0 at a first point
$\tilde\lb(\theta_0,1)$ which belongs to $(0,\lb_0)$ (as can also be numerically checked in our examples).
A contradiction argument with the same ideas as for $\lb=0$ shows that, for any $p\in\mP$,
the two hyperbolic solutions disappear at a first point $\tilde\lb(\theta_0,p)\in(0,\lb_0)$.
In these conditions, we can reason as in Theorems \ref{th:3Cbifur} and \ref{th:3Cbifuruna}
in order to check the absence of bounded solutions for \eqref{eq:concavemodel}$^\lb$ if
$\lb\in(\tilde\lb(\theta_0,p),\lb_0)$. We can also check the absence of bounded solutions for any $\lb\ge\lb_0$:
we assume for contradiction that $b(t)$ is a bounded solution for $\lb_1\ge\lb_0$, take
$\lb_2\in(\tilde\lb(\theta_0,p),\lb_0)$, and apply again \cite[Proposition 3.5(v)]{dno4} to check the
existence of a bounded solution for $\lb_2$, impossible.
Altogether, these conditions ensure the existence of a unique positive bifurcation
value $\tilde\lb(\theta_0,p)$ for each $p\in\mP$, although now the different type of monotonicity in
$\lb$ means the existence of hyperbolic solutions for $0\le\lb<\tilde\lb(\theta_0,p)$ and the
lack of bounded solutions for $\lb>\tilde\lb(\theta_0,p)$.
Consequently, the inequalities of Proposition \ref{prop:3lambda}(ii) will read $\tilde
\lb(\theta_0,-1)\le\tilde\lb(\theta_0,p)\le\tilde\lb(\theta_0,1)$ for all $p\in\mP$.
\begin{figure}[h]
\centering
\includegraphics[width=\textwidth]{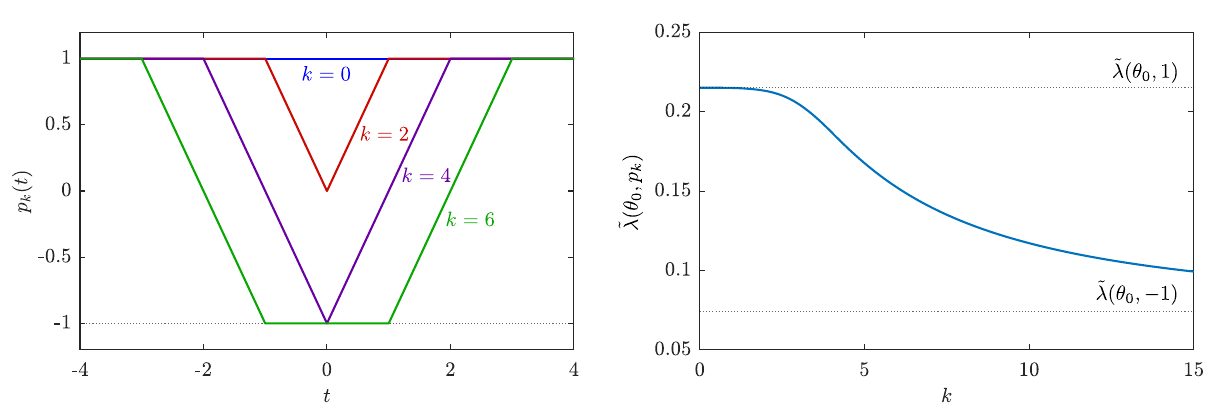}
\caption{In the left panel, $p_k(t):=\max\{-1,\min\{1,1-k/2-t\},\min\{1,1-k/2+t\}\}$
for different choices of $k$ in different colors. They belong to
the set $\mP$ given by \eqref{def:3Deva} with $k_1=1$ and $k_2=1$.
In the right panel, a numerical approximation of the bifurcation curve $\tilde\lb(\theta_0,p_k)$
for \eqref{eq:concavemodelhull}, with $r(t):=0.95+0.02\sin(2\pi t)+0.02\sin(t/10)$,
$K(t):=30+3.2\sin(2\pi t)+\sin(6\pi t)+0.4\cos(5\sqrt{5}\,t)+0.4\cos(40t)$,
$b(t):=400+26\cos^2(40t)$, $a=-3$ and $\G(t):=-160-16\sin(2\pi t)$.}
\label{fig:1}
\end{figure}
\vspace{.2cm}
\noindent\underline{Bifurcations in a quasiperiodic example}.
We take $r(t):=0.95+0.02\sin(2\pi t)+0.02\sin(t/10)$,
$K(t):=30+3.2\sin(2\pi t)+\sin(6\pi t)+0.4\cos(5\sqrt{5}\,t)+0.4\cos(40t)$,
$\G(t):=-160-16\sin(2\pi t)$, $b(t):=400+26\cos^2(40t)$ and $a:=-3$.
So, since $(t,x)\mapsto r(t)\,x\,\big(1-x/K(t)\big)+\lb\,\G(t)\,x^2/(b(t)+x^2)$ is a uniformly almost-periodic map,
the flow $(\Theta,\sigma_\Theta)$ is almost periodic and minimal (see e.g.~\cite{shyi4}),
and hence $p\mapsto\tilde\lb(\theta_0,p)$ also satisfies the last properties of Proposition \ref{prop:3lambda},
with interval $[\tilde\lb(\theta_0,-1),\tilde\lb(\theta_0,p)]$ in (viii).
In particular, Proposition \ref{prop:3lambda}(viii) shows that any
$\lb\in[\tilde\lb(\theta_0,-1),\tilde\lb(\theta_0,1)]\subset(0,\lb_0)$
is the bifurcation point of the $\lb$-family \eqref{eq:3Deva-2} corresponding to at
least one (or even to infinitely many) $p\in\mP$.

To depict this fact, we choose the $k$-parametric family of functions $p_k\in\mP$ defined by
$p_k(t):=\max\{-1,\min\{1,1-k/2-t\},\min\{1,1-k/2+t\}\}$
for $k\ge 0$, some of which are drawn in the left panel of Figure \ref{fig:1}. The right panel
shows an accurate numerical approximation to the map $k\mapsto\tilde\lb(\theta_0,p_k)$. Note that
the monotonicity and continuity of this map on
$[0,\infty)$ is ensured by Proposition \ref{prop:3lambda}(i),(iv),
and that, since $\lim_{k\to\infty}p_k=-1$ in $\mP$, Proposition \ref{prop:3lambda}(vi) guarantees that  $\lim_{k\to\infty}\tilde\lb(\theta_0,p_k)=\tilde\lb(\theta_0,-1)$. (As see in Figure \ref{fig:1}, this
approach is quite slow.) However, despite this continuity in $k$,
Proposition \ref{prop:3lambda}(vi) also shows that any $p_k$ with $\tilde\lb(\theta_0,p_k)>
\tilde\lb(\theta_0,-1)$ is a discontinuity point for the map $p\mapsto\tilde\lb(\theta_0,p)$ defined on $\mP$,
and that $p$ can be written as the limit in $\mP$ of a sequence of (continuity) points
$(p_k^n)$ with $\tilde\lb(\theta_0,p_k^n)=\tilde\lb(\theta_0,-1)$.
We can also check that $(\theta_0,1)$ is a discontinuity point of the
global map $(\theta,p)\mapsto\tilde\lb(\theta,p)$ on $\W$: we choose $k$ with $\tilde\lb(\theta_0,p_k)<\tilde\lb(\theta_0,1)$,
write $\theta_0=\lim_{n\to\infty}\theta_0{\cdot}t_n$ in the minimal set $\Theta$
for a sequence $(t_n)\uparrow\infty$, and note that $\lim_{t\to\infty}p_k{\cdot}t_n=1$ in $\mP$ and that $\tilde\lb(\theta_0,p_k)=\tilde\lb(\theta_0{\cdot}t_n,p_k{\cdot}t_n)$ for all $n\in\N$ (see Proposition \ref{prop:3lambda}(iii)).
So, $\lim_{n\to\infty}(\theta_0{\cdot}t_n,p_k{\cdot}t_n)=(\theta_0,1)$ in $\W$ but $\lim_{n\to\infty}\tilde\lb(\theta_0{\cdot}t_n,p_k{\cdot}t_n)=\tilde\lb(\theta_0,p_k)<\tilde\lb(\theta_0,1)$.

\begin{figure}[h]
\centering
\includegraphics[width=\textwidth]{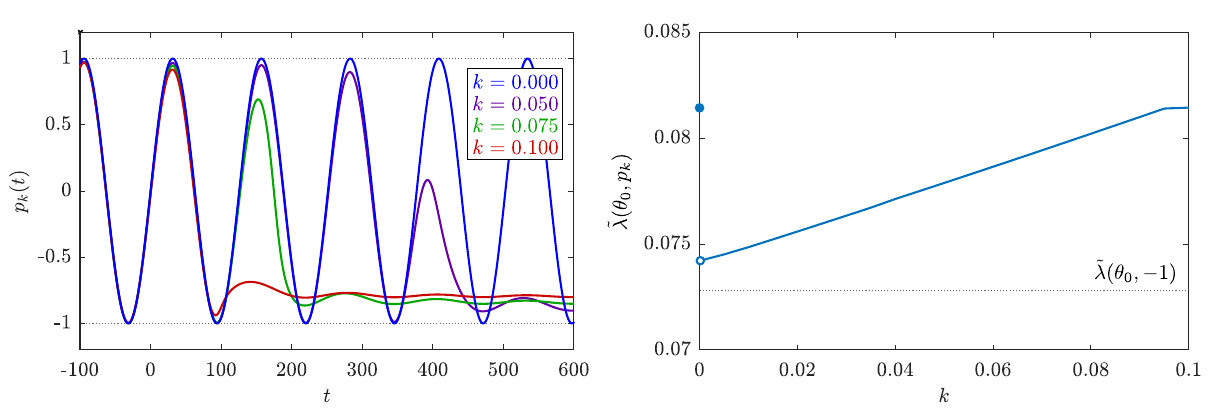}
\caption{In the left panel, $p_0(t):=\sin(t/20)$ and $p_k(t):=(\sin(t/20)-2k+1)(1/2-(1/\pi)\,\arctan(kt-1/k))+2k-1$
for different choices of $k$ in different colors. They belong to
the set $\mP$ given by \eqref{def:3Deva} with $k_1=1$ and $k_2=1$. In the right panel,
a numerical approximation of the bifurcation points $\tilde\lb(\theta_0,p_k)$ for \eqref{eq:concavemodelhull}, with
$r$, $K$, $b$ and $a$ as in Figure \ref{fig:1} and $\G(t):=(-160-16\sin(2\pi t))(1+2\exp(-10\,t^2))$.
A discontinuity  can be observed at $k=0$.
}
\label{fig:2}
\end{figure}

\begin{figure}[h]
\centering
\includegraphics[width=\textwidth]{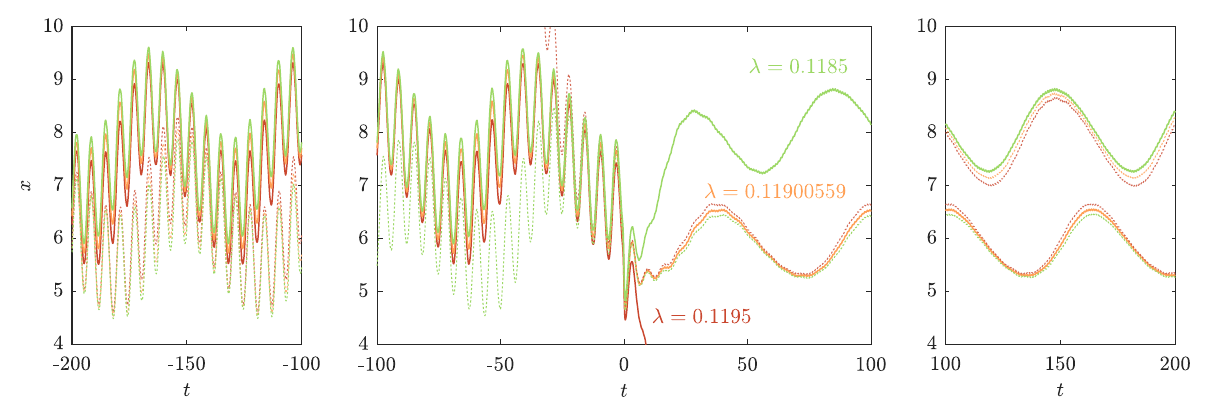}
\caption{Numerical approximations to the representative solutions of a transition equation.
The left panel (resp. right panel) represents solutions of the past equation (resp. future equation),
which is obtained by replacing $\G$ by $\G_-$ (resp. by $\G_+$) and $p_k$ by $p_k^-$ (resp. by $p_k^+$)
in the ``transition" equation \eqref{eq:concavemodel}, whose solutions are represented in the central panel.
The choice of the coefficients is similar to those in Figure \ref{fig:2}, changing $\sin(t/20)$ by $\sin(t)$ in the expression
of $p_k$. These drawings correspond to a fixed $p_k$: we take
$k=0.5$. Each of the three panels depicts two solutions corresponding to three different values of $\lb$:
the numerical approximation $0.11900559$ to $\tilde\lb(\theta_0,p_{0.5})$ (in orange), the smaller values $0.1185$ (in green)
and the larger value $0.1195$ (in red). In the left panel, the two hyperbolic solutions for each of the three values
of the parameter are plotted: in solid lines the attractive ones, in dotted lines the repulsive ones.
These attractive solutions are approximated by the locally pullback attractive solutions $a_\lb$
of the transition equation, plotted in solid lines in the central panel. The middle and right panels show the three possible
outcomes of the transition. In the first one, occurring for $\lb=0.1185$, $a_\lb$
approximates the attractive hyperbolic solution of the future equation, depicted in solid green in the
right panel: this is {\sc Case A}-tracking, $\lb=0.1185$, and it is stable for small variation on $\lb$.
The second one, also stable, and occurring for $\lb=0.1195$,
corresponds to the case in which $a_\lb$ is unbounded as time increases: this is {\sc Case C}-tipping.
The third possibility, which only occurs for $\lb=\lb(\theta_0,p_k)$, is extremely unstable situation in which
$a_\lb$ approximates the repulsive hyperbolic solution of the future equation, depicted in solid
orange in the right panel: this is {\sc Case B}. The dotted lines in the central panel depict the
locally pullback repulsive solutions of the transition equation (which is also $a_\lb$ for $\lb=\tilde\lb(\theta_0,p_{0.5})$).
The dotted lines in the right panel represent those hyperbolic
solutions not approached by $a_\lb$ as time increases (the upper curve being the attractive solution for each color).\hspace{-1cm}~}
\label{fig:tippingconcave}
\end{figure}
\vspace{.2cm}
\noindent\underline{Bifurcations and critical transitions in an asymptotically quasiperiodic case}.
We \linebreak now perform another numerical simulation for the same maps $r$, $K$ and $b$,
but for a not quasiperiodic map $\G$, which
precludes the application of the results of Proposition \ref{prop:3lambda}(v)-(viii).
We will numerically find discontinuities of $p\mapsto\tilde\lb(\theta_0,p)$ in the neighborhood of $p_0(t):=\sin(t/20)$,
and the procedure can be repeated with other functions, leading to the conclusion that the strongly discontinuous behavior
of the semicontinuous map $\tilde\lb$ occurs in more general cases than the almost periodic one.
We take $\G(t):=(-160-16\sin(2\pi t))(1+2\exp(-10\,t^2))$. While the previous choice of $\G$ corresponds
to a seasonal increase in the predator population, this one corresponds to the overlap of the seasonal increase
with a transient increase in predator numbers over a time interval comprising several seasons, which may be
due to a variety of causes, from a change in the attractiveness of the breeding colony to adverse winds hindering
certain migratory routes (see \cite{rappole} for causes of migrant avian population change).
We take $p_k(t):=(\sin(t/20)-2k+1)(1/2-(1/\pi)\,\arctan(kt-1/k))+2k-1$ for each $k\in(0,1]$, and check that they
belong to the set $\mP$ given by \eqref{def:3Deva} with $k_1=k_2=1$ and that $\lim_{t\to\infty}p_k(t)=2k-1$
and $\lim_{t\to-\infty}(p_k(t)-\sin(t/20))=0$. In addition, $(0,1]\to C(\R,\R)$, $k\mapsto p_k$ is continuous in uniform topology,
and $\lim_{k\to0^+}p_k=p_0$ in $\mP$ (but not uniformly on $\R$).
The left panel of Figure \ref{fig:2} depicts some of these functions $p_k$, and
the right panel is a numerical approximation to the map $k\mapsto\tilde\lb(\theta_0,p_k)$
that shows its discontinuity as $k\downarrow 0$. Of course, the map is continuous in $(0,1]$, as Proposition \ref{prop:3lambda}(iv)
guarantees.

Changing $p_k$ in the second example by $p_k(t):=(\sin(t)-2k+1)(1/2-\arctan(kt-1/k)/\pi)+2k-1$, provides the same global
dynamical situation (but the analogue to the left panel in Figure \ref{fig:2} is less clear). With these new data,
\eqref{eq:concavemodel} can be interpreted in the light of
the mathematical theory of critical transitions (see \cite{aspw}, \cite{lnor,lno2,lno3}, \cite{dno3,dno4}).
That is, $p_k$ can be understood as a transition function between $p_k^-(t):=\sin(t)$ and $p_k^+(t):=2k-1$,
and $\G$ can be understood as a transition function with the same past and future $\G_-(t)=\G_+(t):=-160-16\sin(2\pi t)$.
According to \cite[Theorem 4.7]{dno4}, and provided that the past and the future equations have
three hyperbolic solutions if $\lb\in(\tilde\lb(\theta_0,p_k)-\ep,\tilde\lb(\theta_0,p_k)+\ep)$ for an $\ep>0$,
a value $\lb\in(\tilde\lb(\theta_0,p_k)-\ep,\tilde\lb(\theta_0,p_k))$ corresponds to the global behavior
called {\em tracking} or {\sc Case A}, $\lb=\tilde\lb(\theta_0,p_k)$ to {\sc Case B},
and $\lb\in(\tilde\lb(\theta_0,p_k),\tilde\lb(\theta_0,p_k)+\ep)$ to {\em tipping} or {\sc Case C}. Figure \ref{fig:tippingconcave}
depicts the three possibilities, focusing on the limiting behavior of the upper solution $a_\lb$ of
\eqref{eq:concavemodel}$^\lb$ that is bounded as time decreases
(which is locally pullback attractive: see \cite[Theorem 4.6(i)]{dno4}).
So, from this point of view, focused on critical transitions, the discontinuity
properties observed in $p\mapsto\tilde\lb(\theta_0,p)$ mean also that a small variation of $p$ in the topology of $\mP$
may cause a strong variation in the value of the unique critical point $\tilde\lb(\theta_0,0)$ of the
corresponding parametric family of equations.
\section{A global bifurcation diagram in the d-concave in measure case}\label{sec:4}
\begin{figure}[h]
\centering
\includegraphics[width=\textwidth]{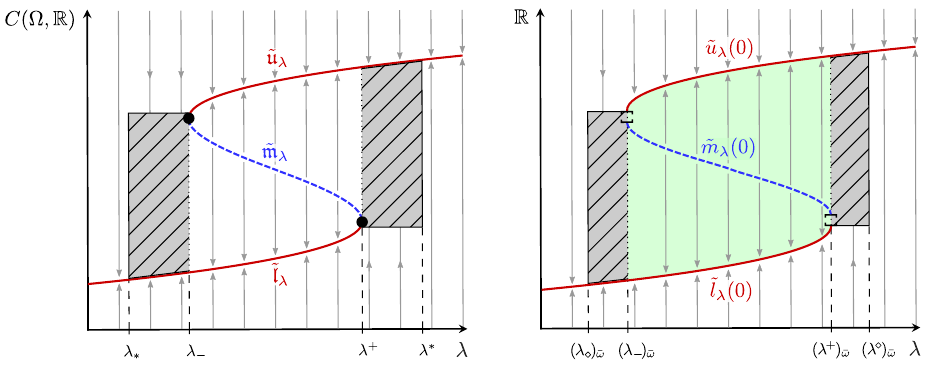}
\caption{In the left panel, the double generalized saddle-node bifurcation diagram
of the $\lb$-parametric family $x'=\mh(\wt,x)+\lb\,\mg(\w)$, $\w\in\W$, described in Theorem \ref{th:4Dbifur}.
The strictly increasing solid red curves represent two $\lb$-families of attractive hyperbolic copies
of the base: $\tilde\ml_\lb$ for $\lb<\lb^+$ and $\tilde\muk_\lb$ for $\lb>\lb_-$.
The strictly decreasing dashed blue curve represents the $\lb$-family of repulsive hyperbolic copies
of the base $\tilde\mm_\lb$ for $\lb\in(\lb_-,\lb^+)$.
The large black point at $\lb_-$ depicts the not necessarily continuous maps
$\muk_{\lb_-}:=\lim_{\lb\downarrow\lb_-}\tilde\muk_\lb$ and $\mm_{\lb_-}:=\lim_{\lb\downarrow\lb_-}\tilde\mm_\lb$,
which are respectively upper and lower semicontinuous and satisfy
$\inf_{\w\in\W}(\muk_{\lb_-}(\w)-\mm_{\lb_-}(\w))=0$. The global attractor $\mA_\lb$
reduces to the hyperbolic copy of the base $\tilde\ml_\lb=\tilde\muk_\lb$ for $\lb<\lb_*\le\lb_-$.
If $\lb_*<\lb_-$ (as drawn here) and $\lb\in[\lb_*,\lb_-)$, then $\inf_{\w\in\W}(\muk_\lb(\w)-\ml_\lb(\w))=0$;
and there exists at least an ergodic measure $m_0\in\merg$ such that $m_0(\{\w\in\W\mid\,\ml_\lb(\w)<\muk_\lb(\w)\})=1$.
In addition, $\lb\to \muk_\lb\gneq \ml_\lb$ is strictly increasing on $(\lb_*,\lb_-]$.
These facts are depicted by the gray-striped region over $[\lb_*,\lb_-)$.
The meanings of the large black point at $\lb^+$ and the gray striped region over $(\lb^+,\lb^*]$ are analogous.
The grey light arrows partly depict the dynamics of the rest of orbits.
\vspace{.15cm}\newline
In the right panel, the double generalized saddle-node bifurcation diagram of the $\lb$-family
$x'=\mh(\bwt,x)+\lb\,\mg(\bwt)$ of nonautonomous equations described by Theorems \ref{th:4Dbifuruna} and
\ref{th:4Dasint}.
The strictly increasing solid red curves represent the values at $t=0$ of the lower and
upper bounded (and hyperbolic) solutions: $\tilde l_\lb$ for $\lb<(\lb^+)_{\bar\w}\in[\lb^+,\lb^*]$
and $\tilde u_\lb$ for $\lb>(\lb_-)_{\bar\w}\in[\lb_*,\lb_-]$.
The strictly decreasing dashed blue curve represents the value at $t=0$ of the repulsive hyperbolic
solution $\tilde m_\lb$ for $\lb\in((\lb_-)_{\bar\w},(\lb^+)_{\bar\w})$.
The vertical intervals in brackets over $(\lb_-)_{\bar\w}$ and $(\lb^+)_{\bar\w}$,
which can be degenerate or not, are analogous to that of Figure \ref{fig:diagramasconcavos}.
There is only one bounded (and hyperbolic) solution $\tilde l_\lb=\tilde u_\lb$ for
$\lb<(\lb_\diamond)_{\bar\w}\in[\lb_*,(\lb_-)_{\bar\w}]$. If $(\lb_\diamond)_{\bar\w}<(\lb_-)_{\bar\w}$
(as drawn here), then $\lb\mapsto u_\lb(0)$ is a strictly increasing right-continuous map on
$[(\lb_\diamond)_{\bar\w},(\lb_-)_{\bar\w})$ with graph contained in the gray-striped region, and
in addition $\lim_{t\to\infty}(u_\lb(t)-\tilde l_\lb(t))=0$ for all $\lb\in[(\lb_\diamond)_{\bar\w},(\lb_-)_{\bar\w})$.
The meanings of the elements to the right of $(\lb^+)_{\bar\w}$ are analogous.
The green-shadowed area represents the initial data of globally bounded solutions (excepting
those filling the regions between $l_\lb$ and $u_\lb$ of the grey areas),
and the light gray arrows give some clues about the dynamics of the rest of the solutions.
}
\label{fig:diagramasdconcavos}
\end{figure}

With a structure very similar to that of Section \ref{sec:3}, this section describes
the bifurcation diagrams for parametric families
of skewproduct flows and of admissible processes determined by scalar ODEs for which the
derivatives with respect to the state variable of the corresponding laws satisfy certain
concavity properties. The main theoretical results are Theorems \ref{th:3Cbifur}, \ref{th:4Dbifuruna} and
\ref{th:4Dasint}, and their conclusions are partly depicted in Figure \ref{fig:diagramasdconcavos}.
We apply our conclusions to the analysis of an electrical circuit in Section \ref{subsec:41} and
to the study of a critical transition in Section \ref{subsec:42}.
\par
So, we work with the family
\begin{equation}\label{eq:4ini}
 x'=\mh(\wt,x)\,,\quad\w\in\W\,,
\end{equation}
assuming now that $\mh\colon\WR\to\R$ satisfies (all or part of) the next conditions:
\begin{enumerate}[leftmargin=20pt,label=\rm{\bf{d\arabic*}}]
\item\label{d1} $\mh\in C^{0,2}(\WR,\R)$,
\item\label{d2}  $\limsup_{x\to\pm\infty} (\pm\mh(\w,x))<0$ uniformly on $\W$,
\item\label{d3} $m(\{\w\in\W\,|\;x\mapsto \mh_x(\w,x) \text{ is concave}\})=1$
for all $m\in\merg$,
\item\label{d4} $m(\{\w\in\W\,|\; x\mapsto \mh_{xx}(\w,x)$ is strictly decreasing
on $\mJ\})>0$ for all compact interval $\mJ\subset\R$ and all $m\in\merg$.
\end{enumerate}
We refer to the case of the equation \eqref{eq:4ini} given by
a map $\mh$ satisfying \ref{d1} and \ref{d3} as the {\em d-concave in measure case}.
Note that \ref{d4} precludes the {\em quadratic case\/}
$\mh(\w,x)=\mc(\w)\,x^2+\md(\w)\,x+\mathfrak e(\w)$.
\begin{nota}\label{rm:4Dtambien}
As explained in Remark \ref{rm:3Ctambien}, if $\mh$ satisfies
\hyperlink{d1}{\bf dj} for a $j\in\{1,2,3,4\}$ and $\W_0\subset\W$ is
a nonempty compact $\sigma$-invariant subset, then also the restriction
$\mh\colon\W_0\times\R\to\R$ also satisfies \hyperlink{d1}{\bf dj}.
\end{nota}
Some of the results of \cite{dno4} will we repeatedly used in the proofs of
this section. Let us mention now two of them, fundamental to understand the rest.
The existence and properties of the global
attractor $\mA$ of the skewproduct flow $\tau$
defined from \eqref{eq:4ini} under \ref{d1} and \ref{d2} are analyzed in
\cite[Proposition 5.5]{dno4}, which proves that any forward $\tau$-semiorbit
is globally defined and that the maps $\ml,\muk\colon\W\to\R$ such that
\[
 \mA=\{(\w,x)\,|\;\w\in\W\text{ and }\ml(\w)\le x\le\muk(\w)\}
\]
are lower and upper semicontinuous equilibria for $\tau$. It also shows that this
global attractor is the set of globally bounded $\tau$-orbits,
and proves some fundamental comparison properties.
On the other hand, \cite[Theorem 5.3]{dno4} shows that, under conditions
\ref{d1}, \ref{d3} and \ref{d4}, there exist
at most three disjoint and ordered $\tau$-invariant compact sets $\mK_1<\mK_2<\mK_3$
projecting onto $\W$, in which case they are hyperbolic copies of the base:
$\mK_1$ and $\mK_3$ attractive,
and $\mK_2$ repulsive. In addition, in this case, $\mK_1$ and $\mK_2$ are the graphs
of the lower and upper bounded equilibria, $\ml$ and $\muk$.

Our bifurcation pattern will be provided by a parametric family
\begin{equation}\label{eq:4Dparam}
 x'=\mh(\wt,x)+\lb\,\mg(\w)
\end{equation}
for $\lb\in\R$, where $\mg\colon\W\to\R$ is continuous and satisfies $\inf_{\w\in\W}\mg(\w)>0$.
We assume that $\mh\in C^{0,i}(\WR,\R)$ for $i=1$ or $i=2$ and satisfies
\begin{list}{}{\leftmargin 25pt}
\item[\hypertarget{d2*}{{\bf d2$^*$}}] $\lim_{x\to\pm\infty}(\pm\mh(\w,x))=-\infty$ uniformly on $\W$.
\end{list}
It is easy to check that $\mh^\lb:=\mh+\lb\,\mg$ belongs to $C^{0,i}(\WR,\R)$ and satisfies \ref{d2}
for all $\lb\in\R$. Let $\tau_\lb\colon\mV_\lb\subseteq\R\times\WR\to\WR,\,(t,\w,x)\mapsto(\wt,v_\lb(t,\w,x))$ be the
(possibly local) corresponding skewproduct flow, let $\mA_\lb\subset\WR$ be the global $\tau_\lb$-attractor,
and let $\ml_\lb,\muk_\lb\colon\W\to\R$ be given by $\mA_\lb=\bigcup_{\w\in\W}(\{\w\}\times[\ml_\lb(\w),\muk_\lb(\w)])$.
This notation is used in Theorem \ref{th:4Dbifur}, which can be understood as an extension
of \cite[Theorem 5.10]{dno1} to the case of a more general map $\mh$ and
of a family of skewproduct flows $(\WR,\tau_\lb)$ defined over a base $(\W,\sigma)$ which
is not necessarily minimal or transitive. Propositions
\ref{prop:4Dsolocoer} and \ref{prop:4Dsuma}
establish some properties required in the main proofs, which hold under less restrictive hypotheses.
\begin{prop}\label{prop:4Dsolocoer}
Let $\mh\in C^{0,1}(\WR,\R)$ satisfy \hyperlink{d2*}{\bf d2$^*$}. Then,
\begin{itemize}[leftmargin=20pt]
\item[\rm(i)] for every $\w\in\W$, the maps $\lb\mapsto\ml_\lb(\w)$ and
$\lb\mapsto \muk_\lb(\w)$ are strictly increasing on $\R$ and they are, respectively, left- and right-continuous.
\item[\rm(ii)] $\lim_{\lb\to\pm\infty} \ml_\lb(\w)=\lim_{\lb\to\pm\infty}\muk_\lb(\w)=\pm\infty$ uniformly on $\W$.
\end{itemize}
\end{prop}
\begin{proof}
(i) The properties of $\lb\to\ml_\lb,\muk_\lb$ can be checked by using \cite[Proposition 5.5(iv),(v)]{dno4}, as
in the proof of \cite[Theorem 5.5(i)]{dno1}.
\smallskip

(ii) We take $n\in\N$ and look for $\lb_n\in\R$ large enough to guarantee
that $\mh(\w,x)+\lb\,\mg(\w)>1$ for all $\w\in\W$ if $x\le n$ and $\lb\ge\lb_n$.
According to \cite[Proposition 5.5(ii)]{dno4}, $\muk_{\lb}(\w)\ge\ml_{\lb}(\w)\ge n$ for
all $\w\in\W$ and $\lb\ge\lb_n$. Similarly, if $\bar\lb_n$ satisfies
$\mh(\w,x)+\lb\,\mg(\w)<-1$ for all $\w\in\W$ if $x\ge-n$ and $\lb\le\bar\lb_n$,
we deduce that $\ml_\lb(\w)\le\muk_\lb(\w)\le-n$ for all $\w\in\W$ and all $\lb\le\bar\lb_n$.
These properties prove (ii).
\end{proof}
\begin{prop} \label{prop:4Dsuma}
Let $\mh$ satisfy \ref{d1}, let us fix $m\in\merg$, and
assume that $m(\{\w\in\W\,|\;x\mapsto \mh_x(\w,x) \text{ is concave}\})=1$. For
$\lb_1<\lb_2$, let $\mb_1,\mb_2\colon\W\to\R$ be two bounded $m$-measurable equilibria
for $\tau_{\lb_1}$ and $\tau_{\lb_2}$ respectively, and such that
$\mb_1(\w)<\mb_2(\w)$ for $m$-a.a.~$\w\in\W$. Then,
\[
\int_\W \mh_x(\w,\mb_1(\w))\, dm+\int_\W \mh_x(\w,\mb_2(\w))\; dm<0\,.
\]
\end{prop}
\begin{proof} The bound $\inf_{\w\in\W}\mg(\w)>0$ allows us to
adapt the proof of \cite[Proposition 5.9]{dno1},
in turn based on that of \cite[Proposition 4.4]{dno1}.
(Observe that the minimality of the base is not used in those proofs.)
\end{proof}
The information provided by Remark~\ref{rm:4Dtambien} is required to check that the
statement Theorem \ref{th:4Dbifur}(v), regarding a restricted flow, makes perfect sense.
\begin{teor}\label{th:4Dbifur}
Let $\mh\colon\WR\to\R$ satisfy \ref{d1}, \hyperlink{d2*}{\bf d2$^*$},
\ref{d3} and \ref{d4}.
Assume that there exists $\lb_0\in\R$ such that $\tau_{\lb_0}$
admits three hyperbolic copies of the base. Then, there exists a bounded interval
$\mI=(\lb_-,\lb^+)$ with $\lb_0\in \mI$ and two real values $\lb_*$ and $\lb^*$ with
$\lb_*\le\lb_-$ and $\lb^+\le\lb^*$ such that
\begin{itemize}[leftmargin=20pt]
\item[\rm(i)] for every $\lb\in\mI$, there exists a continuous equilibrium $\tilde\mm_\lb\colon\W\to\R$
such that $\{\tml_\lb\}$, $\{\tmm_{\lb}\}$ and $\{\tmuk_{\lb}\}$ are three disjoint hyperbolic $\tau_\lb$-copies
of the base, with $\tml_\lb<\tmm_\lb<\tmuk_\lb$. In addition, $\{\tmm_\lb\}$ is repulsive,
$\{\tml_\lb\}$ and $\{\tmuk_\lb\}$ are attractive with
$\mA_\lb=\bigcup_{\w\in\W}(\{\w\}\times[\tml_\lb(\w),\tmuk_\lb(\w)])$, and
$\lb\mapsto\tml_\lb,-\tmm_\lb,\tmuk_\lb$ are continuous and strictly increasing on $\mI$.
In particular, if $\lb_1,\lb_2\in\mI$ and $\lb_1<\lb_2$, then
\begin{equation}\label{des:4Dchain}
 \tml_{\lb_1}<\tml_{\lb_2}<\tmm_{\lb_2}<\tmm_{\lb_1}<\tmuk_{\lb_1}<\tmuk_{\lb_2}\,.
\end{equation}
\end{itemize}
Let us define $\mm_{\lb_-}(\w):=\lim_{\lb\to(\lb_-)^+}\tmm(\w)$ and
$\mm_{\lb^+}(\w):=\lim_{\lb\to(\lb^+)^-}\tmm(\w)$. Then,
\begin{itemize}[leftmargin=20pt]
\item[\rm(ii)] $\{\tml_\lb\}:=\{\ml_\lb\}$
is an attractive hyperbolic $\tau_\lb$-copy of the base for all $\lb<\lb^+$,
the map $(-\infty,\lb_+)\to C(\W,\R),\,\lb\mapsto\tml_\lb$ is continuous and strictly increasing,
$\lim_{\lb\to-\infty}\tml_\lb=-\infty$ uniformly on $\W$,
$\lim_{\lb\to(\lb^+)^-}\tml_\lb=\ml_{\lb^+}$ pointwise on $\W$,
and $\inf_{\w\in\W}(\mm_{\lb^+}(\w)-\ml_{\lb^+}(\w))=0$.
Analogously, $\{\tmuk_\lb\}:=\{\muk_\lb\}$ is
an attractive hyperbolic copy of the base for all $\lb>\lb_-$,
$(\lb_-,\infty)\to C(\W,\R),\,\lb\mapsto\tmuk_\lb$ is continuous and strictly increasing,
$\lim_{\lb\to\infty}\tmuk_\lb=\infty$ uniformly on $\W$,
$\lim_{\lb\to(\lb_-)^+}\tmuk_\lb=\muk_{\lb_-}$ pointwise on $\W$, and
$\inf_{\w\in\W}(\muk_{\lb_-}(\w)-\mm_{\lb_-}(\w))=0$.
\item[\rm(iii)] $\mA_{\lb}$ is a hyperbolic copy of the base if and only if $\lb\notin[\lb_*,\lb^*]$.
\item[\rm(iv)] If $\lb_*<\lb_-$ and $\lb\in[\lb_*,\lb_-)$ or if $\lb^+<\lb^*$ and $\lb\in(\lb^+,\lb^*]$,
then $\inf_{\w\in\W}(\muk_\lb(\w)-\ml_{\lb}(\w))=0$, and there exists at least a measure $m_0\in\merg$
such that $m_0(\{\w\in\W\,|\;\ml_\lb(\w)<\muk_\lb(\w)\})=1$.
\end{itemize}
In addition,
\begin{itemize}[leftmargin=20pt]
\item[\rm(v)] if $\W_\w$ is the closure of the $\sigma$-orbit of $\w\in\W$, and $\mI_\w=((\lb_-)_\w,(\lb^+)_\w)$ is
the interval associated by {\rm (i)} to the restriction of the $\lb$-parametric family \eqref{eq:4Dparam} to $\W_\w$,
then $\lb_*=\inf_{\w\in\W}(\lb_-)_\w$ and $\lb_-=\sup_{\w\in\W}(\lb_-)_\w$. Analogously,
$\lb^+=\inf_{\w\in\W}(\lb^+)_\w$ and $\lb^*=\sup_{\w\in\W}(\lb^+)_\w$.
\item[\rm(vi)] If $\W$ is minimal, then $\lb_-=\lb_*$ and $\lb^+=\lb^*$.
\end{itemize}
\end{teor}
\begin{proof}
We will just work for $\lb\ge\lb_0$. All the arguments are analogous for $\lb\le\lb_0$.
\smallskip

(i) The persistence of hyperbolic copies of the base under small variations of $\lb$
(see Theorem \ref{th:2persistencia}) ensures that $\mI_+:=
\left\{\bar\lb\in\R\,|\; \text{ there exist three hyperbolic}\right.$
$\tau_{\lb}$-copies of the base $\left.\text{for all
$\lb\in[\lb_0,\bar\lb)$}\right\}$ contains $[\lb_0,\lb_0+\delta]$
for a $\delta>0$. For the same reason, $\lb^+:=\sup\mI_+$ (which is
finite, as we will check below) cannot belong to $\mI_+$. This will be
the $\lb^+$ of the statement.

Let us take $\lb\in\mI_+$, and let $\{\tmm_\lb\}$ be the middle hyperbolic
copy of the base. To check that the upper and lower hyperbolic copies of the base,
say $\{\tilde\mc_\lb\}$ and $\{\tilde\ma_\lb\}$, are precisely $\{\muk_\lb\}$ and $\{\ml_{\lb}\}$
(i.e., they are given by the upper and lower equilibria of the global $\tau_{\lb}$-attractor), we consider
the compact sets $\mK^u_\lb:=\{(\w,x)\in\mA_\lb\,|\;x\ge\tilde\mc_\lb(\w)\}$,
$\{\tmm_\lb\}$, and $\mK^l_\lb:=\{(\w,x)\in\mA_\lb\,|\;x\le\tilde\ma_\lb(\w)\}$,
and apply \cite[Theorem 5.3]{dno4} to conclude that $\mK^u_\lb=\{\muk_\lb\}$
and $\mK^l_\lb=\{\ml_\lb\}$, both attractive, and that $\{\tmm_\lb\}$ is repulsive.
We call $\tml_\lb:=\ml_\lb$ and $\tmuk_\lb:=\muk_\lb$, and deduce from
Theorem \ref{th:2persistencia} and Proposition \ref{prop:4Dsolocoer}(i) that
$\lb\to\tml_\lb,\tmuk_\lb$ are continuous and increasing on $\mI_+$.

Theorem \ref{th:2persistencia} also yields the continuity of
$\lb\mapsto\tmm_{\lb}$ on $\mI_+$.
Let us take $\lb_1<\lb_2$ in $\mI_+$. By repeating the argument in the proof of
\cite[Theorem 5.10(i)]{dno1}, we conclude that, for any minimal subset $\mM\subseteq\W$,
$\tmm_{\lb_1}|_{\mM}>\tmm_{\lb_2}|_{\mM}$ (see Remark \ref{rm:4Dtambien}).
Now, we assume for contradiction that $\tmm_{\lb_1}(\w)\le\tmm_{\lb_2}(\w)$ for a point $\w\in\W$,
observe that a standard comparison argument yields $\tmm_{\lb_1}(\wt)\le\tmm_{\lb_2}(\wt)$
for all $t\ge 0$, and deduce that $\tmm_{\lb_1}(\bar\w)\le\tmm_{\lb_2}(\bar\w)$ for
a point $\bar\w$ in a minimal subset of the \upomeg-limit set of the orbit
$\{\wt\,|\;t\in\R\}$, which is impossible. This completes the proof of
\eqref{des:4Dchain} on $\mI_+$.
In turn, \eqref{des:4Dchain} ensures that $\ml_{\lb}$ is bounded while
$\lb\in\mI_+$. Therefore, Proposition \ref{prop:4Dsolocoer}(ii)
ensures that $\lb^+=\sup\mI_+$ is finite.
Analogous arguments to the left of $\lb_0$ complete the proof of (i).
\smallskip

(ii) It follows from \eqref{des:4Dchain} that the pointwise limits
$\hat\muk_{\lb^+}:=\lim_{\lb\to(\lb^+)^-}\tmuk_\lb$,
$\hat\mm_{\lb^+}:=\lim_{\lb\to(\lb^+)^-}\tmm_\lb$ and $\hat\ml_{\lb^+}:=\lim_{\lb\to(\lb^+)^-}\tml_\lb$
are globally defined bounded maps, with $[\hat\ml_{\lb^+}(\w),\hat\mm_{\lb^+}(\w)]\subseteq
[\tml_{\lb}(\w),\tmm_{\lb}(\w)]$ for all $\w\in\W$ and $\lb\in\mI$ and with
$\hat\muk_{\lb^+}\le\muk_{\lb^+}$. In addition,
Proposition \ref{prop:4Dsolocoer}(i) shows that $\hat\ml_{\lb^+}=\ml_{\lb^+}$.

We already know that $\{\tmuk_{\lb}\}$ is a hyperbolic copy of the base if
$\lb\in\mI=(\lb_-,\lb^+)$. Let us fix $\lb\ge\lb^+>\lb_0$. Since
$\muk_{\lb}>\tmuk_{\lb_0}>\tmm_{\lb_0}$ (see Proposition \ref{prop:4Dsolocoer}(i)
and \eqref{des:4Dchain}), the compact closure $\mK^u_{\lb}$ of the graph of $\muk_\lb$
is strictly above $\{\tmm_{\lb_0}\}$. Proposition~\ref{prop:4Dsuma} applied to
$\tmm_{\lb_0}$ (whose graph is a repulsive hyperbolic copy of the base),
any $m\in\merg$ and any $m$-measurable $\tau_\lb$-equilibrium with graph in $\mK^u_\lb$
shows that the upper Lyapunov exponent of $\mK^u_\lb$ is strictly negative
(see Section \ref{subsec:skew}).
Let us check that the section of $\mK^u_\lb$ over a point $\w_0$ belonging to a minimal set
$\mM\subseteq\W$ reduces to a point.
We denote $(\mK^u_\lb)^\mM:=\{(\w,x)\in\mK^u_\lb\,|\;\w\in\mM\}$, which is a
compact invariant set for the restriction of $\tau_\lb$ to $\mM\times\R$.
Clearly, $x\ge\tmuk_{\lb_0}(\w)>\tmm_{\lb_0}(\w)$
for all $(\w,x)\in(\mK^u_\lb)^\mM$.  According to \cite[Theorem 5.10]{dno1}, and since
$\mM$ is minimal, $(\mK^u_\lb)^\mM$ must coincide with the graph of $\muk_\lb|_{\mM}$,
from where the assertion follows. Hence, Theorem \ref{th:2copia} ensures that
$\mK^u_\lb$ is an attractive hyperbolic copy of the base.
Now, an easy contradiction argument shows that
$\muk_\lb$ is continuous, and hence that $\mK^u_\lb=\{\muk_\lb\}$.
The same argument works for $\hat\muk_{\lb^+}$
instead of $\muk_{\lb_+}$. Hence, $\{\hat\muk_{\lb^+}\}$ is an attractive hyperbolic copy
of the base, and so, $\hat\muk_{\lb^+}=\muk_{\lb^+}$:
otherwise, and according to Theorem \ref{th:2persistencia},
there would be four hyperbolic copies of the base for $\lb\in\mI_+$ close enough
to $\lb_+$, which contradicts \cite[Theorem 5.3]{dno4}. Note that this implies that
$(\lb_-,\infty)\to C(\W,\R),\,\lb\to\muk_\lb$ is continuous. The monotonicity and limiting
behavior as $\lb\to(\lb_-)^+$ and as $\lb\to\infty$ follow from
Proposition \ref{prop:4Dsolocoer}. We call $\tmuk_\lb:=\muk_\lb$ for all $\lb\ge\lb^+$.

To unify notation, we denote $\mm_{\lb^+}:=\hat\mm_{\lb^+}$.
Now, let us assume for contradiction that $\delta_1:=\inf_{\w\in\W}(\mm_{\lb^+}(\w)-\ml_{\lb^+}(\w))>0$.
Our goal is to check that the closures
$\mK^l_{\lb^+}$, $\mK^m_{\lb^+}$ and $\mK^u_{\lb^+}=\{\tmuk_{\lb^+}\}$
of the graphs of $\ml_{\lb^+}$, $\mm_{\lb^+}$ and $\tilde\muk_{\lb^+}$
are three different ordered compact sets projecting on $\W$.
We take $(\bar\w,\bar x)\in \mK^m_{\lb^+}$.
By reasoning as in the proof of Theorem \ref{th:3Cbifur}(ii),
we check that $\bar x\ge\ml_{\lb^+}(\bar\w)+\delta_1$. In addition, \eqref{des:4Dchain}
ensures that  $\mm_{\lb^+}<\tilde\muk_{\lb_0}<\tilde\muk_{\lb^+}$,
and hence $\bar x\le\tmuk_{\lb^+}(\bar\w)-\delta_2$
for a $\delta_2$ independent of $(\bar\w,\bar x)$. So, \cite[Theorem 5.2]{dno4} allows us to
assert that all the Lyapunov exponents of $\mK^m_{\lb^+}$ are positive (see Section \ref{subsec:skew}),
and that its upper and lower equilibria coincide on a set $\W_0$ with
$m_0(\W_0)=1$ for all $m_0\in\merg$. Therefore, Theorem \ref{th:2copia}
ensures that $\mK^m_{\lb^+}$ is a repulsive hyperbolic copy of the base.
Now, it is easy to deduce from the previous separation properties that
$\mK^l_{\lb^+}<\mK^m_{\lb^+}<\mK^u_{\lb^+}$, as asserted.
In these conditions, \cite[Theorem 5.3]{dno4}
implies that $\lb^+\in\mI_+$, which is not the case.
This contradiction completes the proof of (ii).
\smallskip

(iii)-(v) Proposition \ref{prop:4Dsolocoer}(ii) and the continuity of $\tmm_{\lb_0}$ ensure that
there exists $\bar\lb$ such that $x>\tmm_{\lb_0}(\w)$ for all
$(\w,x)\in\mA_{\lb}$ if $\lb\ge\bar\lb$. The argument applied to $\mK_\lb^u$ in the proof of (ii)
shows that $\mA_\lb$ is an attractive hyperbolic copy of the base.
Hence, $\mJ_+:=\{\bar\lb\in\R\,|\;\mA_\lb=\{\tmuk_{\lb}\}\text{ for all $\lb\ge\bar\lb$}\}$
is nonempty. We call $\lb^*:=\inf\mJ_+$. Clearly, $\lb^+\le\lb^*\notin\mJ_+$.

Let us define $(\lb^+)_\w$ as in the statement of (v).
It is clear that $\lb^+\le\inf_{\w\in\W}(\lb^+)_\w$. To check that they coincide, we
assume for contradiction that there exists $\bar\lb\in(\lb^+,\inf_{\w\in\W}(\lb^+)_\w)$.
It follows from \eqref{des:4Dchain} applied to each $\W_\w$ that
$\ml_{\lb^+}(\w)<\ml_{\bar\lb}(\w)<\mm_{\lb^+}(\w)$ for all $\w\in\W$.
Working as in the proof of Theorem \ref{th:3Cbifur}(iv),
we prove that $\inf_{\w\in\W}(\ml_{\bar\lb}(\w)-\ml_{\lb^+}(\w))>0$,
and hence $\inf_{\w\in\W}(\mm_{\lb^+}(\w)-\ml_{\lb^+}(\w))>0$, which contradicts (ii).
(When working for values less than $\lb_0$, we must bound from below the difference at $\w{\cdot}(-1)$
instead of at $\w{\cdot}1$.)

The inequality $\lb^*\ge\sup_{\w\in\W}(\lb^+)_\w$ is also clear. Now, we take
$\bar\lb>\sup_{\w\in\W}(\lb^+)_\w$ and check that $\mA_{\bar\lb}=\{\tmuk_{\bar\lb}\}$,
which proves that $\sup_{\w\in\W}(\lb^+)_\w\ge\lb^*$ and hence the equality.
Our goal is to check that, for any $m_0\in\merg$, $\ml_{\bar\lb}(\w)=\tmuk_{\bar\lb}(\w)$
for $m_0$-a.a.~$\w\in\W$. Observe that this ensures that the upper Lyapunov exponent of
$\mA_{\bar\lb}$ is given by $\int_{\W}\mh_x(\w,\tilde\muk_{\bar\lb}(\w))\,dm_0$ for an $m_0\in\merg$:
the graph of any bounded $m_0$-measurable equilibria is delimited by those
of $\ml_{\bar\lb}$ and $\tmuk_{\bar\lb}$. So, this upper Lyapunov exponent
is strictly negative, since $\{\tmuk_{\bar\lb}\}$ is an attractive hyperbolic
copy of the base, which combined with the previous property and the fact that any
minimal subset of $\W$ concentrates a $\sigma$-ergodic measure allow us
to apply Theorem \ref{th:2copia} to prove that $\mA_{\bar\lb}=\{\tmuk_{\bar\lb}\}$.

So, we assume for contradiction the existence of $\delta>0$ such that
$m_0(\{\w\in\W\,|\;\tmuk_{\bar\lb}(\w)\ge\ml_{\bar\lb}(\w)+\delta\})>0$.
Lusin's theorem and the regularity of $m_0$ provide a compact subset
$\Delta\subseteq\W$ with $m_0(\Delta)>0$ such that the
restriction $\ml_{\bar\lb}\colon\Delta\to\R$ is continuous and satisfies
$\ml_{\bar\lb}(\w)+\delta\le\tmuk_{\bar\lb}(\w)$ for all $\w\in\Delta$.
Birkhoff's ergodic theorem ensures the existence of $\bar\w\in\Delta$ and $(t_n)\uparrow\infty$
such that $\bwt_n\in\Delta$. We work in $\W_{\bar\w}$. Since
$\bar\lb>(\lb^+)_{\bar\w}$, we can repeat part of the proof of
\cite[Theorem B.3(iii)]{dlo1} to check that
$\lim_{t\to\infty}(\tmuk_{\bar\lb}(\bwt)-\ml_{\bar\lb}(\bwt))=0$, which yields
$0<\delta\le\lim_{n\to\infty}(\tmuk_{\bar\lb}(\bwt_n)-\ml_{\bar\lb}(\bwt_n))=0$:
this is the sought-for contradiction.

It remains to check the assertions concerning the possible points $\lb\in(\lb^+,\lb^*]$.
First, let us check that
$\inf_{\w\in\W}(\tmuk_\lb(\w)-\ml_\lb(\w))=0$. We fix $\bar\w\in\W$ with
$(\lb^+)_{\bar\w}<\lb$, and restrict ourselves to $\W_{\bar\w}$. As before, as in
\cite[Theorem B.3(iii)]{dlo1}, we have $\inf_{\w\in\W_{\bar\w}}(\tmuk_\lb(\w)-\ml_\lb(\w))=0$,
which proves the assertion. Second, we assume for contradiction that
$m_0(\{\w\in\W\,|\;\ml_\lb(\w)=\tmuk_\lb(\w)\})=1$ for all $m_0\in\merg$, and we deduce
as in the previous paragraph that $\mA_\lb=\{\tmuk_\lb\}$. This ensures that $\lb>(\lb^+)_\w$
for all $\w\in\W$, which combined with (v) (already proved)
and the property $\mA_{\lb^*}\ne\{\tmuk_{\lb^*}\}$ yields the contradiction $\lb>\lb^*$.
\smallskip

(vi) In the minimal case, $\W_\w=\W$ for all $\w\in\W$. This ensures that
$(\lb^+)_\w=\lb^+$  for all $\w\in\W$, and hence (v) yields $\lb^+=\lb^*$.
\end{proof}
The nonautonomous bifurcation diagram described in the previous theorem depicts
a {\em generalized double local saddle-node nonautonomous bifurcation diagram}.
The word ``generalized" refers to the fact that
the transition between the existence of three hyperbolic copies of the base as $\lb$ decreases (or increases)
and an attractor given by a hyperbolic copy of the base may occur during an interval of values of the
parameter instead of at only one value. In these intervals reduce to a point, as in the case of
minimality of $\W$, we have a {\em classic double local saddle-node nonautonomous bifurcation diagram}.

In the line of Theorem \ref{th:3Cbifuruna}, the next result adapts the previous information
to a process instead of a skewproduct flow.
\begin{teor}\label{th:4Dbifuruna}
Assume all the hypotheses of Theorem {\rm\ref{th:4Dbifur}}.
Let us fix $\bar\w\in\W$, and let $l_\lb$ and $u_\lb$ be the lower and
upper bounded solutions of $x'=\mh(\bwt,x)+\lb\,\mg(\bwt)$. Then,
there exist real numbers $\bar\lb_-<\bar\lb^+$ such that
\begin{itemize}[leftmargin=20pt]
\item[-] for $\lb\in(\bar\lb_-,\bar\lb^+)$, this equation has three uniformly separated hyperbolic solutions
$\tilde l_\lb:=l_\lb<\tilde m_\lb<u_\lb=:\tilde u_\lb$, with $\tilde l_\lb$ and $\tilde u_\lb$ attractive
and $\tilde m_\lb$ repulsive;
\item[-]  for $\lb=\bar\lb_-$ (resp.~$\lb=\bar\lb^+$), this equation has
only one hyperbolic solution $\tilde l_{\bar\lb_-}:=l_{\bar\lb_-}$
(resp.~$\tilde u_{\bar\lb^+}:=u_{\bar\lb^+}$)
which is attractive and uniformly separated from $u_{\bar\lb_-}$ (resp.~$l_{\bar\lb^+}$);
\item[-] for $\lb\notin[\bar\lb_-,\bar\lb^+]$, this equation has
at least one hyperbolic solution, namely $\tilde l_\lb:=l_\lb$ for $\lb<\bar\lb_-$
and $\tilde u_\lb:=u_\lb$ for $\lb>\bar\lb^+$, which is attractive, and which is not
uniformly separated from any other bounded solution; more
precisely, $\lim_{t\to\infty}(u_\lb(t)-l_\lb(t))=0$ if $\lb\notin[\bar\lb_-,\bar\lb^+]$.
\end{itemize}
In addition, the maps $\lb\to\tilde l_\lb,\,-\tilde m_\lb,\,\tilde u_\lb$ are continuous in the
uniform topology of $C(\R,\R)$ and strictly increasing on $(-\infty,\bar\lb^+)$, $(\bar\lb_-,\bar\lb^+)$,
$(\bar\lb_-,\infty)$ respectively, with
\begin{itemize}[leftmargin=20pt]
\item[-] $\lim_{\lb\to(\bar\lb^+)^-}\tilde l_\lb(t)=l_{\bar\lb^+}(t)$
and $\lim_{\lb\to(\bar\lb_-)^+}\tilde u_\lb(t)=u_{\bar\lb_-}(t)$ pointwise on $\R$,
\item[-] $\lim_{\lb\to-\infty}\tilde l_\lb(t)=-\infty$ and
$\lim_{\lb\to\infty}\tilde u_\lb(t)=\infty$ uniformly on $\R$,
\item[-] $\inf_{t\in\R}(u_{\bar\lb_-}(t)-m_{\bar\lb_-}(t))=0$
for the solution $m_{\bar\lb_-}(t):=\lim_{\lb\to(\bar\lb_-)^+}\tilde m_\lb(t)$ and
$\inf_{t\in\R}(m_{\bar\lb^+}(t)-l_{\bar\lb^+}(t))=0$
for the solution $m_{\bar\lb^+}(t):=\lim_{\lb\to(\bar\lb^+)^-}\tilde m_\lb(t)$,
\end{itemize}
and the values $\bar\lb_-$ and $\bar\lb^+$ are common for all $\w\in\W$ with
$\W_\w=\W_{\bar\w}$.

Furthermore,
there are two possibilities for $(-\infty,\bar\lb_-)$ (resp.~for $(\bar\lb^+,\infty)$):
\begin{itemize}[leftmargin=20pt]
\item[-] either $\tilde l_\lb$ is the unique bounded solution for all $\lb<\bar\lb_-$
(resp.~$\tilde u_\lb$ is the unique bounded solution for all $\lb>\bar\lb^+$),
\item[-] or there exists $\bar\lb_\diamond<\bar\lb_-$ such that $u_\lb>\tilde l_\lb$
for $\lb\in(\bar\lb_\diamond,\bar\lb_-)$ and $\tilde l_\lb$ is the unique bounded
solution for $\lb\in(-\infty,\bar\lb_\diamond)$
(resp.~there exists $\bar\lb^\diamond>\bar\lb^+$ such that $\tilde u_\lb>l_\lb$ for
$\lb\in(\bar\lb^+,\bar\lb^\diamond)$ and $\tilde u_\lb$ is the unique bounded
solution for $\lb\in(\bar\lb^\diamond,\infty)$);
\end{itemize}
and the first possibility is fulfilled if $\W_{\bar\w}$ is the \upalfa-limit set $\mathbb A_{\bar\w}$
of $\bar\w$ for the flow $\sigma$. More precisely, let us call $\bar\lb_\diamond:=\bar\lb_-$
(resp.~$\bar\lb^\diamond:=\bar\lb^+$) if the first possibility holds, and let
$\mI_{\mathbb A}=(\lb_-^{\mathbb A_{\bar\w}},\lb^+_{\mathbb A_{\bar\w}})$ be the interval associated
by Theorem {\rm\ref{th:4Dbifur}} to the restriction of the family \eqref{eq:4Dparam} to $\mathbb A_{\bar\w}$.
Then, $\bar\lb_-^{\mathbb A_{\bar\w}}\le\bar\lb_\diamond<\bar\lb^\diamond\le\bar\lb^+_{\mathbb A_{\bar\w}}$, and
$\bar\lb_-^{\mathbb A_{\bar\w}}=\bar\lb_\diamond<\bar\lb^\diamond=\bar\lb^+_{\mathbb A_{\bar\w}}$ if
$\W_{\bar\w}=\mathbb A_{\bar\w}$.
\end{teor}
\begin{proof}
Let $\W_{\bar\w}$ be the closure of $\{\bwt\,|\;t\in\R\}$, and $\bar\tau_\lb$ the flow
induced on $\W_{\bar\w}\times\R$ by the family
\begin{equation}\label{eq:4DCproceso}
 x'=\mh(\wt,x)+\lb\,\mg(\wt)\,,\quad \w\in\W_{\bar\w}\,.
\end{equation}
As explained in Remark \ref{rm:4Dtambien}, we can restrict Theorem \ref{th:4Dbifur}
to the family $(\W_{\bar\w},\bar\tau_\lb)$. Note first that
$l_\lb(t)=\ml_\lb(\bwt)$ and $u_\lb(t)=\muk_\lb(\bwt)$ for all $\lb\in\R$, since
the corresponding global attractor is composed by the bounded $\bar\tau_\lb$-orbits.
Theorem \ref{th:4Dbifur} and Proposition \ref{prop:2extiende}
show the existence of an interval $(\bar\lb_-,\bar\lb^+):=((\lb_-)_{\bar\w},(\lb^+)_{\bar\w})$
and a new $\bar\tau_\lb$-copy of the base
$\{\tmm_\lb\}$ for $\lb\in(\bar\lb_-,\bar\lb^+)$ such that, if $\tilde m_\lb(t):=\tilde\mm_\lb(\bwt)$, then
all the assertions of the theorem previous to the word ``Furthermore" hold, excepting the next ones:
$\lim_{t\to\infty}(u_\lb(t)-l_\lb(t))=0$ for $\lb\notin[\bar\lb_-,\bar\lb^+]$, which is proved
by repeating again part of the proof of \cite[Theorem B.3]{dlo1}; and
$\inf_{t\in\R}(u_{\bar\lb_-}(t)-m_{\bar\lb_-}(t))=\inf_{t\in\R}(m_{\bar\lb^+}(t)-l_{\bar\lb^+}(t))=0$,
which follows from the already proved assertions and \cite[Theorem 5.6]{dno4}.

Let us now analyze the situation on $(-\infty,\bar\lb_-)$. The analysis is analogous for $(\bar\lb^+,\infty)$.
Theorem \ref{th:4Dbifur}(iii) ensures the existence of a value $\bar\lb_*$ associated to the family
\eqref{eq:4DCproceso} such that $(-\infty,\bar\lb_*)\subseteq\mI:=
\{\lb<\bar\lb_-\,|\;u_\lb=\tilde l_\lb\}$. Our next goal is to check that
$\mI$ is a negative halfline. We begin by proving that $\mI$ is left-open. So, we assume that
$u_{\lb_1}=\tilde l_{\lb_1}$ for a $\lb_1<\bar\lb_-$.
We also take $\rho>0$ and $\ep>0$ such that $\rho$ is a
radius of uniform exponential stability of the hyperbolic solutions
$\tilde l_\lb$ for all $\lb\in[\lb_1-\ep,\lb_1]$ (see Theorem \ref{th:2persistencia}); and we reduce $\ep$,
if required, to ensure that $\sup_{t\in\R}(\tilde l_{\lb_1}(t)-\tilde l_\lb(t))<\rho$ if $\lb\in[\lb_1-\ep,\lb_1]$:
see Theorem \ref{th:4Dbifur}(ii). Then, for $\lb\in[\lb_1-\ep,\lb_1]$,
$\tilde l_{\lb_1}(t)-\rho<\tilde l_{\lb_1-\ep}(t)\le \tilde l_\lb(t)\le u_\lb(t)<u_{\lb_1}(t)=\tilde l_{\lb_1}(t)$,
and hence $u_\lb(t)-\tilde l_\lb(t)<\rho$. So, for any fixed $t\in\R$,
$u_\lb(t)-\tilde l_\lb(t)\le k_\lb e^{-\beta_\lb(t-s)}(u_\lb(s)-\tilde l_\lb(s))
\le \rho\,k_\lb e^{-\beta_\lb(t-s)}$ for all $s<t$ for certain constants $k_\lb\ge 1$ and $\beta_\lb>0$, which
is only possible if $u_\lb(t)-\tilde l_\lb(t)=0$. This proves the assertion. On the other hand,
the set $\{\lb<\bar\lb_-\,|\;u_\lb>\tilde l_\lb\}$ is right-open: if $u_{\lb_2}(0)>\tilde l_{\lb_2}(0)$, then the
continuity of $\lb\mapsto\tilde l_\lb(0)$ provides $\ep>0$ such that $\tilde l_{\lb_2+\ep}(0)<u_{\lb_2}(0)<u_{\lb_2+\ep}(0)$.
These two properties make it easy to check that $\mI$ is a negative halfline, as asserted.
We define $\bar\lb_\diamond:=\sup\mI$, and note that $\bar\lb_\diamond\le\bar\lb_-$ and
that the definition of $\mI$ shows that there are only the two possibilities described in the statement:
the first one holds if and only if $\bar\lb_\diamond=\bar\lb_-$.

With the notation of the last part of the statement, let us prove that $\lb_-^{\mathbb A_{\bar\w}}\le\bar\lb_\diamond$.
We take $\lb\in(\bar\lb_\diamond,\bar\lb^+)$.
Let $\rho_\lb>0$ be a radius of uniform exponential stability of $\{\tml_\lb\}$.
Then, $\liminf_{t\to-\infty}(u_\lb(t)-\tilde l_\lb(t))\ge\rho_\lb$: otherwise there exists $(s_n)\uparrow\infty$
with $0<u_\lb(-s_n)-\tilde l_\lb(-s_n)<\rho_\lb$, and hence $0<u_\lb(0)-\tilde l_\lb(0)<\rho_\lb\,k_\lb\,e^{-\beta_\lb s_n}$
for all $n\in\N$, which is impossible. Now we take $\w_0\in\mathbb A_{\bar\w}$.
We take any sequence $(t_n)\uparrow\infty$ such that $\w_0:=\lim_{n\to\infty}\bar\w{\cdot}(-t_n)$
and assume without restriction the existence of $\bar l_0:=\lim_{n\to\infty}\tilde l_\lb(-t_n)$ and
$\bar u_0:=\lim_{n\to\infty} u_\lb(-t_n)$. Then, $v_\lb(t,\w_0,\bar u_0)-v_\lb(t,\w_0,\bar l_0)=
\lim_{n\to\infty}(u_\lb(t-t_n)-\tilde l_\lb(t-t_n))\ge\rho_\lb$. The previously proved assertions, applied
to $\w_0$ instead of $\bar\w$, show that $\lb\ge(\lb_-)_{\w_0}$. Hence, Theorem \ref{th:4Dbifur}(v)
ensures that $\lb\ge\sup_{\w\in\mathbb A_{\bar\w}} (\lb_-)_\w=\lb_-^{\mathbb A_{\bar\w}}$.
Therefore, $\bar\lb_\diamond\ge \lb_-^{\mathbb A_{\bar\w}}$, as asserted.
In particular, if $\W_{\bar\w}=\mathbb A_{\bar\w}$, then
$\bar\lb_\diamond\le\bar\lb_-=(\lb_-)_{\bar\w}=\lb_-^{\mathbb A_{\bar\w}}\le
\bar\lb_\diamond$, and hence $\bar\lb_\diamond=\bar\lb_-$.
We can prove in the same way that $\bar\lb^\diamond\le \lb^+_{\mathbb A_{\bar\w}}$ and conclude that
$\bar\lb^\diamond=\bar\lb^+$ if $\mathbb A_{\bar\w}=\W_{\bar\w}$.
\end{proof}

\begin{teor}\label{th:4Dasint}
Assume the conditions of Theorem {\rm\ref{th:4Dbifuruna}} and let
$x_\lb(t,s,x)$ be the maximal solution of $x'=\mh(\bwt,x)+\lb\,\mg(\bwt)$ with $x_\lb(s,s,x)=x$,
defined on $(\alpha_{s,x,\lb},\infty)$. With the notation established in Theorem {\rm\ref{th:4Dbifuruna}}:
for any $\lb\in\R$,
\begin{itemize}[leftmargin=20pt]
\item[-] $\lim_{t\to(\alpha_{s,x,\lb})^+}x_\lb(t,s,x)=-\infty$ if and only if $x<l_\lb(s)$, and
\item[-] $\lim_{t\to(\alpha_{s,x,\lb})^+}x_\lb(t,s,x)=\infty$ if and only if $x>u_\lb(s)$;
\end{itemize}
and, if $\lb\in(\bar\lb_-,\bar\lb^+)$, then
\begin{itemize}[leftmargin=20pt]
\item[-] $\lim_{t\to\infty}(x_\lb(t,s,x)-\tilde u_\lb(t))=0$ if and only if $x>\tilde m_\lb(s)$,
\item[-] $\lim_{t\to\infty}(x_\lb(t,s,x)-\tilde l_\lb(t))=0$ if and only if $x<\tilde m_\lb(s)$, and
\item[-] $\lim_{t\to-\infty}(x_\lb(t,s,x)-\tilde m_\lb(t))=0$ if and only if $x\in(\tilde l_\lb(s),
\tilde u_\lb(s))$.
\end{itemize}
\end{teor}
\begin{proof}
Let us fix $\lb\in\R$. If $x\ge l_\lb(s)$, then $x_\lb(t,s,x)\ge l_\lb(t)$,
which precludes $\lim_{t\to(\alpha_{s,x,\lb})^+}x_\lb(t,s,x)=-\infty$.
On the contrary, if $x<l_\lb(s)$ then \cite[Proposition 5.5(v)]{dno4} proves the existence of $t$ such that
$x_\lb(t,s,x)<\rho$, where $\rho$ satisfies $\mh(\w,x)+\lb>\delta>0$ for all $\w\in\W$ and
$x\le\rho$ (whose existence follows from \hyperlink{d2*}{\bf d2$^*$}). It is easy to deduce
that $\lim_{t\to(\alpha_{s,x,\lb})^+}x_\lb(t,s,x)=-\infty$.
This proves the first assertion. The second one can be checked with similar arguments.
The assertions
concerning $\lb\in(\bar\lb_-,\bar\lb^+)$ follow from \cite[Proposition 5.5(v)-(vi)]{dno4}
(see the end of the proof of Theorem \ref{th:3Casint}).
\end{proof}
\begin{teor}\label{th:4Dbifurproceso}
Let $h\colon\RR\to\R$ and $g\colon\R\to\R$ satisfy: the functions $h$ and $h_x$ are admissible;
$\lim_{x\to\pm\infty}(\pm h(t,x))=-\infty$ uniformly on $\R$;
for any $x_1<x_2$, there exists $\delta_{x_1,x_2}>0$ such that
$\liminf_{t\to\pm\infty}(h_{xx}(t,x_1)-h_{xx}(t,x_2))\ge\delta_{x_1,x_2}$; and $g$ is uniformly
continuous, bounded, and with $\inf_{t\in\R} g(t)>0$. Then, for any value of $\lb$,
the equation $x'=h(t,x)+\lb\,g(t)$ has at most three uniformly separated solutions,
in which case they are hyperbolic. Assume also that there
exists $\lb_0\in\R$ such that $x'=h(t,x)+\lb_0\,g(t)$ has three uniformly separated solutions.
Then, all the conclusions of Theorems {\rm\ref{th:4Dbifuruna}} and {\rm\ref{th:4Dasint}} hold for
$x'=h(t,x)+\lb\,g(t)$.
\end{teor}
\begin{proof}
As in the proof of Theorem \ref{th:3Cbifurproceso}, we check that the hypotheses on
$h$ and $g$ ensure that the corresponding hull extensions of $h(t,x)+\lb\,g(t)$
satisfy conditions \ref{d1}, \hyperlink{d2*}{\bf d2$^*$},
\ref{d3} and \ref{d4}. The
unique significative difference is working with $h_{xx}$ instead of
$h_x$. So, \cite[Theorem 5.6]{dno4} shows that the maximum number of
uniformly separated solutions is three, and that, if they exist, they determine three
hyperbolic copies of the base for the skewproduct defined over the hull $\W_{h,g}$.
Hence, the last assertion follows from the fact that all the hypotheses
of Theorem \ref{th:4Dbifur} are satisfied.
\end{proof}
As in the concave case, we point out that
$h_x$ is not required to be concave in Theorem \ref{th:4Dbifurproceso}: for instance,
$h(t,x)=-x^3+x^4\,e^{-t^2}$ satisfies all the assumed conditions on $h$ in its statement.
This fact broadens the possibilities for the application of the result with respect to
the approach of \cite{dno3}.
\par
We finally point out that, if all the hypotheses of Theorem \ref{th:4Dbifurproceso}
hold and if, in addition, the closure $\W_{h,g}$ of the set $\{(h_s,g_s)\,|\;s\in\R\}$ coincides with
the \upalfa-limit set of $(h,g)$, then the nonautonomous bifurcation diagram corresponding to the process
$x'=h(t,x)+\lb\,g(t)$ is that of a classic double local saddle-node bifurcation:
\begin{itemize}[leftmargin=20pt]
\item[-] for $\lb$ in a bounded open interval $\mI=(\bar\lb_-,\bar\lb^+)$, the set of bounded solutions is bounded by two
attractive hyperbolic ones, whose disjoint domains of stability are separated by a repulsive hyperbolic solution;
\item[-] as $\lb$ approaches $\bar\lb_-$ from the right, the two upper hyperbolic solutions approach each other, and loose their
hyperbolicity and uniform separation at $\bar\lb_-$;
\item[-] for $\lb<\bar\lb_-$, there is a unique bounded solution, which is the continuous hyperbolic continuation of
the upper one for $\lb\in\mI$, and its uniform limit as $\lb$ decreases is $-\infty$;
\item[-] and the situation is symmetric as $\lb$ increases, being $\bar\lb^+$ the second bifurcation
value of the parameter, also of (local) saddle-node type.
\end{itemize}
But the situation can be more complex without the condition on the \upalfa-limit set, with the
occurrence of more than one bounded solution (but not of two uniformly separated ones) for $\lb$ in
a bounded interval $\mJ\supsetneq\mI$. This is what we have called a generalized double saddle-node bifurcation
for the process in the Introduction.
\subsection{An application to an electrical circuit model}\label{subsec:41}
\begin{figure}
         \centering
         \includegraphics[scale=0.75]{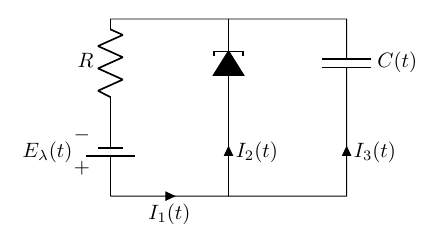}
         \caption{Scheme of the electrical circuit used and modelled in this section.
         It includes a variable capacitor, a resistor, a source and a tunnel diode.}
         \label{fig:realistic}
\end{figure}

In this section, we derive the first order nonlinear differential equation that describes the
time evolution of the charge $Q(t)$ of the capacitor in the scheme of Figure~\ref{fig:realistic}, getting a d-concave equation and
a related bifurcation problem; and we numerically analyze the variation of the bifurcation points provided by the
preceding theory of the d-concave case.

Our electrical circuit is composed of four elements: a tunnel diode, a capacitor, a voltage source
and a resistor. The internal resistance of the source is incorporated in the resistor, to which it is connected
in series: see Figure \ref{fig:realistic}. According to \cite{heinrich} and \cite{nagumo},
the current-voltage characteristic of the tunnel diode
can be suitably approximated by a cubic polynomial: we assume that the current moving through the tunnel
diode at time $t\in\R$ is $I_2(t)=f(V_d(t))$, where $V_d(t)$ stands for the voltage drop across the tunnel diode
at time $t\in\R$ and
\begin{equation}\label{eq:tunneldiodeIVcharacteristic}
 f(V)=I_0-\alpha\,(V-V_0)+\beta\,(V-V_0)^3\,.
\end{equation}
The constants $I_0$, $V_0$, $\alpha>0$ and $\beta>0$ are fixed by the type and materials of the tunnel diode,
which is assumed to operate at room temperature (see \cite{sze}), and they satisfy $0=I_0+\alpha\,V_0-\beta\,V_0^3$, since $f(0)=0$.
Relying on \cite{dellantonio}, \cite{gupta} and \cite{wedzicha}, we assume that the capacitance $C(t)$ of the
capacitor is time-dependent either because one of the capacitor plates is an oscillating membrane, because
it is attached to a vibrating solid, or because it degrades over time.
The function $C(t)$ is assumed to be uniformly continuous and positively bounded from below.
The source produces a time-dependent voltage difference $E_\lb(t)$, which is also uniformly continuous and
can be modified according to the bifurcation parameter $\lb$, and the resistor has resistance $R>0$.
We denote by $I_1(t)$, $I_2(t)$ and $I_3(t)$ the currents  at time $t\in\R$ through the resistor, the tunnel
diode and the capacitor respectively, with the directions shown in Figure \ref{fig:realistic}.

The application of Kirchhoff's first law ensures that $I_1(t)=I_2(t)+I_3(t)$.
The application of Kirchhoff's second law to the right-hand loop including the capacitor and the tunnel diode
ensures that $Q(t)/C(t)=V_d(t)$, and its application to the left-hand loop including the resistor, the voltage
source and the tunnel diode ensures $E_\lambda(t)=V_d(t)+I_1(t)R$. We also know that $Q'(t)=I_3(t)$.
So, joining the previous equations, we obtain
\[
 Q'(t)=I_3(t)=I_1(t)-I_2(t)=\frac{E_\lambda(t)-V_d(t)}{R}-f(V_d(t))\,.
\]
From now on we assume that the form of the parametric variation of the source voltage
is $E_\lambda(t)=E_0(t)+\lambda R$, which can be easily generated with a power supply.
Assuming that $Q$ is the unknown of our differential equation, we rewrite the above
equation using $V_d(t)=Q(t)/C(t)$ and $-I_0=\alpha V_0-\beta V_0^3$,
\begin{equation}\label{eq:circuito}
 Q'=\lambda+\frac{E_0(t)}{R}+\left(\alpha-\frac{1}{R}-3\,\beta\,V_0^2\right)\frac{Q}{C(t)}+
 3\,\beta\,V_0\,\frac{Q^2}{C(t)^2}-\beta\,\frac{Q^3}{C(t)^3}\,.
\end{equation}
To reflect the effect of possible mechanical oscillations and degradation
on the capacitance $C(t)$, we take $C(t)=C_0(t)+A\,p(t)$, where $p\in\mP$ (see \eqref{def:3Deva})
with $k_1=k_2=1$, where $A>0$ is constant, and where $C_0$ may either take a constant value
if there is no time degradation of the capacitor or be a monotonic decreasing function if there is.

\subsubsection{Numerical simulations in an almost periodic example}\label{subsubsec:4.1.1}
Let us define the set $\Theta:=\mathrm{closure}_{C(\R,\R^2)}\{(E_0{\cdot}s,C_0{\cdot}s)\mid s\in\R\}$
(with $(f{\cdot}s)(t):=f(t+s)$) (the joint hull of $E_0$ and $C_0$), let $\W:=\Theta\times\mP$,
let $\theta_0:=(E_0,C_0)\in\Theta$, and
let $\mf_i\colon\Theta\to\R$ and $\mpp\colon\mP\to\R$ be defined by $\mf_i(x_1,x_2):=x_i(0)$ for $i=1,2$
and $\mpp(p):=p(0)$. Then, taking $\theta=\theta_0$ in the varying-in-$\W$ family of equations
\begin{equation}\label{eq:circuitohull}
 Q'=\lambda+\frac{\mf_1(\theta{\cdot}t)}{R}+\rho\,\frac{Q}{\mathfrak{C}(\theta{\cdot}t,p{\cdot}t)}+
 3\,\beta\,V_0\,\frac{Q^2}{\mathfrak{C}(\theta{\cdot}t,p{\cdot}t)^2}-\beta\,\frac{Q^3}{\mathfrak{C}(\theta{\cdot}t,p{\cdot}t)^3}
\end{equation}
we recover the family \eqref{eq:circuito}. Here,
$\rho:=\alpha-1/R-3\,\beta\,V_0^2$ and $\mathfrak{C}(\theta,p):=\mf_2(\theta)+A\,\mpp(p)$.
Since it is a cubic polynomial with a strictly negative leading coefficient,
it is easy to check that \eqref{eq:circuitohull} satisfies the regularity, coercivity and d-concavity
conditions required in Theorem \ref{th:4Dbifur}.
To get an example of applicability of this result and visualize some characteristics of the behavior of
the bifurcation function, we choose $E_0(t):=\sin(t)+\sin(\sqrt{2}\,t)$, $C_0(t):=1$, $R:=1$,
$\alpha:=19/3$, $\beta:=1$, $V_0=:1/3$, and $A:=0.2$, and check that the resulting
equation \eqref{eq:circuitohull}$_{\lb_0}$ has three hyperbolic solutions for all
$\w\in\W$ if $\lb_0:=-0.5$. Let us rewrite it as
$Q'=\mh(\wt,Q)$. It is easy to check that, for all $\w\in\W$, $\mh(\w,Q)\ge \mh_l(Q)$ if $Q\ge 0$
and $\mh(\w,Q)\le\mh_u(Q)$ if $Q\le 0$,
with $\mh_l(Q):=\lb_0+\inf E_0/R+(\rho/(\sup C_0+A))\,Q+(3\,\beta\,V_0/(\sup C_0+A)^2)\,Q^2-(\beta/(\inf C_0-A)^3)\,Q^3$
and $\mh_u(Q):=\lb_0+\sup E_0/R+(\rho/(\sup C_0+A))\,Q+(3\,\beta\,V_0/(\inf C_0-A)^2)\,Q^2-
(\beta/(\inf C_0-A)^3)\,Q^3$. We obtain the local maximum $Q^+_l$ of $\mh_l$ and the local minimum
$Q^-_u$ of $\mh_u$, and check that $Q_u^-<0<Q_l^+$ and $\mh_u(Q_u^-)<0<\mh_l(Q_l^+)$. So,
$\mh(\w,Q_u^-)<0<\mh(\w,Q_l^+)$ for all $\w\in\W$. In addition, we can find $Q_0>\max(Q_l^+,-Q_u^-)$ large enough to ensure
that $\mh(\w,-Q_0)>0>\mh(\w,Q_0)$ for all $\w\in\W$. Standard comparison results allows us to deduce the
existence of three uniformly separated bounded solutions of \eqref{eq:circuitohull}$_{\lb_0}$ for all $\w\in\W$,
and hence \cite[Theorem 5.6]{dno4} ensures that it has three hyperbolic solutions, as asserted. The same result
ensures that the corresponding skewproduct flow has three different hyperbolic copies of the base.

Therefore, all the hypotheses of Theorems \ref{th:4Dbifur} and \ref{th:4Dbifuruna} are satisfied,
and thus there exist two bifurcation points $\tilde\lb_-(\theta,p)<-0.5<\tilde\lb^+(\theta,p)$
such that \eqref{eq:circuitohull} has three hyperbolic solutions for $\lambda\in(\tilde\lb_-(\theta,p),\tilde\lb^+(\theta,p))$,
only one hyperbolic solution for $\lambda=\tilde\lb_-(\theta,p)$ or $\lambda=\tilde\lb^+(\theta,p)$, and at
least one hyperbolic solution which is not uniformly separated from any other bounded solution for
$\lambda\not\in[\tilde\lb_-(\theta,p),\tilde\lb^+(\theta,p)]$.

\begin{figure}[h]
\centering
\includegraphics[width=\textwidth]{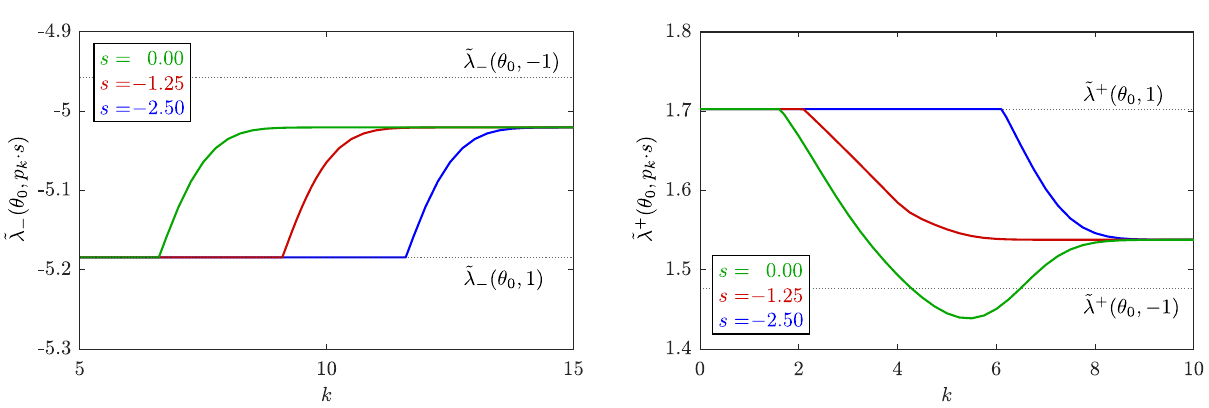}
\caption{
In the left panel (resp. right panel), numerical approximation of the bifurcation curve $\tilde\lb_-(\theta_0,p_k)$ (resp. $\tilde\lb^+(\theta_0,p_k)$) for \eqref{eq:circuitohull}, with
$E_0(t):=\sin(t)+\sin(\sqrt{2}\, t)$, $C_0(t):=1$, $R:=1$, $\alpha:=19/3$, $\beta:=1$, $V_0:=1/3$, $A:=0.2$, and $(p_k{\cdot}s)(t):=\max\{-1,\min\{1,1-k/2-t-s\},\min\{1,1-k/2+t+s\}\}$ (see the left panel of Figure \ref{fig:1}),
for different choices of $s$ in different colors.}
\label{fig:5}
\end{figure}

In order to make some representations, we first take one of the families of functions $p_k$ used
in Section \ref{subsec:32}: $p_k(t):=\max\{-1,\min\{1,1-k/2-t\},\min\{1,1-k/2+t\}\}$
(see the left panel of Figure \ref{fig:1}). In Figure \ref{fig:5}, we compare the variation of
$k\mapsto\tilde\lb_-(\theta_0,p_k{\cdot}s),\tilde\lb^+(\theta_0,p_k{\cdot}s)$ for different time shifts
$s\in\mathbb{R}$. The left and right panels of Figure \ref{fig:5} respectively
show the lower and upper bifurcation curves $k\mapsto\tilde\lb_-(\theta_0,p_k{\cdot}s)$
and $k\mapsto\tilde\lb^+(\theta_0,p_k{\cdot}s)$ for $s=0,-1.25,-2.5$.
Is is interesting to observe the lack of monotonicity of $k\mapsto\tilde\lb^+(\theta_0,p_k{\cdot}0)$:
although $k\mapsto p_k(t)$ and hence $k\mapsto C_k(t):=C_0(t)+A\,p_k(t)$ are nonincreasing, this is not
the case for the right-hand side of \eqref{eq:circuito} with respect to $C(t)$.
This lack of monotonicity is also evident in the fact that $k\mapsto\tilde\lb^+(\theta_0,p_k)$
intersects the line $\tilde\lb^+(\theta_0,-1)$: in this case there is no analogue to
Proposition \ref{prop:3lambda}(i)-(ii). Another point to highlight is that the variations of the
bifurcation points $\tilde\lb_-(\theta_0,p_k)$ and $\tilde\lb^+(\theta_0,p_k)$ with respect to $s$ allow
us to consider bifurcations as a function of the phase $s$.
To explain this assertion, we first point out that, if $k$ is fixed, then $s\mapsto p_k{\cdot}s$
is continuous in the uniform topology of
$C(\R,[-1,1])$, and hence analogous arguments to those of Proposition~\ref{prop:3lambda}(iv) show the continuity of
$s\mapsto\tilde\lb^+(\theta_0,p_k{\cdot}s)$ for a fixed value of $\lb$. Now, for instance, we fix $k=5$ and choose
$\lb_0\in(\,\tilde\lb^+(\theta_0,p_5{\cdot}0),\,\tilde\lb^+(\theta_0,p_5{\cdot}(-1.25))\,)$ (see Figure \ref{fig:5}).
Then, the equation \eqref{eq:circuitohull}$_{\lb_0}$ corresponding to $p_5{\cdot}(-1.25)$ has three hyperbolic solutions,
while that given by $p_5$ has only one. Combined with the continuity, this shows the existence of at least a bifurcation
value of $s\in(-1.25,0)$: changing the phase may also induce a global change in the dynamics.

\begin{figure}[h]
\centering
\includegraphics[width=\textwidth]{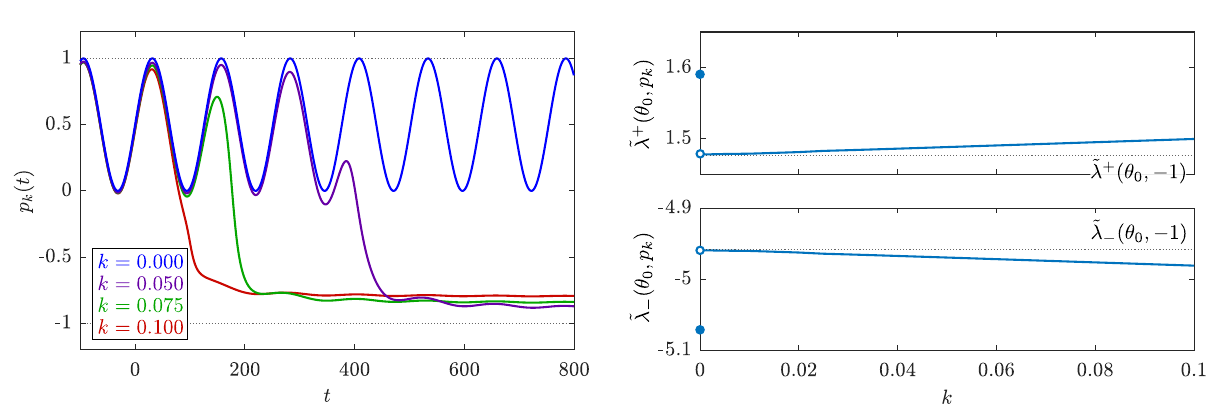}
\caption{In the left panel, $p_0(t):=(\sin(t/20)+1)/2$ and $p_k(t):=(p_0(t)-2k+1)(1/2-(1/\pi)\,\arctan(kt-1/k))+2k-1$
for different choices of $k$ in different colors.
They belong to
the set $\mP$ given by \eqref{def:3Deva} with $k_1=1$ and $k_2=1$.
In the bottom right panel (resp.~top right panel),
a numerical approximation of the bifurcation points $\tilde\lb_-(\theta_0,p_k)$
(resp. $\tilde\lb^+(\theta_0,p_k)$) for \eqref{eq:circuitohull}, with
$E_0$, $C_0$, $R$, $\alpha$, $\beta$, $V_0$, and $A$ as in Figure \ref{fig:5}.
A discontinuity can be observed at $k=0$ in both cases.}
\label{fig:6}
\end{figure}

To conclude, we show numerically the presence of discontinuities of the maps $p\mapsto\tilde\lb^+_-(\theta_0,p)$, as
in Section \ref{subsec:32}. To this end, we take $p_k(t):=((\sin(t/20)+1)/2-2k+1)(1/2-\arctan(kt-1/k)/\pi)+2k-1$ for each
$k\in(0,1]$ and $p_0(t)=(\sin(t/20)+1)/2$. These maps belong to the set $\mP$ given by \eqref{def:3Deva} with $k_1=k_2=1$,
and they satisfy $\lim_{t\to-\infty}(p_k(t)-p_0(t))=0$, $\lim_{t\to\infty}p_k(t)=2k-1$ and $\lim_{k\to0^+}p_k=p_0$ in
$\mP$ (but not uniformly on $\R$). The left panel of Figure \ref{fig:6} depicts some of these functions $p_k$, and
the right panels depict numerical approximations to the maps $k\mapsto\tilde\lb_-(\theta_0,p_k)$
and $k\mapsto\tilde\lb^+(\theta_0,p_k)$ that show their discontinuity as $k\downarrow 0$.
The continuity of these applications in $(0,1]$ follows from the continuity in the uniform topology
of $k\mapsto p_k$ in this interval, reasoning as in the proof of Proposition \ref{prop:3lambda}(iv).
\subsection{Critical transitions in a d-concave in measure population model}\label{subsec:42}
We finish with a numerical simulation of a critical extinction in a ecological
community described by the d-concave equation \eqref{eq:huntingintrod} of the Introduction,
d-concave in measure but not asymptotically d-concave. We will briefly study numerically
the bifurcation with respect to the magnitude of predation during its occurrence intervals:
\begin{equation}\label{eq:lastsimulation}
x'=\frac{1}{6}\,x\,(1-x)(x-2)-\lambda\,b(t)\,\frac{x^2}{1+x^2}\,,
\end{equation}
which we write for short as $x'=h(t,x,\lb)$.
Recall that the function $b$ represents predation starting at a certain point in time
and occurring at increasingly spaced intervals (see the left panel of Figure \ref{fig:tophull}).
So, a higher $\lb$ means a stronger predation.
As explained in the Introduction, and in the line of approach posed in \cite{dno4},
the existence of the compact hull $\W_b$ of the bounded and
uniformly continuous map $b$ allows us to identify \eqref{eq:lastsimulation}$_\lb$
(or any of its time-shifts, obtain by changing $b$ by $b_s$) as a transition
equation between the two families of equations defined over the \upalfa-limit set $\W_b^\alpha$ and
the \upomeg-limit set $\W_b^\omega=\W_\G$ of $b$ on its hull (see the right panel of
Figure \ref{fig:tophull}). However, the results of \cite{dno4} are not
applicable here: the map $h_-(t,x):=h(t,x,0)$, strictly d-concave and giving rise to
three hyperbolic solutions, provides the ``past equation" required in
that work; but there is no way to define a ``future equation". The results of this paper
make this condition unnecessary: the existence of the
family of equations over $\W_b^\omega=\W_\G$ suffices for the analysis.

\begin{figure}[h]
\centering
\includegraphics[width=\textwidth]{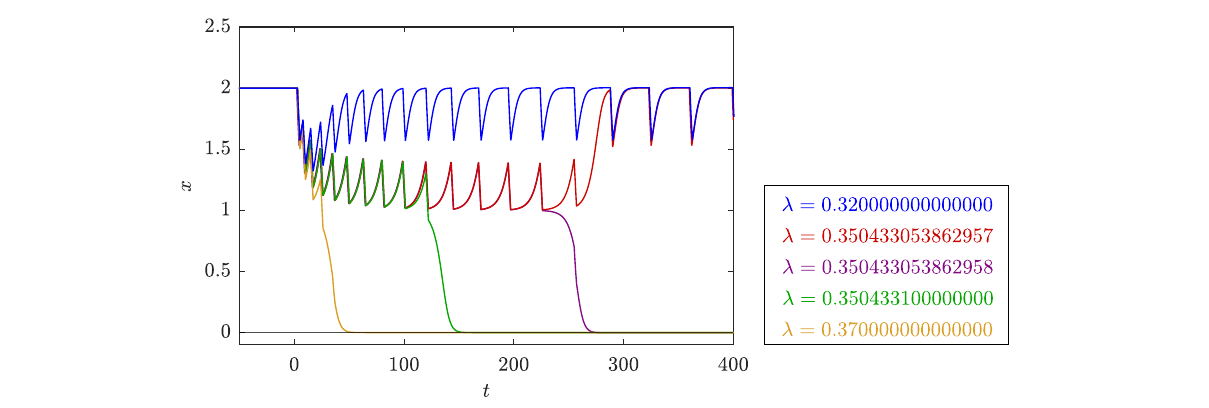}
\caption{
Numerical approximation of the attractive hyperbolic solutions of equation
\eqref{eq:lastsimulation}$_\lb$ for different values of the parameter $\lb$.
In black, the constant solution $0$, which is an attractive hyperbolic solution
for all the parameter values. In the colors indicated in the legend, the attractive
hyperbolic solution corresponding to the upper bounded solution of the equation.}
\label{fig:lastfigure}
\end{figure}
Note that the parametric perturbation in \eqref{eq:lastsimulation} does not match with that
of Theorem \ref{th:4Dbifurproceso}. (In fact, also the corresponding bifurcation diagram is different.)
Let us list some facts easy to establish. For $\lb=0$, the equation has three hyperbolic solutions:
the three equilibrium points $0$, $1$ and $2$ of $x'=(1/6)\,x\,(1-x)(2-x)$. In addition, $0$ is an
attractive hyperbolic solution of \eqref{eq:lastsimulation} for all $\lb$, since $h_x(t,0,\lb)=-1/3<0$
for all $\lb$. It is easy to check that $0$ is the lower bounded solution $l_\lb$ for all $\lb\in\R$,
and that the upper bounded solution $u_\lb$ takes the value 2 for all $t\le 0$. This precludes the
existence of a unique bounded solution for any $\lb$ (which precludes the bifurcation diagram
of Theorem \ref{th:4Dbifurproceso}). In addition, Theorem \ref{th:2persistencia} shows
that \eqref{eq:lastsimulation} has three hyperbolic solutions is $\lambda>0$ is small enough, with two
positive hyperbolic solutions, $u_\lb$ (attractive) and $m_\lb$ (repulsive), uniformly separated
from~$0$.

Now, a numerical analysis shows that this situation persists until $\lb$ reaches a
bifurcation value of $\lb$, corresponding to a critical transition in the model:
$\lb_0\approx0.350433053862957$, partly depicted in Figure \ref{fig:lastfigure}
(where the repulsive solution is not drawn).
It can be numerically proved that, for all the depicted values of $\lb$,
the equation $x'=(1/6)\,x\,(1-x)(x-2)-\lambda\,\G(t)\,x^2(1+x^2)$ (which is one of the equations of the
\upomeg-limit family) has three hyperbolic solutions, and that, for $\lb<\lb_0$, these three hyperbolic
solutions are approached uniformly on compact sets
through suitable increasing sequences $(t_n)\uparrow\infty$ by the three
hyperbolic solutions of the transition equation:
in the words of \cite{dno4}, there is tracking of local attractors,
and also of the global attractor. This means the survival of an initially not too small (over $m_\lb(0)$)
population. However, for $\lb>\lb_0$, all the bounded solutions of \eqref{eq:lastsimulation} (and hence all the solutions)
converge to the attractive hyperbolic solution 0 as time increases: using again the words of \cite{dno4},
the connection of global attractors is lost after the critical transition at $\lb_0$, and
tipping occurs. For these values of $\lb$, any initial population gets extinct. A standard
comparison argument shows that this extinction situation persists for any $\lb>\lb_0$: the monotonicity
with respect to $\lb$ provides a unique critical transition in this case.

\end{document}